\documentclass[12pt]{amsart}
\usepackage{amssymb, amsthm}
\usepackage[utf8]{inputenc}
\usepackage{comment}
\usepackage{xcolor}
\usepackage[all]{xy}
\usepackage{tikz}
\usepackage{tkz-euclide}
\usepackage{tikz-cd}
\usetikzlibrary{calc, positioning, decorations.pathreplacing}
\usepackage{pgfplots}
\pgfplotsset{compat=1.18}
\usepackage{fullpage}
\usepackage{mathdots} 

\usepackage{listings} 
\usepackage{tcolorbox}
\tcbuselibrary{listings} 
\usepackage{upquote} 
\definecolor{codegreen}{rgb}{0,0.6,0}
\definecolor{darkblue}{rgb}{0.1, 0.2, 0.6}
\newtcblisting{inputbox}[1]{
    colback=gray!5,       
    colframe=gray!40,     
    colbacktitle=gray!15, 
    coltitle=black,        
    title=\textbf{In [#1]:}, 
    boxsep=1pt, 
    left=2pt, right=2pt, top=0pt, bottom=0pt,
    toptitle=1pt, bottomtitle=1pt,
    before skip=6pt, after skip=2pt, 
    listing only,
    listing options = {
        language=Python, 
        basicstyle=\ttfamily\footnotesize, 
        showstringspaces=false, 
        breaklines=true, 
        breakatwhitespace=true, 
        postbreak=\mbox{\normalfont\textcolor{gray}{$\hookrightarrow$}\space}, commentstyle=\color{codegreen}, keywordstyle=\color{darkblue}\bfseries}}
\newtcolorbox{outputbox}[1]{
    colback=white,
    colframe=gray!30,
    colbacktitle=white,      
    coltitle=black,          
    titlerule=0.4pt,         
    title=\textbf{Out [#1]:},
    boxsep=1pt,
    left=2pt, right=2pt, top=2pt, bottom=2pt,
    toptitle=1pt, bottomtitle=1pt,
    before skip=2pt, after skip=6pt}

\newtheorem{thm}{Theorem}[section]
\newtheorem{cor}[thm]{Corollary}
\newtheorem{lem}[thm]{Lemma}
\newtheorem{prop}[thm]{Proposition}
\theoremstyle{definition}
\newtheorem{definition}[thm]{Definition}
\newtheorem{remark}[thm]{Remark}
\newtheorem*{remark*}{Remark}
\newtheorem{ex}[thm]{Example}

\newtheorem{thmAlph}{Theorem}
\newtheorem{propAlph}[thmAlph]{Proposition}

\newcommand{\adm}{\operatorname{Adm}}
\newcommand{\dadm}{\operatorname{Adm}^{\mathrm{OP}\bullet}}
\newcommand{\bP}{\mathbb{P}}
\newcommand{\renzo}[1]{{\color{blue}{Renzo: #1}}}

\newcommand{\ain}{\mathbf{a}}
\newcommand{\e}{\varepsilon}
\newcommand{\h}{\hbar}
\newcommand{\Tr}[1]{\mathcal{T}_{#1}}
\newcommand{\Trz}[1]{\mathsf{T}_{#1}}

\newcommand{\sumn}[2]{\sum_{#1\in\mathbb{N}^{#2}}}

\newcommand{\pb}{p_{\vec{b}}}

\newcommand{\mgbar}[1]{{\overline{\mathcal{M}}}_{#1}}

\title{$\lambda_g$-Integrals on Triple, One-part Admissible Covers}
\author{Renzo Cavalieri, David Klompenhouwer and Paolo Rossi}
\date{}

\begin{document}

\begin{abstract}
    We effectively compute all intersection numbers of the form $\lambda_g\psi^d$ on moduli spaces of genus $g$ admissible covers with one point of full ramification and two additional branch points with arbitrary profile. We encode these invariants in appropriate generating functions, which we relate to the Hamiltonians of the KdV hierarchy.
\end{abstract}

\maketitle

\tableofcontents

\section{Introduction}

We study a class of Hodge integrals on moduli spaces of admissible covers, for which we provide both an effective reconstruction scheme as well as  structure results. We begin this introduction by stating the results and then give a motivating discussion for the problem and the techniques used.

\subsection{Summary of results}

For a positive integer $\delta$, let $\nu_1, \ldots, \nu_s$ denote  integer partitions of $\delta$. We denote by $\adm_{g}(\nu_1, \ldots, \nu_s)$ the moduli space of genus $g$ admissible covers of a rational curve of degree $\delta$, with marked fibers over $s$ special points and ramification specified by the partitions, and simple branching at $r$ other points (see Definition \ref{def:ADM}). There is a natural branch morphism $br:\adm_{g}(\nu_1, \ldots, \nu_s)\to \overline{M}_{0, s+r}$. \newline

Our first main result is a reconstruction scheme for certain Hodge integrals on spaces of admissible covers with three not necessarily simple branch points, where the cover is fully ramified over one of the branch points.
\begin{thmAlph}[\bf{Reconstruction}]\label{thm:reconstruct}
    For all nonnegative integers $g$, positive integers $\delta$, and partitions $\nu_2, \nu_3$ of $\delta$, let $d = g-\delta-1+\ell(\nu_2)+\ell(\nu_3)\ge 0$. Any invariant of the form
    \begin{equation}\label{eq:inv}
        \int_{\adm_{g}( \delta,\nu_2,\nu_3)} \lambda_g\, br^\ast(\psi_2^d)
    \end{equation}
    may be effectively reconstructed from the initial conditions
    \begin{equation}\label{eq:incond}
    \int_{\adm_{0}(\delta,\delta)} 1 := \frac{1}{\delta}, \ \ \ \ \ \
        \int_{\adm_{0}(\delta,(i, \delta-i))} 1 = 1, \ \ \ \ \ \  \int_{\adm_{1}( \delta,\delta)} \lambda_1 = \frac{\delta^2-1}{12},
    \end{equation}
    via WDVV recursions and the string equation (Lemma \ref{lem:string}).
\end{thmAlph}
The recursive relations used to reconstruct the invariants also involve Hodge integrals on moduli spaces of disconnected admissible covers, where over one of the branch points  is a point of full ramification in each connected component of the cover (see Definition \ref{def:ADMOP}). These invariants are not necessarily related to the product of invariants relative to the covers from each connected component, but are also reconstructed from the initial conditions \eqref{eq:incond}. \newline

The building blocks of the recursive structure are a smaller class of connected invariants of type \eqref{eq:inv}, where the power of the $\psi$-class is zero and the preimage of the second branch point has exactly one point of non-trivial ramification $k+1$, i.e.
\begin{equation*}
    \nu_2 = (k+1 ,1^{\delta-k-1}), \quad k\ge1.
\end{equation*}
We write these invariants as
\begin{equation}\label{eq:bblock}
    x_{k;g}(\ain) := \frac{1}{({\delta-k-1})!}\int_{\adm_g\left(\delta,\,(k+1 ,1^{\delta-k-1}),\,\ain\right)}\lambda_g, \quad 0\le g\le k, \quad \ain\in\mathbb{N}^{k+1-g}.
\end{equation}
In the main body of the paper we define these invariants (and all the subsequent ones that we introduce in this section) without the factorial sign in front, since we will use a version of admissible cover space where only the nontrivial points associated to $\nu_2$ are marked (see Definition \ref{def:ADM}). We encode these in a generating series $\Trz{k}$ (Definition \ref{def:Trz}). The first series $\Trz{1}$ encodes the second and third initial conditions \eqref{eq:incond}. The first step towards the reconstruction theorem is showing that the $x_{k;g}(\ain)$ (or, equivalently, their generating series $\Trz{k}$) can be effectively reconstructed from the initial conditions, via WDVV relations pulled back via the branch morphism. The quantum $\star$-product \cite{Moyal_1949} provides algebraic language in which to state this result. The $\star$-product is defined, with examples, in Section \ref{subsec:star}.
\begin{propAlph}\label{prop:Tk}
    For any $k>1$, the generating function $\Trz{k}$ is reconstructed from $\mathsf{T}_1$ and the unique evaluation $x_{k;0}(1^{k+1}) = k!$ via the commutation relation involving the $\star$-product
    \begin{equation}\label{eq:commT}
        [\Trz{k}, \Trz{1}] := \Trz{k}\star\Trz{1} - \Trz{1}\star\Trz{k} = 0.
    \end{equation}
\end{propAlph}
The commutator encodes WDVV relations pulled back via the branch morphism involving invariants encoded by $\Trz{k}$. This first reconstruction result can be readily implemented to compute the invariants \eqref{eq:bblock}. We demonstrate this in Appendix \ref{app:python}. For example, for $k=2$,
\begin{equation*}
    x_{2;0}(a_1,a_2,a_3) = 2, \quad x_{2;1}(a_1,a_2) = \frac{a_1^2+a_2^2}{6} + \frac{a_1a_2}{6} - \frac12, \quad x_{2;2}(a)=\frac{(a^2-1)(a^2-4)}{60}.
\end{equation*}
In Table \ref{tab:xkg} in Section \ref{sec:recTk} we write down the invariants for $k=3$ and $k=4$. By degenerating the target rational curve \cite{op_trilogy1} and applying a wall-crossing formula from \cite{nss2025}, we show that these invariants are related to Hodge integrals on the double ramification cycle.
\begin{propAlph}[Proposition \ref{prop:Adm2DR}] \label{prop:C}
    For every ramification degree $k\ge1$, genus $0\le g\le k$, and profile $\ain=(a_1,\dots,a_{k+1-g})\in\mathbb{N}^{k+1-g}$,
    \begin{equation}\label{eq:admdrinv}
        \frac{x_{k;g}(\ain)}{k!g!} = \sum_{h=0}^g\frac{(k+1)!}{2^h(k+1-h)!}c_{k,h} \int_{\operatorname{DR}_{g-h}(-\sum_ia_i,0,a_1,\dots,a_{k+1-g})}\lambda_{g-h}\psi_2^{k-h},
    \end{equation}
    where $c_{k,h}$ is the coefficient of $t^{2h}$ in $\left(\frac{\operatorname{arcsinh}(t/2)}{t/2}\right)^{2(k+1-h)}$.
\end{propAlph}
The $g=k$ case of the result above yields (Lemma \ref{cor:xgg})
\begin{equation*}
    x_{g;g}(a) = \frac{(g!)^2}{2^{2g}(2g+1)!!}\prod_{n=1}^g(a^2-n^2),
\end{equation*}
which should be compared with the corresponding DR cycle intersection number
\begin{equation*}
    \int_{\operatorname{DR}_g(-a,0,a)}\lambda_g\psi_2^g = \frac{a^{2g}}{2^{2g}(2g+1)!!}.
\end{equation*}
The DR integrals in the proposition appear in the DR hierarchy associated to the trivial Cohomological Field Theory, which is equivalent to the KdV hierarchy of partial differential equations \cite{bur15}. Therefore Proposition \ref{prop:C} above can be expressed as a relation between the generating series $\Trz{k}$ and the Hamiltonian flows of the KdV hierarchy (Corollary \ref{cor:Tk_KdV}). The more general commutation relations 
\begin{equation*}
    [\Trz{k_1},\Trz{k_2}] = 0, \quad k_1,k_2\ge0.
\end{equation*}
hold. This corresponds to the Poisson commutativity of the Hamiltonians of the KdV hierarchy. \newline

Next, we consider invariants of type \eqref{eq:inv} where we allow the ramification over the second branch point to spread out over more than one ramification point, while still requiring the power of the $\psi$-class to be zero. They are indexed by nonempty integer partitions $\mu=(\mu_1,\dots,\mu_\ell)$, satisfying $\delta\ge\lvert\mu\rvert+\ell$ and
\begin{equation*}
    \nu_2= \left(\mu_1+1,\dots,\mu_\ell+1,1^{\delta-\lvert\mu\rvert-\ell}\right).
\end{equation*}
We write these invariants as
\begin{equation*}
    x_{\mu;g}(\ain) := \frac{1}{(\delta-\lvert\mu\rvert-\ell)!}\int_{\adm_g\left(\delta,\nu_2,\,\ain\right)}\lambda_g, \quad 0\le g\le\lvert\mu\rvert, \quad \ain\in\mathbb{N}^{\lvert\mu\rvert+1-g}.
\end{equation*}
This generalizes \eqref{eq:bblock} in the sense that $x_{(k);g}(\ain)=x_{k;g}(\ain)$. We encode these in generating series $\Tr{\mu}^\bullet$ (see Definition \ref{def:Tmu}) which are defined so that $\Tr{(k)}^\bullet = \Trz{k}$. The $\bullet$ indicates that it is now necessary to include contributions from invariants arising from disconnected curves. In this case, these are just the products of the associated connected invariants. The extended system of invariants $x_{\mu;g}(\ain)$ can be expressed in terms of the building blocks $x_{k;g}(\ain)$. 
\begin{propAlph}[Proposition \ref{prop:TmuTk}]
    For every partition $\mu=(\mu_1,\dots,\mu_\ell)$, the generating series $\Tr{\mu}^\bullet$ is a polynomial in the $\Trz{k}$, with multiplication given by the $\star$-product, with leading term 
    \begin{equation*}
        \Tr{\mu}^\bullet = \frac{1}{\lvert\operatorname{Aut}(\mu)\rvert}\Trz{\mu_1}\star\cdots\star\Trz{\mu_\ell} + \cdots.
    \end{equation*}
\end{propAlph}
The first two non-trivial examples of this are
\begin{equation*}
    \Tr{(1,1)}^\bullet = \frac12\Trz{1}\star\Trz{1} - \frac32\Trz{2} - \frac14\h\e^2v\Trz{1}, \quad \Tr{(2,1)}^\bullet=\Trz{2}\star\Trz{1}-4\Trz{3}-2\h\e^2v\Trz{2}.
\end{equation*}
In Appendix \ref{appsec:TmuTk} we write the explicit relations for all $\lvert\mu\rvert\le4$. \newline

We then study invariants where the power of the $\psi$-class is  positive. 
With $\nu_2$ as above and $d\ge0$,
\begin{equation}\label{eq:geninv}
    x_{\mu,d;g}(\ain) := \frac{1}{(\delta-\lvert\mu\rvert-\ell)!} \int_{\adm_g\left(\delta,\nu_2,\,\ain\right)}\lambda_g\,br^*(\psi_2^d), \quad 0\le g\le\lvert\mu\rvert+d, \quad \ain\in\mathbb{N}^{\lvert\mu\rvert+d+1-g}.
\end{equation}
Clearly, $x_{\mu,0;g}(\ain)=x_{\mu;g}(\ain)$. We encode these in a generating series $\Tr{\mu,d}^\bullet$ (Definition \ref{def:Tmud}), defined so that $\Tr{\mu,0}^\bullet = \Tr{\mu}^\bullet$. In this case, the disconnected invariants are generally not the product of the associated connected invariants. By expressing the $\psi$-class as a sum of boundary divisors, we deduce that the larger system of invariants $x_{\mu,d;g}(\ain)$ can be expressed in terms of the $d=0$ invariants $x_{\nu;g}(\ain)$.
\begin{propAlph}[Proposition \ref{prop:Tds}]
    For every $d\ge1$, the generating series $\Tr{\mu,d}^\bullet$ is written as a linear combination of $\star$-products of the $\Tr{\nu}^\bullet$, with coefficients determined explicitly by $\Tr{\mu,d-1}^\bullet$.
\end{propAlph}
Two examples of this are
\begin{align*}
    \Tr{(1,1),1}^\bullet = & \; v\Trz{1}\star\Tr{(1,1)}^\bullet - \h\e^2\frac{v^2}{2!}\Tr{(1,1)}^\bullet, \\
    \Tr{(2),2}^\bullet = & \; \frac{v}{2}\Trz{1}\star\Tr{(2),1}^\bullet - \h\e^2\frac{v^3}{3!}\left(\frac72\Tr{(2,1)}^\bullet + 24\Tr{(3)}^\bullet\right) - \h^2\e^4\frac{v^4}{4!}12\Tr{(2)}^\bullet.
\end{align*}
In Appendix \ref{appsec:TmudTmuprop} we write all such relations for $\lvert\mu\rvert+d\le5$. An immediate consequence of the two propositions above is that all the invariants $x_{\mu,d;g}(\ain)$, and also their disconnected counterparts (see Definition \ref{def:Smud}), can be expressed in terms of the building blocks $x_{k;g}(\ain)$, as described in the following theorem.
\begin{thmAlph}[Theorem \ref{thm:Tmud_in_Tk}] \label{thm:F}
    For every partition $\mu=(\mu_1,\dots,\mu_\ell)$ and integer $d\ge 0$, the generating series $\Tr{\mu,d}^\bullet$ can be written as a linear combination of $\star$-products of the $\Trz{k}$, with leading term
    \begin{equation*}
        \Tr{\mu,d}^\bullet = \frac{1}{\lvert\operatorname{Aut}(\mu)\rvert}\frac{v^d}{d!}\Trz{1}^{\star d}\star\Trz{\mu_1}\star\cdots\star\Trz{\mu_\ell} + \cdots.
    \end{equation*}
\end{thmAlph}
The two examples above become
\begin{align*}
    \Tr{(1,1),1}^\bullet = & \; v\left(\frac12\Trz{1}^{\star3} - \frac32\Trz{2}\star\Trz{1}\right) - \h\e^2\frac{v^2}{2!}\left(\Trz{1}^{\star2} - \frac32\Trz{2}\right) + \h^2\e^4\frac{v^3}{3!}\frac{3}{4}\Trz{1}, \\
    \Tr{(2),2}^\bullet = & \; \frac{v^2}{2!}\Trz{2}\star\Trz{1}^{\star2} - \h\e^2\frac{v^3}{3!}\left(\frac{13}{2}\Trz{2}\star\Trz{1} + 10\Trz{3}\right) + \h^2\e^4\frac{v^4}{4!}16\Trz{2}.
\end{align*}
In Appendix \ref{appsec:Tmud_in_Tk} we write all such relations for $\lvert\mu\rvert+d\le5$. Theorem \ref{thm:F}, together with Proposition \ref{prop:Tk}, proves the overall reconstruction statement of Theorem \ref{thm:reconstruct}. 
Proposition \ref{prop:C} and Theorem \ref{thm:F}, together with  the polynomiality of the double ramification cycle \cite{pix23, spe24}, allow us to prove
that all the invariants \eqref{eq:geninv} are polynomial. 
\begin{thmAlph}[\bf{Polynomiality}]\label{thm:poly}
    Fix a nonempty partition $\mu$ and integer $d\ge0$. For every $0\le g\le\lvert\mu\rvert+d$, the family of invariants $x_{\mu, d;g}(a_1,\dots,a_{\lvert\mu\rvert+d+1-g})$ is encoded by a symmetric polynomial in the $a_1,\dots,a_{\lvert\mu\rvert+d+1-g}$ of degree $2g+d +\ell(\mu)-1$.
\end{thmAlph}
We also give a more general statement in Theorem \ref{thm:disc_poly} regarding polynomiality of invariants arising from disconnected admissible covers.
\subsection{Context and motivation}

The origin of this work lies at the intersection of two distinct, but closely related, mathematical theories: the double ramification (DR) hierarchy and Hurwitz theory. On the one hand, the double ramification locus in the moduli space $\mathcal M_{g,n}$ of smooth pointed curves consists of curves admitting a map to $\mathbb P^1$ with prescribed ramification over $0$ and $\infty$. Equivalently, for integers $a_1,\ldots,a_n$ summing to zero, it is the locus of $(C,x_1,\ldots,x_n)$ such that $\mathcal O_C(\sum_i a_i x_i)\cong\mathcal O_C$. Motivated by Symplectic Field Theory \cite{Intro2SFT}, the problem of extending this locus to obtain a cycle in $\overline{\mathcal M}_{g,n}$ led to the DR cycle, initially defined using stable relative maps to a rubber target \cite{li02,gv05}; see also \cite[Section~1.6]{hs21} for other approaches. Pixton's formula for this cycle, proved in \cite{jppz17}, implies in particular its polynomiality in the ramification multiplicities \cite{pix23,spe24}. Buryak used the properties of the DR cycle to associate an integrable Hamiltonian system of PDEs to any Cohomological Field Theory \cite{bur15}, as part of the broader program initiated by Witten's conjecture \cite{wit90}. The resulting DR hierarchy and its quantization have been studied extensively \cite{br16,BR2016q,BDGR18,brs21}. For the trivial CohFT, the DR hierarchy is equivalent to the KdV hierarchy and is governed by intersection numbers of the form
\[
  \int_{\mathrm{DR}_g(-\sum_i a_i,0,a_1,\ldots,a_n)}\lambda_g\psi_2^d.
\]

On the other hand, Hurwitz spaces parameterize covers of $\mathbb P^1$ with prescribed discrete invariants and ramification data. They admit a smooth, proper orbifold compactification by admissible covers \cite{hm:kd,acv03,Jarvis_Kaufmann_Kimura}. Hurwitz numbers, the degrees of the corresponding branch morphisms, are natural enumerative invariants arising from these spaces. The case with one special branch point, usually called simple Hurwitz numbers, is  well understood; notably, the ELSV formula relates these counts to Hodge integrals on moduli spaces of curves \cite{elsv:ohnahi,elsv:hnaiomsoc}. Double Hurwitz numbers, counting covers with two special branch points, have also been intensely studied, beginning with \cite{gjv:ttgodhn}. They are piecewise polynomial in the ramification multiplicities, of degree  $4g-3+\ell(\alpha)+\ell(\beta)$ for profiles $\alpha$ and $\beta$ of fixed lengths; their chambers and wall crossings have been investigated in \cite{ssv:g0,cjm:thn,cjm:wcfdhn}. Triple Hurwitz numbers are substantially more complicated, as they encode  the full character theory of symmetric groups;  degeneration expresses Hurwitz numbers for covers with additional special branch points in terms of such three-point data. Hodge integrals, and more generally tautological intersections, on spaces of admissible covers---often called Hurwitz--Hodge integrals---were initially investigated in connection with orbifold Gromov--Witten theory \cite{bgp:crc,cc:c3z3}; they are also of intrinsic interest and arise naturally in the study of the tautological intersection theory of moduli spaces of curves \cite{r:dd2,jpt:elsv}.

The point of confluence of these two theories, and the initial motivation for the collaboration of the authors, is the study of $\lambda_g$-integrals on spaces of covers with two distinguished ramification profiles, which may be compactified either by rubber relative stable maps, as in the construction of the DR cycle, or by admissible covers. There is a familiar trade-off between these two compactifications: admissible covers have a better geometric behavior, while the stable-map compactification has better formal properties, most notably the polynomiality that is crucial for the DR hierarchy; see \cite[Section~2.3]{bssz15}. The natural admissible-cover analogues of the DR intersection numbers above led us initially to consider
\[
  x_{(1),d;g}(\mathbf a)
   = \frac{1}{(\delta-2)!}\int_{\operatorname{Adm}_g(\delta,(2, 1^{\delta-2}),\mathbf a)}
       \lambda_g\,br^*(\psi_2^d),
  \qquad \delta=\sum_i a_i,
\]
where the cover is fully ramified over $0$ and has ramification profile $\mathbf a$ over $\infty$. The insertion of $\lambda_g$, natural from the perspective of the DR hierarchy, also simplifies the combinatorics of the recursions: it vanishes on strata whose source curves are not of compact type. The choice of a one-part ramification profile is motivated by double Hurwitz theory, where it singles out a fundamental chamber of polynomiality  \cite{gjv:ttgodhn}. This seemed the natural starting point for our analysis.

In the process of reconstructing these invariants using WDVV relations, or equivalently the commutativity of the corresponding Hamiltonian flows, we realized that the same structure allows one to prescribe an arbitrary, but fixed, ramification profile over a third special branch point. This considerably enlarges the class of invariants accessible to our methods and, from the perspective of Hurwitz theory, opens the door to the more mysterious geometry of Hodge integrals on triple Hurwitz spaces. The change in viewpoint also revealed that the natural building blocks for the reconstruction are not the descendants $x_{(1),d;g}(\mathbf a)$, but the integrals with no $\psi$-class and only one nontrivial ramification cycle over the third special point:
\[
  x_{k;g}(\mathbf a)
     =\frac{1}{(\delta-k-1)!}\int_{\operatorname{Adm}_g(\delta,(k+1, 1^{\delta-k-1}),\mathbf a)}\lambda_g,
  \qquad \mathbf a\in\mathbb N^{k+1-g},\quad 0\leq g\leq k.
\]
These invariants arise from admissible covers of connected curves, since the cover has a point of full ramification over $0$. More importantly, the WDVV recursions used for their reconstruction do not involve  contributions from spaces of disconnected covers. This greatly simplifies the proof that the recursions determine them uniquely from the initial conditions. For fixed $k$ and degree, we organize the resulting linear system of principal parts into blocks according to genus: the distinguished square blocks have positive diagonal coefficients, while the nonzero coefficients in the neighboring blocks have negative sign. An analysis of these signs shows that the columns of the matrix are linearly independent. This argument becomes substantially more complicated for the general descendant invariants, where disconnected contributions can no longer be ignored.

The correct encoding in generating functions is a crucial algebraic language for carrying out this reconstruction efficiently. We encode the $x_{k;g}(\mathbf a)$ in the series $\Trz{k}$, and the more general invariants, including the necessary contributions from disconnected admissible covers, in the series $\mathcal{T}^\bullet_{\mu,d}$. The quantum $\star$-product \cite{Moyal_1949}, used in the context of Symplectic Field Theory and quantum DR hierarchies \cite{Intro2SFT,BR2016q}, is an efficient algebraic model for gluing components of covers arising when degenerating the target. In this language, geometric relations such as WDVV and the expression of $\psi$-classes in terms of boundary divisors become algebraic relations among generating functions. We obtain, in particular, an effective expression for every $\Tr{\mu,d}^\bullet$ as a polynomial in the $\Trz{k}$'s, with multiplication given by the $\star$-product. The coefficients of these expressions are determined by explicit recursions, which are implemented in \texttt{Python} with the help of the \texttt{pyBoLaNO} package \cite{pybolano} in Appendix \ref{app:python}. 

There is a further feature of this description worth emphasizing. The generating functions satisfy the much more general commutation relations
\[
 [\Tr{\mu,d}^\bullet,\Trz{1}]=0,
 \qquad
 [\Tr{\mu_1,d_1}^\bullet,\Tr{\mu_2,d_2}^\bullet]=0.
\]
The first of these relations suggests a second route to reconstruction: for a fixed, but arbitrary, profile $\mu$, one can use the commutation with $\Trz{1}$ to reconstruct the family $\Tr{\mu,d}^\bullet$ for all $d$, recursively using only the invariants associated with that fixed profile, together with the initial data encoded by $\Trz{1}$. This is precisely the kind of reconstruction we initially sought in the case $\mu=(1)$. In computations, these commutation relations do appear sufficient to determine the invariants. Establishing their sufficiency directly, however, requires a much more intricate analysis of the corresponding linear systems. For this reason, we found it more natural to first reconstruct the $\Trz{k}$'s and then express all the other generating functions as $\star$-polynomials in these building blocks.

We  return to the original relationship with the DR hierarchy. We compare the building blocks $x_{k;g}(\mathbf a)$ with intersection numbers on the DR cycle by applying the tooth-and-comb degeneration of  \cite{op_trilogy1}, followed by the wall-crossing formula  relating admissible covers and stable maps \cite{nss2025}. This yields an explicit relation between the $x_{k;g}(\mathbf a)$ and specific Hodge integrals on DR cycles, and hence a relation between the series $\Trz{k}$ and the Hamiltonians of the KdV hierarchy. Since DR cycles are polynomial in their ramification multiplicities \cite{pix23,spe24}, the relation immediately implies polynomiality of the building blocks. The reconstruction by $\star$-products then propagates this property to all the invariants studied in this paper, giving a polynomiality result for $\lambda_g$-Hodge integrals on  one-part admissible-cover spaces.

\subsection{Future work}

There are several natural directions that we intend to pursue, which we expect will strengthen the ties between Hurwitz theory and the DR hierarchy. First, we would like to allow insertions of tautological classes other than the top Hodge class $\lambda_g$, and ultimately to remove the assumption of full ramification over one of the branch points. In this setting the WDVV recursions can involve nonseparating nodes and a richer system of Hodge classes. Some preliminary computations suggest that the resulting families of invariants might satisfy structural relations naturally connected to the quantum DR hierarchy \cite{BR2016q}.

We also expect that these more general Hurwitz--Hodge integrals exhibit a wall-and-chamber structure analogous to that of double Hurwitz numbers. More precisely, as the ramification multiplicities vary, we expect the walls to lie in the resonance hyperplane arrangement defined by the vanishing of sums of subsets of the signed ramification data. It would be interesting to extend the picture developed for double Hurwitz numbers \cite{ssv:g0,cjm:wcfdhn} to arbitrary Hodge integrals on spaces of admissible covers.

We expect that the $\star$-product will be a natural tool in this broader context. As the example in Section~4.2 shows, once the one-part  condition is relaxed, the $\star$-product of generating functions with polynomial coefficients can produce piecewise-polynomial expressions. A particularly suggestive example is the tropical algorithm for double Hurwitz numbers of \cite{cjm:thn,cjm:wcfdhn}:  its successive gluings can be naturally interpreted as successive $\star$-products of a generating function $\widehat{\Trz{1}^{0}}$ collecting invariants for both left- and right-pointing genus-zero tripods. Understanding this interpretation more systematically could explain the appearance of resonance walls directly in the algebra of generating functions and provide a useful framework for wall-crossing formulas for Hurwitz--Hodge integrals.

\subsection{Organization of the paper}

In Section~2 we introduce the moduli spaces of admissible covers, the tautological classes and generating functions used throughout the paper. We also recall the relevant WDVV relations and introduce the $\star$-product as an algebraic model for gluing.
Section~3 proves the reconstruction results. We first reconstruct the building blocks $\Trz{k}$ and then express the generating functions with arbitrary ramification profile and descendant insertions as $\star$-polynomials in the $\Trz{k}$.
In Section~4 we relate the building blocks to integrals on DR cycles and to the KdV hierarchy, and deduce the polynomiality results. Section~5 collects some closed formulas and consequences of the reconstruction.

Appendix \ref{app:python} describes an implementation of the recursions, while Appendix \ref{app:ex} records explicit examples of generating functions and reconstruction formulas.
\subsection{Disclosure of AI usage}
The mathematical framework and the computations in this project have been conducted without the assistance of LLM's. ChatGPT has been used to express in appropriate symmetric function basis some of the  polynomial expressions obtained, to write the code in input box [9] in Appendix \ref{app:python}, to tighten the exposition in parts of the introduction and for an overall editorial check of the manuscript.

\subsection{Acknowledgements}
The first author thanks Universit\`a degli Studi di Padova and INDAM for their hospitality and support during his sabbatical, and acknowledges support from the Simons foundation through MPS-TSM-00007937. The second and third authors are supported by Universit\`a degli Studi di Padova and are affiliated with the INFN under the national project MMNLP. The second author thanks Colorado State University for the generous hospitality in the fall of 2025. The \texttt{pyBoLaNO} package \cite{pybolano} was of great use when performing numerical experiments. We also wish to thank Alessio Cela, Aaron Pixton and Dhruv Ranganathan for conversations related to the project.

\section{Background}

\subsection{Admissible Covers}

Hurwitz schemes parameterize covers of smooth curves with fixed discrete invariants: genera of the curves, degree, and ramification profiles over all branch points. Moduli spaces of admissible covers compactify Hurwitz schemes insisting that all discrete data be preserved. We don't give here the exact definition of admissible cover, which can be found in the references we are about to provide. There is a plethora of variants, which we briefly recall to orient the reader. We then carefully describe the spaces of admissible covers we are working with. 

In their original incarnation \cite{hm:kd}, admissible covers are singular (in fact not normal) coarse moduli spaces. In \cite{acv03}, the authors interpret the normalization as the fine moduli stacks of twisted stable maps from orbi-curves to the classifying stack $BS^d$, and show these are smooth and proper DM stacks. In this perspective, branch points are marked but their fibers are not. Shortly after, \cite{Jarvis_Kaufmann_Kimura} construct a version of the space where all points in the fibers of branch points are marked. In particular all ramification is marked, but so are unramified points lying above a branch point. By taking appropriate stack quotients, one can choose to ``unmark" some points, obtaining for example a case in which only honest ramification points are marked.

 When the base curve of the covers has genus zero, an orthogonal direction of variants deals with the notion of isomorphism of covers. On the one end of the spectrum one has admissible covers of a parameterized $\bP^1$ (\cite{r:tqft}), where in the open locus two covers are isomorphic when they agree up to an isomorphism of the source curve; on the other we have admissible covers of a rational curve, where one allows also for composing with an isomorphism of the base curve. It's a well-known folklore fact that admissible covers of a rational curve may be viewed inside the parameterized version as the slice where three distinguished branch points are fixed, say to $0,1, \infty$ (see \cite{r:adm} for a more extended discussion). When two branch points play a special role in the theory, it is sometimes advantageous to fix them at $0, \infty$ and have a residual $\mathbb{C}^\ast$ worth of automorphisms of the base to quotient by: this version is often called {\it rubber}, see \cite{OkounkovPandharipande2006Virasoro}. 
 
 Finally, the admissible cover spaces that make most naturally contact with $S_d$-equivariant geometry and representation theory parameterize covers from possibly disconnected source curves. These moduli spaces have multiple connected components, and one can restrict to just a proper subset of them satisfying some further geometric property.  A natural example is to look only at components parameterizing covers from connected source curves, which tend to be more natural in the study of curve geometry. 

We spell out carefully the spaces we are working with here.

\begin{definition}\label{def:ADM}
Let
${\delta}$ be a positive integer, $g$ a non-negative integer. Let $\{\nu_i\}_{i\in [b]}$ be a list of partitions  $\nu_i \vdash \delta'\leq \delta$ and $\ain$ a vector of positive integers adding to $\delta$. By
\[
\adm_g(\delta, \nu_1, \ldots, \nu_b, \ain)
\]
we denote a moduli space of degree $\delta$ admissible covers of a rubber $\bP^1$ by a (connected) curve of genus $g$ such that:
\begin{itemize}
    \item the cover is fully ramified over $0$;
    \item the ramification profile over $\infty$ is given by $\ain$;
    \item there are $b$ additional points on the base where the ramification profile is 
    $(\nu_i, 1^{\delta-|\nu_i|})$, but only the points corresponding to the parts of $\nu_i$ are marked;
    \item all other branch points are required to be simple and their fibers are unmarked.
\end{itemize}

\end{definition}
Informally, the spaces we consider are hybrid between the \cite{acv03} and the \cite{Jarvis_Kaufmann_Kimura} spaces, where all branch points are marked but only part of their fibers are. The covers are automatically connected because of the point of full ramification requested. There are a tautological source morphism $sr:\adm_g(\delta, \nu_1, \ldots, \nu_b, \ain)\to \mgbar{g,1+\sum\ell(\nu_i)+\ell(\ain)}$ and a finite branch morphism $br: \adm_g(\delta, \nu_1, \ldots, \nu_b, \ain)\to \mgbar{0, b+2+r}$, 
where $r=2g-1-\sum_{i=1}^b(|\nu_i|-\ell(\nu_i))+\ell(\ain)$ is the number of remaining simple branch points.

We introduce  moduli spaces of disconnected admissible covers which appear in the sequel.

\begin{figure}[tb]
     \centering
     \resizebox{0.6 \textwidth}{!}{\begin{tikzpicture}[
  x=1cm,y=1cm,
  line cap=round,line join=round,
  dot/.style={circle,draw=black,fill=black,inner sep=1.7pt},
  open dot/.style={circle,draw=black,fill=white,thick,inner sep=1.7pt},
  guide/.style={draw=black,densely dashed,thin},
  strand/.style={draw=black,thick},
  lab/.style={draw=none,fill=white,inner sep=1.5pt},
  brace/.style={draw=black,decorate,decoration={brace,amplitude=5pt}},
  >=Latex
]

\def\xzero{0}
\def\xtwo{5.4}
\def\xinf{8.6}

\newcommand{\ram}[5]{%
  \ifnum#1=1
    \draw[strand] ({#2-#4/2},#3) -- (#2,#3) -- ({#2+#4/2},#3);
    \node[open dot] at (#2,#3) {};
  \else
    \node[dot] at (#2,#3) {};
    \foreach \k in {1,...,#1}{
      \pgfmathsetmacro{\ry}{#3-#5/2+(\k-1)*#5/(#1-1)}
      \draw[strand]
        ({#2-#4/2},\ry) .. controls ({#2-#4/4},\ry) and
        ({#2-#4/5},#3) .. (#2,#3)
        .. controls ({#2+#4/5},#3) and ({#2+#4/4},\ry) ..
        ({#2+#4/2},\ry);
    }
  \fi
}

\newcommand{\torus}[2]{%
  \draw[thick] ({#1-.64},#2)
    .. controls ({#1-.34},{#2-.48}) and ({#1+.34},{#2-.48}) ..
    ({#1+.64},#2);
  \draw[thick] ({#1-.34},{#2-.18})
    .. controls ({#1-.17},{#2+.08}) and ({#1+.17},{#2+.08}) ..
    ({#1+.34},{#2-.18});
}

\draw[thick] (-0.55,0) -- (9.25,0);
\foreach \x in {\xzero,\xtwo,\xinf} \node[dot] at (\x,0) {};
\node[lab,below=7pt] at (\xzero,0) {$0$};
\node[lab,below=7pt] at (\xtwo,0) {$p_2$};
\node[lab,below=7pt] at (\xinf,0) {$\infty$};
\node[lab,right=13pt] at (9.25,0) {$\bP^1$};

\draw[guide] (\xzero,0.18) -- (\xzero,6.90);
\draw[guide] (\xtwo,0.18) -- (\xtwo,7.24);
\draw[guide] (\xinf,0.18) -- (\xinf,8.10);

\draw[->,thick] (4.30,1.45) -- (4.30,0.35);
\node[lab,right=3pt] at (4.30,0.92) {$\varphi$};

\draw[thick]
  (\xzero,3.08)
  .. controls (0.15,4.10) and (1.05,4.28) .. (1.70,4.28)
  .. controls (3.80,4.46) and (\xinf,4.30) .. (\xinf,4.08)
  .. controls (\xinf,3.88) and (7.88,3.50) .. (\xinf,3.25)
  .. controls (\xinf,3.04) and (7.94,2.57) .. (\xinf,2.32)
  .. controls (\xinf,2.10) and (1.75,1.92) .. (0.72,2.22)
  .. controls (0.28,2.36) and (0.05,2.69) .. (\xzero,3.08);
\node[dot] at (\xzero,3.08) {};
\node[lab,left=7pt] at (\xzero,3.08) {$\delta_m$};
\node[dot] at (\xinf,3.25) {};

\ram{4}{\xtwo}{3.30}{1.85}{0.82}
\ram{1}{\xtwo}{2.48}{1.85}{0}
\draw[brace] (6.58,3.75) -- (6.58,2.84);
\node[lab,right=5pt] at (6.58,3.30) {$ \nu_{1,m}$};
\torus{2.70}{3.03}
\node[dot] at (\xinf,4.08) {};
\node[dot] at (\xinf,2.32) {};
\draw[brace] (9.00,4.08) -- (9.00,2.32);
\node[lab,above right=3pt ] at (9.00,4.08) {$C_{g_m}$};
\node[lab,right=6pt] at (9.00,3.20) {$ \ain_{m}$};

\node[draw=none] at (4.30,5.25) {$\vdots$};

\draw[thick]
  (\xzero,6.90)
  .. controls (0.15,7.88) and (1.05,8.22) .. (1.75,8.28)
  .. controls (3.85,8.48) and (\xinf,8.34) .. (\xinf,8.10)
  .. controls (\xinf,7.88) and (7.88,7.56) .. (\xinf,7.30)
  .. controls (\xinf,7.08) and (7.94,6.58) .. (\xinf,6.25)
  .. controls (\xinf,6.02) and (1.75,5.78) .. (0.72,6.05)
  .. controls (0.28,6.20) and (0.05,6.53) .. (\xzero,6.90);
\node[dot] at (\xzero,6.90) {};
\node[lab,left=7pt] at (\xzero,6.90) {$\delta_1$};

\ram{3}{\xtwo}{6.62}{1.55}{0.42}
\ram{3}{\xtwo}{7.24}{1.55}{0.42}
\draw[brace] (6.45,7.50) -- (6.45,6.35);
\node[lab,right=5pt] at (6.45,6.92) {$ \nu_{1,1}$};
\torus{2.25}{6.90}
\torus{3.75}{6.90}

\node[dot] at (\xinf,8.10) {};
\node[dot] at (\xinf,7.30) {};
\node[dot] at (\xinf,6.25) {};
\draw[brace] (9.00,8.10) -- (9.00,6.25);
\node[lab,above right=3pt] at (9.00,8.10) {$C_{g_1}$};
\node[lab,right=6pt] at (9.00,7.18) {$ \ain_{1}$};

\end{tikzpicture}} 
     \caption{Illustration of  covers parameterized by points of $\dadm_{\vec{g}}(\vec{\delta}, \vec{\nu}, \vec{\ain})$. Black dots denote marked points and empty circles points in the fiber of $p_2$ that are not marked.}
     \label{fig:discpic}
 \end{figure}
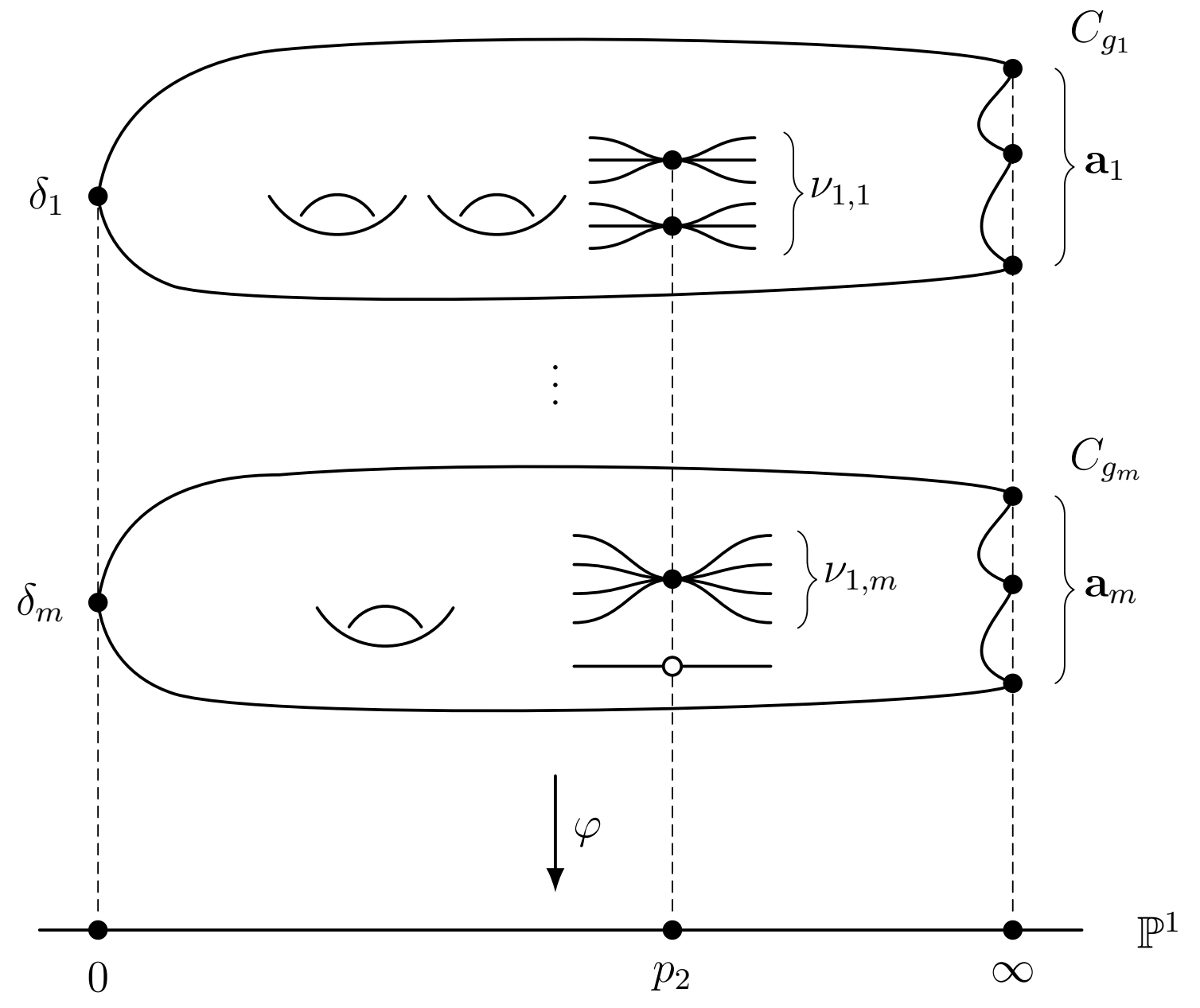

\begin{definition}\label{def:ADMOP}
    Let $\vec{\delta} = (\delta_1, \ldots, \delta_m)$ be a vector of positive integers, $\vec{g} = (g_1, \ldots, g_m)$  a vector of non-negative integers. Let $\vec{{\nu}_i} = ({\nu}_{i,1}, \ldots, {\nu}_{i,m})$ be  vectors of partitions such that  for all $j$, $\nu_{i,j}$ is a partition of $\delta'_j\leq \delta_j$.
    Let $\vec{\ain} = ({\ain}_{1}, \ldots, {\ain}_{m})$ be a vector of partitions such that  for all $j$, $\ain_{j}$ is a partition of $\delta_j$.
We denote by 
\[
\dadm_{\vec{g}}(\vec{\delta}, \vec{\nu}_1,\ldots, \vec{\nu}_b, \vec{\ain})
\]
the union of connected components of the moduli space of disconnected admissible covers $\adm^\bullet_{\vec{g}}(\vec{\delta}, \vec{\nu}_1,\ldots, \vec{\nu}_b, \vec{\ain})$ parameterizing admissible covers where the source curve consists of $m$ connected components; the $j$-th component has genus $g_j$ and the restriction of the cover to the $j$-th component is a connected admissible cover in $\adm_g(\delta_j, {\nu}_{1,j}, \ldots, {\nu}_{b,j}, \ain_j)$, see Figure \ref{fig:discpic} for an illustration. 
\end{definition}

This space admits a morphism
\begin{equation}\label{eq:cmap}
    c:\dadm_{\vec{g}}(\vec{\delta}, \vec{\nu}_1,\ldots, \vec{\nu}_b, \vec{\ain})\to \prod_{j=1}^m \adm_{g_j}(\delta_j, {\nu}_{1,j},\ldots,  {\nu}_{b,j},\ain_j),
\end{equation}
which remembers the individual connected components of the admissible cover and forgets they are part of one disconnected cover. 
\begin{remark}\label{rem:partiallylabeleddisconn} We restrict our attention to the case $b=1$ (and we denote $\vec{\nu}_1$ simply by $\vec{\nu}$), i.e. covers with three special branch points.
If none of the partitions $\nu_{j}$ are empty, then the map $c$ is an isomorphism. This is a consequence of the  fact that there is a unique automorphism of $\bP^1$ sending any three points to any other three points. Any generic admissible cover with $3$ special branch points can be rigidified, i.e. can be made into a cover of a parameterized $\bP^1$ with the three special branch points fixed at $0,1$ and $\infty$.
If however some of the $\nu_{j}$'s are empty, that means that on the corresponding connected components  of the cover none of the points in the fibers are marked. Then such components cannot be rigidified (they don't know of the second branch point) and the map $c$ has positive dimensional fibers. While this statement is correct, there is a subtlety that deserves being carefully discussed: assume that  one of the factors, say the last one,  is a dimension $-1$ stack of the form $\adm_0(\delta_\ell, \delta_\ell)$; denote by $\vec{g'},(\vec{\delta'}, \vec{{\nu'}}, \vec{{\ain'}})$ the datum obtained by dropping the last entries of each of the vectors. Chasing the diagram
\begin{equation}\label{eq:chat}
    \xymatrix{
    \dadm_{\vec{g}}(\vec{\delta}, \vec{{\nu}}, \vec{{\ain}}) \ar[rr]^c \ar[d] \ar[drr]^{\hat{c}}& & \prod_{j=1}^m \adm_{g_j}(\delta_j, {\nu}_{j},  {\ain}_{j}) \ar[d]\\
    \dadm_{\vec{g'}}(\vec{\delta'}, \vec{{\nu'}}, \vec{{\ain'}}) \ar[rr]^{c'}& & \prod_{j=1}^{m-1} \adm_{g_j}(\delta_j, {\nu}_{j},  {\ain}_{j})
    }
\end{equation}
one may see that the fibers of $c'$ have one dimension less than the fibers of $c$ and that the fibers of the diagonal map $\hat{c}$ are  $\nu_d$-gerbes over the corresponding fibers of $c'$.

\end{remark}

\subsection{Tautological classes}

In this section we define the tautological classes on admissible covers we wish to intersect, and recall some of their properties.

\begin{definition}
 Consider the space of admissible covers $\adm_g(\delta, \nu, \ain)$, together with the tautological morphisms $sr,br$. \begin{itemize}
     \item by $\psi_i\in A^1(\adm_g(\delta, \nu, \ain))$ we mean the class $br^\ast(\psi_i)$, i.e. the pullback of the homonymous class on $\overline{\mathcal{M}}_{0,3+r}$ via the branch morphism;
     \item by $\lambda_j\in A^j(\adm_g(\delta, \nu, \ain))$ we mean the class $sr^\ast(\lambda_j)$, i.e. the pullback of the homonymous class on $\overline{\mathcal{M}}_{g,1+\ell(\nu)+\ell(\ain)}$ via the source morphism.
     \item let $\Delta$ denote a (closure of a) stratum in $\adm_g(\delta, \nu, \ain)$; $\Delta$ is isomorphic to a fiber product of moduli spaces of admissible covers, one for each irreducible component of the target curves of the covers generically parameterized by $\Delta$, over some zero-dimensional stacks (see \cite[Section 2.3]{cc:c3z3} for a discussion); we denote by $i_\Delta$ the morphism mapping this fiber product to $\adm_g(\delta, \nu, \ain)$. The class of the stratum  is 
     \begin{equation}
         [\Delta]:= i_{\Delta \ \ast} ([1]).
     \end{equation}
     In practical terms, one thinks of the stratum as the corresponding product of moduli spaces of admissible covers, weighted by a product of gluing factors: every node of the source curves of the covers parameterized contributes its ramification order. 
 \end{itemize} 
 \end{definition}

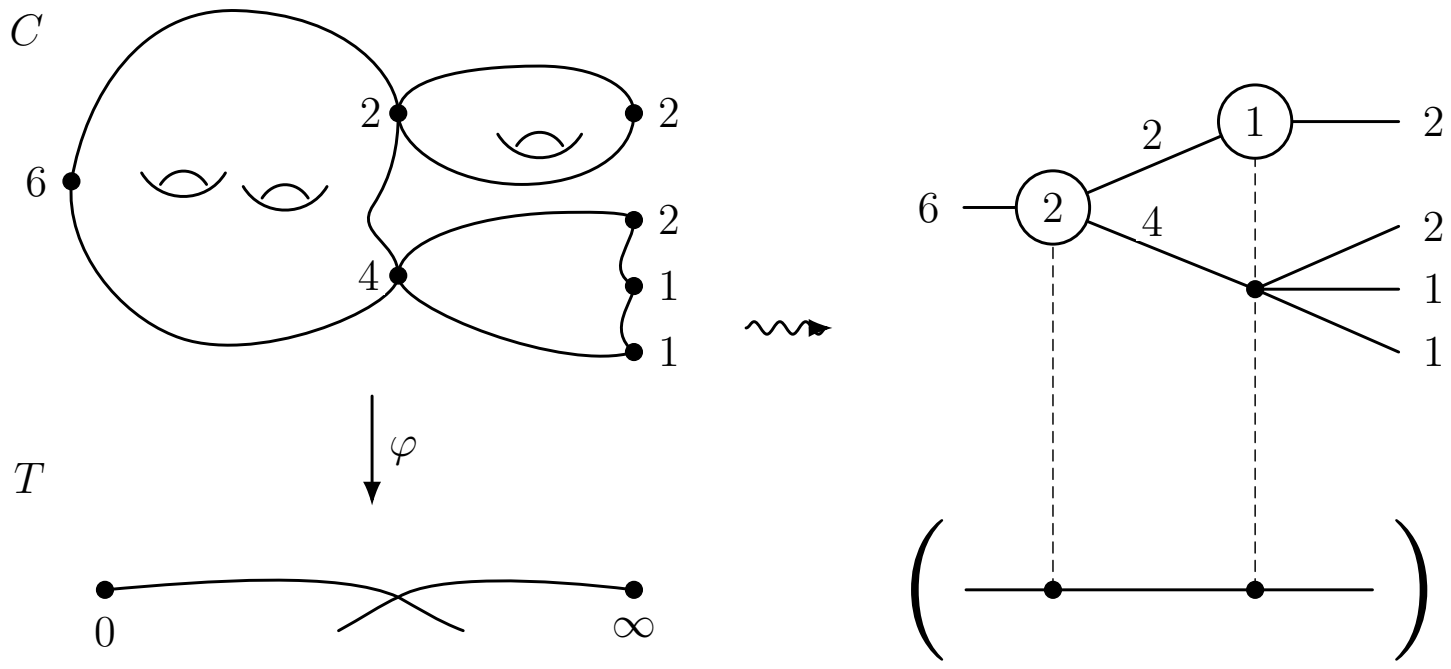
\begin{figure}[tb]
    \centering
    \resizebox{0.6 \textwidth}{!}{\begin{tikzpicture}[
  x=1cm,y=1cm,
  line cap=round,line join=round,
  every path/.style={draw=black},
  dot/.style={circle,draw=black,fill=black,inner sep=1.6pt},
  open dot/.style={circle,draw=black,fill=white,thick,inner sep=1.7pt},
  vertex/.style={circle,draw=black,fill=white,thick,minimum size=7mm,inner sep=0pt},
  guide/.style={densely dashed,thin},
  lab/.style={draw=none,fill=white,inner sep=1.4pt},
  >=Latex
]

\newcommand{\torus}[2]{%
  \draw[thick] ({#1-.40},#2)
    .. controls ({#1-.22},{#2-.30}) and ({#1+.22},{#2-.30}) ..
    ({#1+.40},#2);
  \draw[thick] ({#1-.22},{#2-.11})
    .. controls ({#1-.11},{#2+.06}) and ({#1+.11},{#2+.06}) ..
    ({#1+.22},{#2-.11});
}

\node[lab] at (-0.55,7.90) {$C$};

\coordinate (Ntop) at (3.00,7.10);
\coordinate (Nbot) at (3.00,5.55);
\coordinate (Lsix) at (-0.12,6.45);
\coordinate (Rtop) at (5.25,7.10);
\coordinate (Rtwo) at (5.25,6.08);
\coordinate (RoneA) at (5.25,5.45);
\coordinate (RoneB) at (5.25,4.82);

\draw[thick]
  (Ntop)
  .. controls (3.00,7.42) and (2.55,8.05) .. (1.45,8.08)
  .. controls (0.55,8.10) and (0.02,7.28) .. (Lsix)
  .. controls (-0.18,5.78) and (0.35,5.18) .. (0.85,4.98)
  .. controls (1.60,4.66) and (3.00,5.23) .. (Nbot)
  .. controls (3.00,5.87) and (2.58,5.95) .. (2.76,6.22)
  .. controls (2.94,6.50) and (3.00,6.78) .. (Ntop);

\draw[thick]
  (Ntop)
  .. controls (3.00,7.42) and (3.62,7.55) .. (4.35,7.55)
  .. controls (4.88,7.55) and (5.25,7.35) .. (Rtop)
  .. controls (5.25,6.70) and (4.78,6.42) .. (4.18,6.42)
  .. controls (3.52,6.42) and (3.00,6.78) .. (Ntop);

\draw[thick]
  (Nbot)
  .. controls (3.00,5.87) and (3.70,6.12) .. (4.48,6.15)
  .. controls (4.86,6.16) and (5.25,6.16) .. (Rtwo)
  .. controls (5.25,5.84) and (4.96,5.63) .. (RoneA)
  .. controls (5.25,5.30) and (4.98,5.04) .. (RoneB)
  .. controls (4.60,4.62) and (3.00,5.23) .. (Nbot);

\node[dot] at (Lsix) {};
\node[dot] at (Ntop) {};
\node[dot] at (Rtop) {};
\node[dot] at (Nbot) {};
\node[dot] at (Rtwo) {};
\node[dot] at (RoneA) {};
\node[dot] at (RoneB) {};

\node[lab,left=5pt] at (Lsix) {$6$};
\node[lab,left=3pt] at (Ntop) {$2$};
\node[lab,right=5pt] at (Rtop) {$2$};
\node[lab,left=3.5pt] at (Nbot) {$4$};
\node[lab,right=5pt] at (Rtwo) {$2$};
\node[lab,right=5pt] at (RoneA) {$1$};
\node[lab,right=5pt] at (RoneB) {$1$};

\torus{0.95}{6.53}
\torus{1.92}{6.40}
\torus{4.35}{6.90}

\draw[->,thick] (2.75,4.40) -- (2.75,3.35);
\node[lab,right=3pt] at (2.75,3.88) {$\varphi$};
\node[lab] at (-0.52,3.62) {$T$};

\coordinate (Tzero) at (0.20,2.55);
\coordinate (Tinf)  at (5.25,2.55);
\coordinate (Tint)  at (3.00,2.48);
\draw[thick]
  (Tzero) .. controls (1.20,2.64) and (2.45,2.73) .. (Tint)
  .. controls (3.18,2.38) and (3.38,2.25) .. (3.62,2.16);
\draw[thick]
  (2.43,2.16) .. controls (2.65,2.28) and (2.84,2.41) .. (Tint)
  .. controls (3.35,2.69) and (4.35,2.66) .. (Tinf);
\node[dot] at (Tzero) {};
\node[dot] at (Tinf) {};
\node[lab,below=6pt] at (Tzero) {$0$};
\node[lab,below=6pt] at (Tinf) {$\infty$};

\draw[->,thick,decorate,
  decoration={snake,amplitude=1.7pt,segment length=7pt}]
  (6.32,5.05) -- (7.15,5.05);

\coordinate (Vtwo) at (9.25,6.20);
\coordinate (Vone) at (11.18,7.02);
\coordinate (Vbranch) at (11.18,5.42);

\draw[thick] (8.40,6.20) -- (Vtwo);
\draw[thick] (Vtwo) -- (Vone);
\draw[thick] (Vone) -- (12.55,7.02);
\draw[thick] (Vtwo) -- (Vbranch);
\draw[thick] (Vbranch) -- (12.55,6.02);
\draw[thick] (Vbranch) -- (12.55,5.42);
\draw[thick] (Vbranch) -- (12.55,4.82);

\node[vertex] at (Vtwo) {$2$};
\node[vertex] at (Vone) {$1$};
\node[dot] at (Vbranch) {};

\node[lab,left=5pt] at (8.40,6.20) {$6$};
\node[lab,above=2pt] at (10.20,6.63) {$2$};
\node[lab,right=5pt] at (12.55,7.02) {$2$};
\node[lab,above=2pt] at (10.20,5.80) {$4$};
\node[lab,right=5pt] at (12.55,6.02) {$2$};
\node[lab,right=5pt] at (12.55,5.42) {$1$};
\node[lab,right=5pt] at (12.55,4.82) {$1$};

\coordinate (Qtwo) at (9.25,2.55);
\coordinate (Qbranch) at (11.18,2.55);
\draw[guide] (Qtwo) -- (9.25,5.84);
\draw[guide] (Qbranch) -- (11.18,6.66);
\draw[thick] (8.42,2.55) -- (12.30,2.55);
\node[dot] at (Qtwo) {};
\node[dot] at (Qbranch) {};

\node[draw=none,scale=3.2] at (8.02,2.55) {$($};
\node[draw=none,scale=3.2] at (12.70,2.55) {$)$};

\end{tikzpicture}} 
    \caption{An example of the dual graph notation for strata, in the space $\adm_3((6), (2^2,1^2))$: if a stratum generically parameterizes covers as in the left hand side of the picture, the dual graph of the cover  is depicted on the right. When the target degeneration is just a chain degeneration, we  denote the corresponding stratum just by the dual graph of the cover curve, where the  vertical alignment of the vertices keeps track of the target degeneration. The left part of the graph consists of the vertex of genus $2$ with a left end of weight $6$ and two right ends of weights $2$ and $4$. We omit labeling the leafs of the cover graph to avoid burdening the notation, but remember that indeed points are marked; so there are actually two strata corresponding to this dual graph, corresponding to the two non-isomorphic ways of labeling the points of ramification order $2$ over $\infty$.}
    \label{fig:stratum}
 \end{figure}

 \begin{definition}
     We denote strata by using their {\bf dual graph}, as shown in Figure \ref{fig:stratum}. The dual graph of an admissible cover consists of the natural map from the dual graph of the source curve to the dual tree of the target. We denote any special ramification with additional half edges and label them with their ramification orders. We also decorate compact edges of the dual graph of the cover with the ramification order of the corresponding shadows of the node. We call either a segment or a two-valent tree a {\bf trivial graph}. The dual graph of a boundary divisor, which we call a {\bf divisor graph}, covers a tree with two vertices. By convention, we have the flag corresponding to the first branch point ($0$ in the DR/rubber notation) incident to the left vertex, and the flag corresponding to the third branch point (i.e. $\infty$) incident to the right vertex. By {\bf left part} (resp. {\bf right}) of a divisor graph we mean the inverse image of the left (resp. right) vertex with all its incident half-edges; it identifies the factor of the boundary divisor corresponding to the covers of the $\bP^1$ containing the point $0$ (resp. $\infty$).  
 \end{definition}
The definitions are analogous for the disconnected spaces $\dadm_{\vec{g}}(\vec{\delta}, \vec{{\nu}}, \vec{{\ain}})$. The Hodge bundle on the moduli space of disconnected curves equals the direct sum of the Hodge bundles on the factors corresponding to the connected components of the curves, hence $\lambda$ classes on disconnected moduli spaces may be expressed in terms of $\lambda$ classes on the factors. In particular, for the top Chern class of the Hodge bundle on $\dadm_{\vec{g}}(\vec{\delta}, \vec{{\nu}}, \vec{{\ain}})$ we have
\begin{equation}
    \lambda_{|\vec{g}|} = \prod_{j=1}^{\ell(\vec{g}) }\lambda_{g_j}
\end{equation}

We recall now some well known properties of these classes:
\begin{description}
    \item[boundary expression for $\psi$ classes] fix three distinct marks $i,j,k$ in the set $[3+r]$; the class $\psi_i$ may be represented by the sum of boundary divisors separating the mark $i$ from the marks $j,k$ (\cite[Eq. (9)]{Getzler1998}):
    \begin{equation}\label{eq:psibound}
        \psi_i = \sum_{K\subset[3+r]\smallsetminus\{i,j,k\}} D(\{i\}\cup K^c| K\cup\{j,k\}) \in A^1(\overline{\mathcal{M}}_{0,3+r});
    \end{equation}
    pulling back \eqref{eq:psibound} via $b$ one may represent the class $\psi_i \in A^1(\adm_g(\delta, \nu, \ain))$ as a linear combination of boundary divisors.
    \item[vanishing of $\lambda_g$ off compact type] denoting by $\delta_0$ the irreducible divisor in $\overline{\mathcal{M}}_{g,1+\ell(\nu)+\ell(\ain)}$, we have (\cite[\S 1.6, p.340]{p:kappa})
   \begin{equation}
       \lambda_g\cdot \delta_0 = 0;
   \end{equation}
   by the projection formula one sees that $sr^\ast(\lambda_g)$ vanishes when restricted to a cycle of admissible covers supported on source curves not of compact type. 
    \item[splitting of $\lambda_g$ on compact type strata] for a stratum $\Delta$ generically parameterizing covers where the source curve is a nodal curve of compact type with irreducible components of genus $g_1, \ldots, g_m$. The restriction of the Hodge bundle splits \cite[p.183]{fp:hiagwt}
\[
\mathbb{E}_{g|\Delta}\cong \mathbb{E}_{g_1}\boxplus, \ldots \boxplus \mathbb{E}_{g_m},
\]
and therefore $\lambda_g$
splits multiplicatively:
\begin{equation}\label{eq:lasplit}
    \lambda_g \cdot [\Delta] = i_{\Delta \ \ast}\left(\prod_{i =1}^m \lambda_{g_i}\right)
\end{equation}
\end{description}
 
\begin{remark}
    Because of how they split along boundary strata, lambda classes are recorded on the dual graph indexing a stratum as vertex decorations, whereas psi classes as edge decorations (see \cite[\S 0.3]{jppz17}).
\end{remark}


\subsection{Invariants and generating functions}\label{subsec:inv}
In this section we introduce the generating functions that encode the invariants we wish to study. We choose to present them in order of increasing generality, and in the same order that they will be reconstructed, rather than making one general definition and viewing everything as appropriate specializations. While this strategy is somewhat less efficient, we believe it is more digestible.

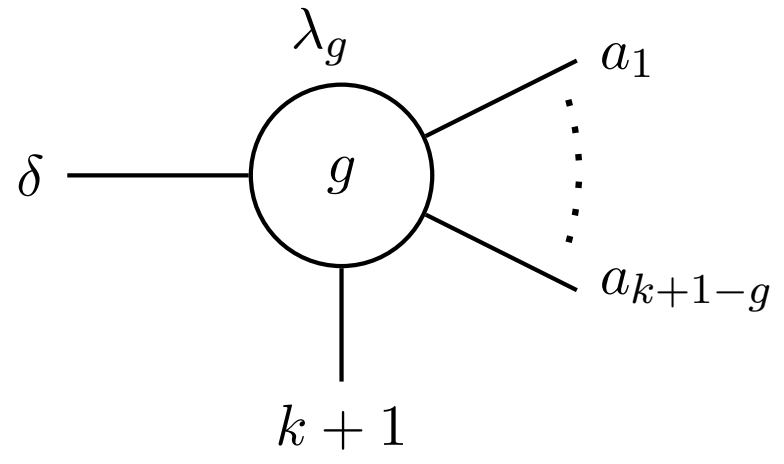
\begin{figure}
    \centering
    \begin{tikzpicture}[thick, font=\small]
        \node[draw, circle, minimum size=1.2cm] (g) at (0,0) {$g$};
        \node[above left] at (g.430) {$\lambda_{g}$};
        \draw (g.180) -- ++(-1.2, 0) node[left] {$\delta$};
        \draw (g.25) -- ++(1, 0.5) node[right] {$a_1$};
        \draw (g.-25) -- ++(1, -0.5) node[right] {$a_{k+1-g}$};
        \draw[loosely dotted, very thick] (1.5, 0.5) to[bend left=15] (1.5, -0.45);
        \draw (g.-90) -- ++(-0.0, -0.75) node[below] {$k+1$};
    \end{tikzpicture}
    \caption{The graph representing the  invariant $x_{k;g}(\mathbf{a})$; the marked ramification points are represented by half edges, with their orientations suggesting which branch point they lay over.}
    \label{fig:dualofinv}
\end{figure}

We first introduce $\lambda_g$ integrals on admissible cover spaces of the form $\adm_g(\delta,(k+1),\ain)$. We write these as
\begin{equation}\label{eq:xkg}
    x_{k;g}(\ain) = \int_{\adm_g\left(\delta,\,(k+1),\,\ain\right)}\lambda_g.
\end{equation}
Notice that $\delta\ge k+1$ is required for \eqref{eq:xkg} to be nonzero. The dimension of the space $\adm_g(\delta,(k+1),\ain)$ is the number of additional simple branch points of the cover, which by the Riemann-Hurwitz formula is $r=2g-1+n-k$. In order for $x_{k;g}(\ain)$ to not vanish, this must equal $g$, so $n=k+1-g$ and $0\le g\le k$. The invariant is represented by the decorated graph in Figure \ref{fig:dualofinv}.
The $k=1$ invariants are
\begin{equation}\label{eq:k1inv}
    x_{1;0}(a_1,a_2) =
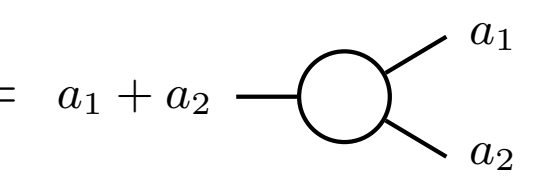
    \begin{tikzpicture}[baseline=-0.5ex, thick, scale=0.6, font=\scriptsize, vertex/.style={draw, circle, minimum size=0.6cm, inner sep=0pt}]
        \node[vertex] (v) at (0,0) {};
        \draw (v.west) -- ++(-0.67, 0) node[left] {$a_1+a_2$};
        \draw (v.30) -- ++(0.67, 0.4) node[right] {$a_1$};
        \draw (v.-30) -- ++(0.67, -0.4) node[right] {$a_2$};
    \end{tikzpicture}
    = 1, \quad x_{1;1}(a) =
    
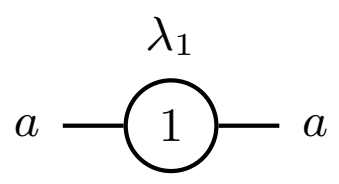
\begin{tikzpicture}[baseline=-0.5ex, thick, scale=0.6, font=\scriptsize, vertex/.style={draw, circle, minimum size=0.6cm, inner sep=0pt}]
        \node[vertex] (v) at (0,0) {$1$};
        \node[above] at (v.90) {$\lambda_1$};
        \draw (v.west) -- ++(-0.67, 0) node[left] {$a$};
        \draw (v.east) -- ++(0.67, 0) node[right] {$a$};
    \end{tikzpicture}
    = \frac{a^2-1}{12}.
\end{equation}

We encode the invariants \eqref{eq:xkg} in a generating function. To this end, we introduce a set of formal variables $\{p_a:a\in\mathbb{N}\}$ which keep track of the ramification profile over $\infty$, and another set of variables $\{q_a:a\in\mathbb{N}\}$ which keep track of the ramification profile over $0$. For a ramification profile $\ain=(a_1,\dots,a_n)\in\mathbb{N}^n$, we use the short-hand notation
\begin{equation*}
    p_{\ain} := \prod_{i=1}^np_{a_i}.
\end{equation*}
We also introduce three more variables $\h,\e^2,v$. The variable $\h$ keeps track of (a multiple of) the Euler characteristic of the cover, $\e^2$ keeps track of the geometric genus, and $v$ keeps track of the number of additional simple branch points of the cover.
\begin{definition}\label{def:Trz}
    For $k\ge1$, define $\Trz{k}\in\mathbb{C}[\h,\h^{-1},\e,v,q_*][[p_*]]$ by
    \begin{equation*}
        \Trz{k} = \sum_{g=0}^{k}\h^{g-1}\e^{2g}\frac{v^g}{g!}\sumn{\ain}{k+1-g}x_{k;g}(\ain)\frac{p_{\ain}}{(k+1-g)!}q_{\lvert\ain\rvert}.
    \end{equation*}
\end{definition}
Two of the three variables $\h,\e^2,v$ are redundant, but it is useful to keep them for when we consider more general invariants. For $k=1$,
\begin{equation}
    \Trz{1} = \frac{1}{\h}\sum_{\ain\in\mathbb{N}^2}\frac{p_\ain}{2!}q_{\lvert\ain\rvert} + \e^2v\sum_{a\in\mathbb{N}}\frac{a^2-1}{12}p_aq_a.
\end{equation}
Next we  allow the ramification over the second branch point to spread out over more points of the fiber. For a partition $\mu=(\mu_1,\dots,\mu_\ell)$ we denote 
\begin{equation*}
    (\mu+1) = (\mu_1+1,\dots,\mu_\ell+1).
\end{equation*}
We define the invariant 
\begin{equation}
    x_{\mu;g}(\ain) := \int_{\adm_g\left(\delta,\,(\mu+1),\,\ain\right)}\lambda_g,
\end{equation}
where $0\le g\le\lvert\mu\rvert$ and $\ain\in\mathbb{N}^{\lvert\mu\rvert+1-g}$. Notice that $\delta\ge\lvert\mu\rvert+\ell$ is required for the invariant to be nonzero. We represent graphically this invariant by a one-vertex graph, which in this case has $\ell$ downward ends labeled by the parts of $\mu+1$.

We are also interested in invariants on the spaces of disconnected admissible covers introduced in Definition \ref{def:ADMOP}, which we define by
\begin{equation}\label{disc}
   x_{\vec{\mu};\vec{g}}(\vec\ain):= \int_{\dadm_{\vec{g}}(\vec{\delta},\,(\vec{\mu}+1),\,\vec{\ain})}\lambda_g.   
\end{equation}

\begin{definition}\label{def:Tmu}
We define the generating function for connected invariants $x_{\mu;g}(\ain)$:
\begin{equation*}
    \Tr{\mu}^\circ = \sum_{g=0}^{\lvert\mu\rvert}h^{g-1}\e^{2g}\frac{v^g}{g!}\sum_{\ain\in\mathbb{N}^{\lvert\mu\rvert+1-g}}\frac{x_{\mu;g}(\ain)}{\lvert\operatorname{Aut}(\mu)\rvert}\frac{p_\ain}{(\lvert\mu\rvert+1-g)!}q_{\lvert\ain\rvert},
\end{equation*}
and for disconnected invariants $x_{\vec{\mu};\vec{g}}(\vec\ain)$:
    \begin{align*}
        \Tr{\mu}^\bullet =  & \; \sum_{m=1}^{\ell(\mu)}\frac{1}{m!}\sum_{\substack{\vec{\mu} =\mu_1,\dots,\mu_m\\ \coprod_i\mu_i=\mu}}\sum_{\substack{\vec{g}=g_1,\dots,g_m\ge0\\g_1+\cdots+g_m =:g}}\h^{g-m}\e^{2g}\frac{v^g}{g!}\binom{g}{g_1,\dots,g_m}\sum_{\ain_1\in\mathbb{N}^{\lvert\mu_1\rvert+1-g_1}}\cdots\sum_{\ain_m\in\mathbb{N}^{\lvert\mu_m\rvert+1-g_m}} \\
        & \;x_{\vec\mu;\vec{g}}(\vec\ain) \prod_{i=1}^m\frac{p_{\ain_i}}{\lvert\operatorname{Aut}(\mu_i)\rvert(\lvert\mu_i\rvert+1-g_i)!}q_{\lvert\ain_i\rvert}.
    \end{align*}
\end{definition}

We give some structural results about these invariants.
\begin{lem}\label{lem:disceval}
    Consider discrete data  $(\vec{\mu}, \vec{g},\vec\ain )=\{(\mu_i,g_i,\ain_i)\}_{i=1}^m$ defining a disconnected invariant above. 
    \begin{enumerate}
        \item If $\mu_i=\emptyset$ and $(g_i,\ain_i)\neq(0,(\delta_i))$ for some $i$, 
        \begin{equation*}
            \int_{\dadm_{\vec{g}}(\vec{\delta},\,(\vec{\mu}+1),\,\ain)}\lambda_g = 0.
        \end{equation*}
        \item  If $(\mu_i,g_i,\ain_i)=(\emptyset,0,(\delta_i))$ for some $i$, denote by $\{(\mu_j^\prime,g_j^\prime,\ain_j^\prime)\}_{j=1}^{m-1}$ the datum obtained by dropping the $i$-th entry. Then
        \begin{equation*}
            \int_{\dadm_{\vec{g}}(\vec{\delta},\,(\vec{\mu}+1),\,\ain)}\lambda_g = \frac{1}{\delta_i}\int_{\dadm_{\vec{g}}(\vec{\delta}^\prime,\,(\vec{\mu}^\prime+1),\,\ain^\prime)}\lambda_g.
        \end{equation*}
        \item If $\mu_i\neq\emptyset$ for all $i=1,\dots,m$, \begin{equation*}
            \int_{\dadm_{\vec{g}}(\vec{\delta},\,(\vec{\mu}+1),\, \vec{\ain})}\lambda_g = \prod_{i=1}^m \int_{\adm_{g_i}(\delta_i,\mu_i,\ain_i)}\lambda_{g_i}.
        \end{equation*}
    \end{enumerate}
\end{lem}
\begin{proof}
    The results stated here follow  from the observations in  Remark \ref{rem:partiallylabeleddisconn}. Consider the map
    $c:\dadm_{\vec{g}}(\vec{\delta},(\vec{\mu}+1), \vec{\ain})\to \prod_{j=1}^m \adm_{g_i}(\delta_i,\mu_i,\ain_i)$ from \eqref{eq:cmap}. Since the Hodge bundle on the space of disconnected covers is pulled back via $c$, by the projection formula we have:
    \begin{equation}\label{eq:prefubini}
        \int_{\dadm_{\vec{g}}(\vec{\delta},\,(\vec{\mu}+1),\,\vec{\ain})}\lambda_g =  \int_{\prod_{i=1}^m \adm_{g_i}(\delta_i,\mu_i,\ain_i)}c_\ast([1])\cdot \boxtimes_{i=1}^m\lambda_{g_i}.
    \end{equation}
    Case $(1)$ follows from \eqref{eq:prefubini} since $c$ has positive dimensional fibers and therefore $c_\ast([1]) = 0$; in case $(3)$ the map $c$ is an isomorphism, and therefore the result follows by applying Fubini's theorem.
    Case $(2)$ is complicated by the presence of some negative dimensional  stacks among the factors of the target space. Given any one such factor, we can drop it and apply pushforward via the map $\hat{c}$ from \eqref{eq:chat}, thus obtaining the desired equality.
\end{proof}

As a consequence of Lemma \ref{lem:disceval}, we may assume that no partition $\mu_i\in \vec{\mu}$ is empty, and the generating functions $\Tr{\mu}^\bullet$ can be obtained from the $\Tr{\mu}^\circ$ via standard exponentiation. 
\begin{cor}\label{cor:expT}
   Let $\mathbb{T}_0^\circ=\sum_{\mu}\Tr{\mu}^\circ$ and $\mathbb{T}_0^\bullet=\sum_{\mu}\Tr{\mu}^\bullet$. Then:
   \begin{equation*}
    \mathbb{T}_0^\bullet = \exp(\mathbb{T}_0^\circ),
\end{equation*}
For every partition $\mu$, this means:
    \begin{equation*}
        \Tr{\mu}^\bullet = \Tr{\mu}^\circ + \sum_{m\ge2}\frac{1}{m!}\sum_{\substack{\mu_1,\dots,\mu_m\\ \coprod_i\mu_i=\mu}}\prod_{i=1}^m\Tr{\mu_i}^\circ.
    \end{equation*}
\end{cor}

We broaden the system of invariants that we consider, by allowing a descendant insertion on the second branch point.  For $d\ge0$, we define
\begin{equation}
    x_{\mu,d;g}(\ain) := \int_{\adm_g\left(\delta,(\mu+1),\,\ain\right)}\lambda_g\psi_2^d.
\end{equation}
We remind the reader that the $\psi$-class is pulled back from the branch morphism. By the Riemann-Hurwitz formula we have $0\le g\le\lvert\mu\rvert+d$ and $\ain\in\mathbb{N}^{\lvert\mu\rvert+d+1-g}$, and by dimension counting the total number of additional branch points is $r=d+g$.
We also consider disconnected invariants; adopting notation from before.
  \begin{equation}
        x_{\vec\mu,d;\vec{g}}(\vec{\ain}) = \int_{\dadm_{\vec{g}}(\vec{\delta},(\vec{\mu}+1),\vec{\ain})}\lambda_g\psi_2^d.
    \end{equation}
  
   The invariant $x_{\vec\mu,d;\vec{g}}(\vec{\ain})$ can be represented by $m$ decorated dual graphss similar to Figure \ref{fig:dualofinv} stacked on top of each other. In this case, the analogues of statements (1) and (3) in Lemma \ref{lem:disceval} are not true. Namely, the disconnected invariants are not the product of the associated connected invariants, and some of the connected components are allowed to be unramified over the second branch point (i.e. the corresponding partition $\mu_i = \emptyset$).  

We now write the generating functions for these invariants. 

\begin{definition}\label{def:Tmud}
   For a partition $\mu$ and an integer $d\ge0$, let
\begin{equation*}
    \Tr{\mu,d}^\circ = \sum_{g=0}^{\lvert\mu\rvert+d}\h^{g-1}\e^{2g}\frac{v^{d+g}}{(d+g)!}\sum_{\ain\in\mathbb{N}^{\lvert\mu\rvert+d+1-g}}\frac{x_{\mu,d;g}(\ain)}{\lvert\operatorname{Aut}(\mu)\rvert}\frac{p_{\ain}}{(\lvert\mu\rvert+d+1-g)!}q_{\lvert\ain\rvert},
\end{equation*}

The generating function for invariants from possibly disconnected sources is:
 
    \begin{align*}
        \Tr{\mu,d}^\bullet =  & \; \sum_{m\ge1}\frac{1}{m!}\sum_{\substack{\mu_1,\dots,\mu_m\\ \coprod_i\mu_i=\mu}}\sum_{g_1,\dots,g_m\ge0}\h^{g-m}\e^{2g}\sum_{\substack{n_1,\dots,n_m\ge1\\(g_i,n_i)\neq(0,1)}}\frac{v^r}{r!}\binom{r}{r_1,\dots,r_m} \\
        & \; \sum_{\ain_1\in\mathbb{N}^{n_1}}\cdots\sum_{\ain_m\in\mathbb{N}^{n_m}}\frac{x_{\vec\mu,d;\vec{g}}(\vec{\ain})}{\prod_i\lvert\operatorname{Aut}(\mu_i)\rvert}\prod_{i=1}^m\frac{p_{\ain_i}}{n_i!}q_{\lvert\ain_i\rvert},
    \end{align*}
\end{definition}

Since the disconnected generating function $\Tr{\mu,d}^\bullet $ is the most general we introduce in this work (all others could in principle be obtained by specialization from this one), we spell out carefully the indexing data for the potentially non-zero coefficients.

\begin{definition}\label{def:Smud}
    For a partition $\mu$ and an integer $d\ge0$, let $S_{\mu,d}$ be the set whose elements $\sigma$ consist of
    \begin{itemize}
        \item an integer $m\ge1$;
        \item for every $i=1,\dots,m$, a \emph{possibly empty} partition $\mu_i$, such that $\coprod_{i=1}^m\mu_i = \mu$;
        \item for every $i=1,\dots,m$, an integer $g_i\ge0$;
        \item for every $i=1,\dots,m$, a positive integer $n_i$ such that $(g_i,n_i)\neq(0,1)$;
        \item for every $i=1,\dots,m$, a vector of positive integers $\ain_i\in\mathbb{N}^{n_i}$.
    \end{itemize}
    For every such $\sigma=(\mu_i,g_i,n_i,\ain_i)_{i=1}^m\in S_{\mu,d}$, let
    \begin{itemize}
        \item $r_i=2g_i-1+n_i-\lvert\mu_i\rvert $ for every $i=1,\dots,m$, and $r=\sum_{i=1}^mr_i$;
        \item $\vec{\mu}=(\mu_1,\dots,\mu_m)$;
        \item $\vec{g}=(g_1,\dots,g_m)$ and $g=\sum_{i=1}^mg_i$;
        \item $n=\sum_{i=1}^mn_i$;
        \item $\vec{\ain}=(\ain_1,\dots,\ain_m)$;
        \item $\delta_i=\lvert\ain_i\rvert$ and $\vec{\delta}=(\delta_1,\dots,\delta_m)$;
        \item $(\vec{\mu}+1)=\bigl((\mu_1+1),\dots,(\mu_m+1)\bigr)$.
    \end{itemize}
    \end{definition}
    The data $\sigma=(\mu_i,g_i,n_i,\ain_i)_{i=1}^m$ corresponds to an admissible cover with $m$ connected components in the same way as in Definition \ref{def:ADMOP}, where the $i$-th component has $r_i$ other simple branch points (by the Riemann-Hurwitz formula). 
    We use the shorthand
\begin{equation} \label{eq:discacshorthand}
    \adm^\bullet_\sigma = \dadm_{\vec{g}}(\vec{\delta}, \vec{\mu}+1, \vec{\ain}).
\end{equation}
    The number $n_i$ of ramification points over $\infty$ is not fixed by $\mu_i$ and $g_i$, and there are a total of $n$ marked points over $\infty$ and a total of $r$ non-special simple branch points. We can now write more compactly:
 \begin{equation}\label{eq:tmudbull}
        \Tr{\mu,d}^\bullet =  \; \sum_{\sigma\in S_{\mu,d}}\frac{1}{m!}\h^{g-m}\e^{2g}\frac{v^r}{r!}\binom{r}{r_1,\dots,r_m}\frac{x_{\vec\mu,d;\vec{g}}(\vec{\ain})}{\prod_i\lvert\operatorname{Aut}(\mu_i)\rvert}\prod_{i=1}^m\frac{p_{\ain_i}}{n_i!}q_{\lvert\ain_i\rvert}
    \end{equation}

\begin{lem}\label{lem:conditions}
    For every $\sigma=(\mu_i,g_i,n_i,\ain_i)_{i=1}^m\in S_{\mu,d}$ the intersection number $x_{\vec\mu,d;\vec{g}}(\vec{\ain})$ vanishes unless all of the following conditions hold:
    \begin{enumerate}
        \item $n=\lvert\mu\rvert+d-g+m$,
        \item $g\le\lvert\mu\rvert + d$,
        \item $\lvert\mu_i\rvert + \ell(\mu_i) \le \delta_i$ for every $i=1,\dots,m$,
        \item $r=d+g$,
        \item $r_i\ge g_i$ for every $i=1,\dots,m$,
        \item $m\le \ell(\mu)+d$.
    \end{enumerate}
\end{lem}
\begin{proof}
    The first condition follows from applying the Riemann-Hurwitz formula to every connected component of the cover. If $g>\lvert\mu\rvert+d$, then $m\le n = \lvert\mu\rvert+d-g+m<m$, so (2) must hold. For the ramification profile $(\mu_i+1)=(\mu_{i,1}+1,\ldots,\mu_{i,\ell(\mu_i)}+1,1,\ldots,1)$ to be valid, the degree $\delta_i$ must be at least $\lvert\mu_i\rvert+\ell(\mu_i)$, hence (3). Part (4) follows from dimension counting. For (5), pushing forward along the morphism $c:\adm_{\sigma}^\bullet\rightarrow\prod_{i=1}^m\adm_{g_i}(\delta_i,(\mu_i+1),\ain_i)$ introduces a $\lambda_{g_i}$ term in the $i$-th factor. Therefore one must have $r_i=\dim \adm_{g_i}(\delta_i,(\mu_i+1),\ain_i)\ge g_i$. Lastly, let $s_i:=g_i-1+n_i-\lvert\mu_i\rvert$, so that $\sum_{i=1}^ms_i=d$. Let $I$ be the set of those $i\in\{1,\dots,m\}$ for which $\mu_i\neq\emptyset$, and let $I_\emptyset$ be the complementary set. Then $\lvert I\rvert\le\ell(\mu)$, and when $i\in I_\emptyset$ we have $s_i\ge1$ since $(g_i,n_i)\neq(0,1)$. Therefore $m=\lvert I\rvert+\lvert I_\emptyset\rvert \le \ell(\mu)+\sum_{i\in I_\emptyset}s_i\le\ell(\mu)+d$, proving (6).
\end{proof}

\subsubsection{String Equation} In this section we introduce one of the recursions needed to reconstruct invariants, commonly known as the string equation (see \cite{k:pc} for a friendly introduction in the context of moduli spaces of curves).

\begin{lem}[String]\label{lem:string}
    Assume the space $\adm_{g}(\delta, \nu, 1^\delta)$ has  dimension $r>0$, and let $g+d = r$. Then
    \[
    \int_{\adm_{g}(\delta, \nu, (1^\delta))} \lambda_g \psi_2^d = \delta! \int_{\adm_{g}(\delta, \nu)} \lambda_g \psi_2^{d-1}
    \]
\end{lem}
\begin{proof}
Consider $\adm_{g}(\delta,\nu,\emptyset)$ the moduli space of admissible covers with two special branch points of profiles $\delta$ and $\nu$ and an additional marked point on the target curve  which is not a branch point  - in other words, we are unlabeling the points corresponding to the partition $(1^\delta)$.
Consider the diagram:
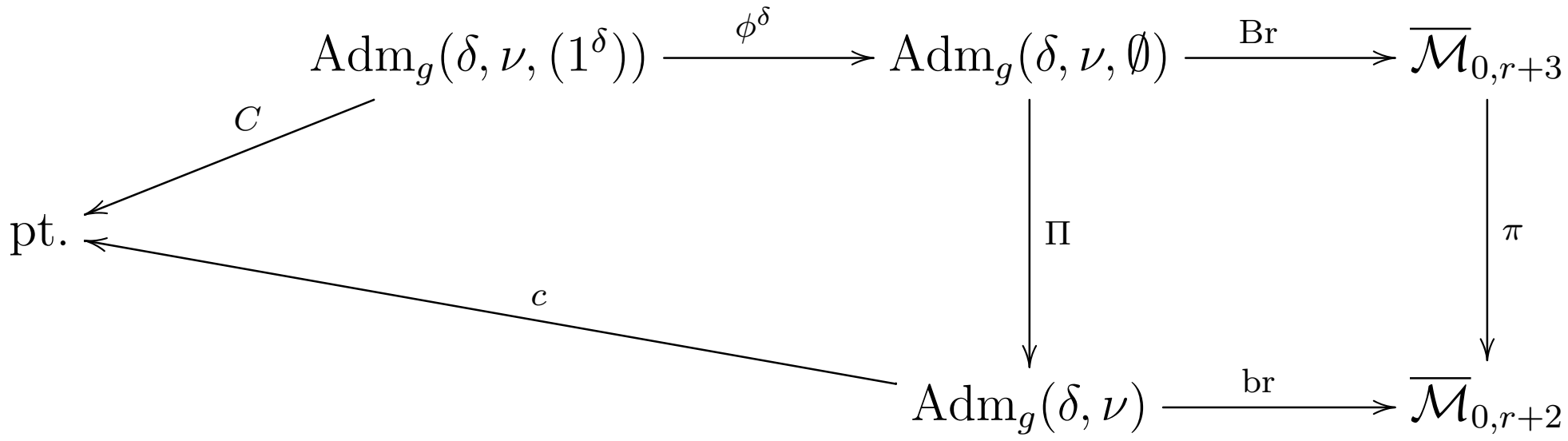
\[    \xymatrix@C=4.2em@R=1.9em{
  & \adm_{g}(\delta, \nu, (1^\delta)) \ar[dl]_{C} \ar[r]^{\phi^\delta}
    & \adm_{g}(\delta, \nu, \emptyset) \ar[r]^{\operatorname{Br}} \ar[dd]^{\Pi}
    & \overline{\mathcal{M}}_{0,r+3} \ar[dd]^{\pi}
  \\
  \mathrm{pt.} & & &
  \\
  &  
    & \adm_{g}(\delta, \nu) \ar[ull]_{c} \ar[r]^{\operatorname{br}}
    & \overline{\mathcal{M}}_{0,r+2}
}\]
We observe:
\begin{description}
    \item[(1)] the map $\phi^\delta$ is finite etale of degree $\delta!$; 
    \item[(2)] the class $\lambda_g$ is stable under pullback via $\Pi$ and $\phi^\delta$;
    \item[(3)] the class we call $\psi_2$ on the spaces of admissible covers is pulled back  via $Br, br$ from the corresponding moduli space of curves;
    \item[(4)] the right square is Cartesian.
\end{description}
Then we have:
\[\begin{array}{lll}
 \int_{\adm_{g}(\delta, \nu, (1^\delta))} \lambda_g \psi_2^d & = C_\ast(\lambda_g \psi_2^d) & =   c_\ast \Pi_\ast \phi^\delta_\ast(\lambda_g \psi_2^d) \\
  &\stackrel{(1,2)}{=}  c_\ast \Pi_\ast \phi^\delta_\ast(\phi^\delta)^\ast(\lambda_g \psi_2^d)&\stackrel{(1)}{=} \delta! c_\ast \Pi_\ast (\lambda_g \psi_2^d)\\
  &\stackrel{(2)}{=} \delta! c_\ast \lambda_g \Pi_\ast ( \psi_2^d)&\stackrel{(3)}{=} \delta! c_\ast \lambda_g \Pi_\ast ( Br^\ast\psi_2^d)\\
  &\stackrel{(4)}{=} \delta! c_\ast \lambda_g br^\ast ( \pi_\ast\psi_2^d) &  \stackrel{(st)}{=} \delta! c_\ast \lambda_g br^\ast ( \psi_2^{d-1})\\
  & \stackrel{(3)}{=} \delta! c_\ast \lambda_g  \psi_2^{d-1} & =\delta! \int_{\adm_{g}(\delta, \nu)} \lambda_g \psi_2^{d-1},
\end{array}\]
where each nontrivial equality is labeled by the motivating fact listed above, with $(st)$ denoting the classical string equation on moduli spaces of stable curves (\cite[equation (1), p. 191]{Looijenga1993}).
\end{proof}
In terms of invariants, Lemma \ref{lem:string} says:
\begin{equation}
    x_{\mu,d;g}(1^\delta) = \delta(\delta-1)  x_{\mu,d-1;g}(2,1^{\delta-2}).
\end{equation}

There is a disconnected version of the string equation, which is formally identical with multiplication factor equal to $\prod_{i=1}^m \delta_i!$; we point out a subtlety on how to apply the  recursion to invariants corresponding to covers with source curves with $m$ connected components. For every $i=1,\dots,m$ for which $r_i>0$ (i.e., there is at least one additional simple branch point $p_i\in\mathbb{P}^1$ associated with the $i$-th connected component), one can ``promote'' $p_i$ to the position previously occupied by $\infty\in\mathbb{P}^1$. Therefore the string equation is equivalent to
\begin{equation}\label{eq:disc_string_1}
    x_{\vec{\mu},d;\vec{g}}\bigl((1^{\delta_1}),\dots,(1^{\delta_m})\bigr) = \delta_i(\delta_i-1) \, x_{\vec\mu,d-1;\vec{g}}\bigl((1^{\delta_1}),\dots,(2,1^{\delta_i-2}),\dots,(1^{\delta_m})\bigr).
\end{equation}
However if $(g_i,\delta_i)=(0,2)$, after forgetting the two marked points and moving the ramification point over to $\infty$,  the $i$-th connected component is a trivial cover of degree $2$, hence:
\begin{equation}\label{eq:disc_string_2}
    x_{\vec\mu,d;\vec{g}}((1^{\delta_1}),\dots,(1^{\delta_m})) = x_{\vec\mu,d-1;\vec{g}}((1^{\delta_1}),\dots,(1^{\delta_{i-1}}),(1^{\delta_{i+1}})\dots,(1^{\delta_m})).
\end{equation}
At the level of generating functions we have the following result.
\begin{lem}\label{lem:Tmu_string}
    The string equation implies
    \begin{equation*}
        \left.\frac{\partial}{\partial v}\Tr{\mu,d}^\bullet\right\rvert_{p_{\ge2=0}} = \left.\left(p_1^2\frac{\partial}{\partial p_2} + \frac{1}{\h}\frac{p_1^2}{2}q_2\right)\Tr{\mu,d-1}^\bullet\right\vert_{p_{\ge2=0}}.
    \end{equation*}
\end{lem}
\begin{proof}
    Setting $p_{\ge2=0}$ on the left hand side selects the invariants in $\Tr{\mu,d}^\bullet$  the string equation applies to. The derivative with respect to $v$ corresponds to the fact that applying the string equation one general simple branch point gets ``promoted" to  becoming $\infty$. The two terms on the right hand side correspond  to the substitutions \eqref{eq:disc_string_1} and \eqref{eq:disc_string_2}.
    \end{proof}

\subsection{WDVV relations and graph calculus}\label{subsec:wdvv}
In this section we recall, set  notation and describe a graph-calculus for WDVV relations; to minimize notational burden, we do this for spaces of connected covers - everything goes through analogously in the disconnected case. Consider a space $\adm_{g}(\delta, \nu,(2),\ain)$ of dimension $r\geq1$.  There is a forgetful morphism:
\begin{equation*}
    f: \adm_{g}(\delta, \nu, (2),\ain) \to \overline{\mathcal{M}}_{0,4},
\end{equation*}
which remembers only the stabilization of the target curve of a cover with the four marked points corresponding to the three special branch points plus $b$. For any Chow class $\alpha$ of degree $r-1$, we have the vanishing of the zero-dimensional cycle class:
\begin{equation}\label{eq:WDVVcyc}
    \alpha\cdot f^\ast(D(\{1,3\}|\{2,4\})-D(\{1,2\}|\{3,4\}))=0.
\end{equation}
Integrating \eqref{eq:WDVVcyc} one obtains a numerical equation among invariants that we refer to as a WDVV relation. We give a combinatorial description of the WDVV relations used to reconstruct invariants in terms of dual graphs.
\begin{lem}[Graphic description of WDVV]\label{lem:graphic_wdvv}
    For an integer $k\geq 2$, consider the space $\adm_{g}(\delta, \nu,(2),\ain)$ with $\ain\in\mathbb{N}^{|
    \nu|+d+2-g}$. Choosing $\alpha = \lambda_g\psi_2^d$ in \eqref{eq:WDVVcyc}, the relation is a sum of terms indexed by divisor graphs $\Gamma$ where the only non-trivial component of either the left or right hand side of the graph consists of a $\Trz{1}$-graph \eqref{eq:k1inv}. For each graph, the contribution is obtained as the product of:
    \begin{itemize}
        \item a factor of $(-1)$ when the $\Trz{1}$-graph is on the right side of the graph;
        \item a factor $\frac{\lvert\operatorname{Aut}(\ain)\rvert}{\lvert\operatorname{Aut}(\Gamma)\rvert}$;
        \item the product of weights of all compact edges of the graph;
        \item a factor of $(g+d)\cdot\frac{a^2-1}{12}$ when the $\Trz{1}$-graph has genus one with ends of weight $a$;
        \item the invariant corresponding to the side of the graph not containing the $\Trz{1}$-graph. 
    \end{itemize}
\end{lem}

\begin{proof}
It may be worth referring to Examples \ref{ex:wdvv1} and \ref{ex:wdvv2} while reading this proof, as they may clarify the general description.
Since the second factor in relation \eqref{eq:WDVVcyc} is a pullback of a linear combination of boundary divisors,  the relation is expressed as a linear combination of strata classes indexed by divisor graphs in the space of admissible covers. Multiplying with the class $\lambda_g$ forces all possibly contributing graphs to be trees. Divisor graphs index strata that are isomorphic to a (stacky fiber) product of two (possibly disconnected) moduli spaces of admissible covers; by dimension reasons, the factor where the branch point $b$ lies parameterizes admissible covers with at most two nonsimple branch points. The only non-trivial connected component of the corresponding dual graph 
 is a $\Trz{1}$-graph. Having narrowed down the class of possibly contributing graphs, we now explain the coefficients. Fix a possibly contributing divisor graph $\Gamma$. The minus sign picks up  which of the boundary divisors of $\overline{\mathcal{M}}_{0,4}$ the stratum corresponding to $\Gamma$ is in the inverse image of. The factor $\lvert\operatorname{Aut}(\ain)\rvert/\lvert\operatorname{Aut}(\Gamma)\rvert$ keeps track of the number of  strata parameterized by the graph $\Gamma$ due to the fact that the inverse images of the third branch point are marked. The product of weights of compact edges is the contribution by ghost automorphisms. The factor of $(g+d)$ is because there besides $b$ there are $(g+d)$ marked simple branch points, one of which is on the $\Trz{1}$ side, and therefore there are $(g+d)$ distinct isomorphic strata. Finally we have the evaluation of the class $\lambda_1$ on one of the factors (when the $\Trz{1}$-graph has genus one), and the evaluation of the class on the second factor. 
\end{proof}

\begin{figure}
    \hspace{-3cm} \makebox[\textwidth][c]{\resizebox{1.2\textwidth}{!}{\definecolor{edgeblue}{RGB}{0,105,210}
\definecolor{markred}{RGB}{210,25,25}
\definecolor{markgreen}{RGB}{98,180,35}
\definecolor{markteal}{RGB}{0,125,105}

\tikzset{
  every path/.style={line cap=round,line join=round},
  fat/.style={circle,draw=black,fill=white,line width=.8pt,
              minimum size=7mm,inner sep=0pt,font=\small},
  point/.style={circle,fill=black,inner sep=1.55pt},
  coloured point/.style={circle,inner sep=2.2pt},
  graph/.style={line width=1.05pt},
  blue graph/.style={graph,draw=edgeblue},
  tick/.style={line width=1.25pt}
}

\newcommand{\darkflag}[1]{%
  \draw[tick,markteal] (#1.south west) -- ++(-.17,-.34);}
\newcommand{\redflag}[1]{%
  \draw[tick,markred] (#1.south) -- ++(0,-.38);}
\newcommand{\lightflag}[1]{%
  \draw[tick,markgreen] (#1.south east) -- ++(.17,-.34);}
\newcommand{\redlightflags}[1]{\redflag{#1}\lightflag{#1}}

\newcommand{\ProfileFour}{%
\begin{tikzpicture}[baseline=-.55ex,x=.75cm,y=.75cm]
  \node[fat] (v) at (0,0) {$4$};
  \node[above left] at (v.north) {$\lambda_4$};
  \draw[graph] (-1.2,0)--(v.west) node[above,xshift=-5mm] {$5$};
  \draw[graph] (v.east)--(.98,.42) node[right=-1pt] {$3$};
  \draw[graph] (v.east)--(.98,-.42) node[right=-1pt] {$2$};
  \redlightflags{v}
  \node[markred,below=13pt] at (v.south) {$5$};
\end{tikzpicture}}

\newcommand{\CrossingCurve}[4]{%
\begin{tikzpicture}[baseline=-.45ex,x=.58cm,y=.58cm]
  \coordinate (L) at (0,0);
  \coordinate (A) at (.95,.48);
  \coordinate (C) at (3.35,.48);
  \coordinate (R) at (4.2,0);
  \draw[graph] (L)--(2.80,1.42);
  \draw[graph] (1.40,1.42)--(R);
  \node[point] at (L) {};
  \node[point] at (R) {};
  \node[coloured point,fill=#1] at (A) {};
  \node[coloured point,fill=#2] at (C) {};
  \node[text=#1,below=5pt] at (A) {$#3$};
  \node[text=#2,below=5pt] at (C) {$#4$};
  \node[below=3pt] at (L) {$0$};
  \node[below=3pt] at (R) {$\infty$};
\end{tikzpicture}}

\newcommand{\Bracketed}[1]{\left[\,#1\,\right]}

\newcommand{\mathpic}[1]{\mathord{\vcenter{\hbox{#1}}}}

\newcommand{\ProfileOneThree}{%
\begin{tikzpicture}[baseline=-.65ex,x=.68cm,y=.68cm]
  \node[fat,draw=edgeblue,text=black] (a) at (0,0) {$1$};
  \node[above=2pt,text=edgeblue] at (a.north) {$\lambda_1$};
  \node[fat] (b) at (2.0,0) {$3$};
  \node[above=2pt] at (b.north) {$\lambda_3$};
  \draw[blue graph] (-1.2,0)--(a.west)
     node[pos=.18,above=1pt,text=edgeblue] {$5$};
  \coordinate (abmid) at ($(a.east)!.5!(b.west)$);
  \draw[blue graph] (a.east)--(abmid);
  \draw[graph] (abmid)--(b.west);
  \node[above=1pt] at (abmid) {$5$};
  \draw[graph] (b.east)--(3.05,.45) node[right=-1pt] {$3$};
  \draw[graph] (b.east)--(3.05,-.45) node[right=-1pt] {$2$};
  \lightflag{a}
  \redflag{b}
  \node[markred,below=12pt] at (b.south) {$5$};
\end{tikzpicture}}

\newcommand{\ProfileThreeFork}{%
\begin{tikzpicture}[baseline=-.7ex,x=.64cm,y=.64cm]
  \node[fat] (v) at (0,0) {$3$};
  \node[above=2pt] at (v.north) {$\lambda_3$};
  \draw[graph] (-1.05,0)--(v.west) node[pos=.18,above=1pt] {$5$};
  \node[fat,draw=edgeblue] (u) at (1.65,.55) {$1$};
  \node[above=1pt,text=edgeblue] at (u.north) {$\lambda_1$};
  \coordinate (vumid) at ($(v.25)!.5!(u.west)$);
  \draw[graph] (v.25)--(vumid);
  \draw[blue graph] (vumid)--(u.west);
  \node[above=1pt,text=edgeblue] at (vumid) {$3$};
  \draw[blue graph] (u.east)--(2.75,.55) node[right=-1pt,text=edgeblue] {$3$};
  \coordinate (vturn) at (.82,-.48);
  \draw[graph] (v.-28) .. controls (.42,-.28) and (.62,-.43) .. (vturn);
  \draw[blue graph] (vturn) .. controls (1.25,-.66) and (1.82,-.70) ..
     (2.75,-.58) node[right=-1pt,text=edgeblue] {$2$};
  \redflag{v}
  \lightflag{u}
  \node[markred,below=12pt] at (v.south) {$5$};
\end{tikzpicture}}

\newcommand{\ProfileThreeOne}{%
\begin{tikzpicture}[baseline=-.7ex,x=.64cm,y=.64cm]
  \node[fat] (v) at (0,0) {$3$};
  \node[above=2pt] at (v.north) {$\lambda_3$};
  \draw[graph] (-1.05,0)--(v.west) node[pos=.18,above=1pt] {$5$};
  \node[fat,draw=edgeblue] (u) at (1.75,-.02) {$1$};
  \node[above=2pt,text=edgeblue] at (u.north) {$\lambda_1$};
  \coordinate (vumid) at ($(v.east)!.5!(u.west)$);
  \draw[graph] (v.east)--(vumid);
  \draw[blue graph] (vumid)--(u.west);
  \node[above=1pt,text=edgeblue] at (vumid) {$2$};
  \coordinate (vturn) at (.95,1.45);
  \draw[graph] (v.38) .. controls (.48,.82) and (.68,1.28) .. (vturn);
  \draw[blue graph] (vturn) .. controls (1.42,1.84) and (2.15,1.62) ..
     (2.75,.72) node[right=-1pt,text=edgeblue] {$3$};
  \draw[blue graph] (u.east)--(2.75,-.42) node[right=-1pt,text=edgeblue] {$2$};
  \redflag{v}
  \lightflag{u}
  \node[markred,below=12pt] at (v.south) {$5$};
\end{tikzpicture}}

\newcommand{\ProfileFourBlue}{%
\begin{tikzpicture}[baseline=-.7ex,x=.67cm,y=.67cm]
  \node[fat] (v) at (0,0) {$4$};
  \node[above=2pt] at (v.north) {$\lambda_4$};
  \draw[graph] (-1.05,0)--(v.west) node[pos=.18,above=1pt] {$5$};
  \coordinate (s) at (.82,0);
  \coordinate (vsmid) at ($(v.east)!.5!(s)$);
  \draw[graph] (v.east)--(vsmid);
  \draw[blue graph] (vsmid)--(s);
  \node[above=1pt] at (vsmid) {$5$};
  \draw[blue graph] (s)--(1.9,.48) node[right=-1pt,text=edgeblue] {$3$};
  \draw[blue graph] (s)--(1.9,-.48) node[right=-1pt,text=edgeblue] {$2$};
  \fill[edgeblue] (s) circle (1.5pt);
  \redflag{v}
  \draw[tick,markgreen] (s)--++(0,-.42) coordinate (sb);
  \node[markred,below=12pt] at (v.south) {$5$};
\end{tikzpicture}}

\newcommand{\ProfileThree}{%
\begin{tikzpicture}[baseline=-.55ex,x=.60cm,y=.60cm]
  \node[fat] (v) at (0,0) {$3$};
  \node[above=2pt] at (v.north) {$\lambda_3$};
  \draw[graph] (-1.02,0)--(v.west) node[pos=.18,above=1pt] {$5$};
  \draw[graph] (v.east)--(.95,.38) node[right=-1pt] {$3$};
  \draw[graph] (v.east)--(.95,-.38) node[right=-1pt] {$2$};
  \redflag{v}
  \node[markred,below=11pt] at (v.south) {$5$};
\end{tikzpicture}}

\newcommand{\ProfileFourSimple}{%
\begin{tikzpicture}[baseline=-.55ex,x=.60cm,y=.60cm]
  \node[fat] (v) at (0,0) {$4$};
  \node[above=2pt] at (v.north) {$\lambda_4$};
  \draw[graph] (-1.02,0)--(v.west) node[pos=.18,above=1pt] {$5$};
  \draw[graph] (v.east)--(.95,0) node[right=-1pt] {$5$};
  \redflag{v}
  \node[markred,below=11pt] at (v.south) {$5$};
\end{tikzpicture}}

\begin{minipage}{22.5cm}
\centering
\setlength{\jot}{1.6em}
\begin{gather*}
  \mathpic{\ProfileFour}\,\cdot f^{*}\!\left(
    \Bracketed{\mathpic{\CrossingCurve{markgreen}{markred}{b}{5}}}
    -
    \Bracketed{\mathpic{\CrossingCurve{markred}{markgreen}{5}{b}}}
  \right)
  ={}
  \\[1.2em]
  \mathpic{\ProfileOneThree}
  \quad - \hspace{0.5em}
  \left(
      \mathpic{\ProfileThreeFork}
      \quad + \hspace{0.5em} 
      \mathpic{\ProfileThreeOne}
      \quad + \hspace{0.5em} 
      \mathpic{\ProfileFourBlue}
    \right)
  {}={}
  \\[1.2em]
  {\color{markteal}4}\cdot\left(
    5\cdot\frac{5^{2}-1}{12}
    -3\cdot\frac{3^{2}-1}{12}
    -2\cdot\frac{2^{2}-1}{12}
  \right)\cdot\mathpic{\ProfileThree}
  \quad - \hspace{0.5em} 
  5\cdot\mathpic{\ProfileFourSimple}
  {}
\end{gather*}
\end{minipage}}}
    \caption{The WDVV relation from Example \ref{ex:wdvv1}.}
    \label{fig:wdvv1}
\end{figure}
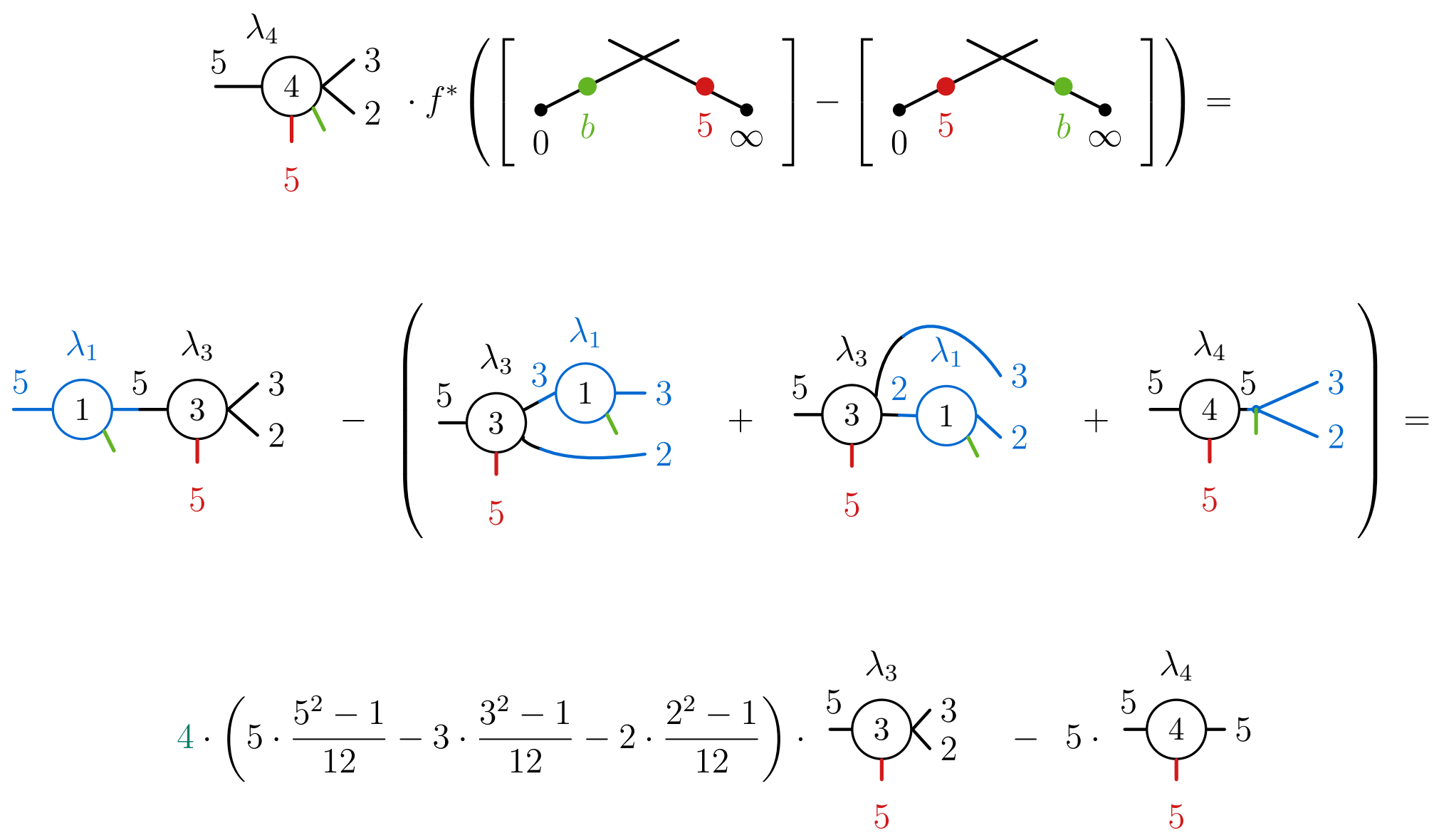
\begin{ex} \label{ex:wdvv1}
    We consider the WDVV relation on the moduli space $\adm_{4}((5),(5),(3,2))$ corresponding to $\alpha=\lambda_4$. There is only one possible divisor graph which has a $\Trz{1}$-graph on the left, and three that have it on the right, as depicted  in the second line of Figure \ref{fig:wdvv1}. The third line collects the various multiplicities and graphically depicts the relation
    \begin{equation*}
        6x_{4;3} (3,2)-x_{4;4}(5)=0.
    \end{equation*}
\end{ex}

\begin{figure}
    \hspace{-3cm} \makebox[\textwidth][c]{\resizebox{1.35\textwidth}{!}{\input{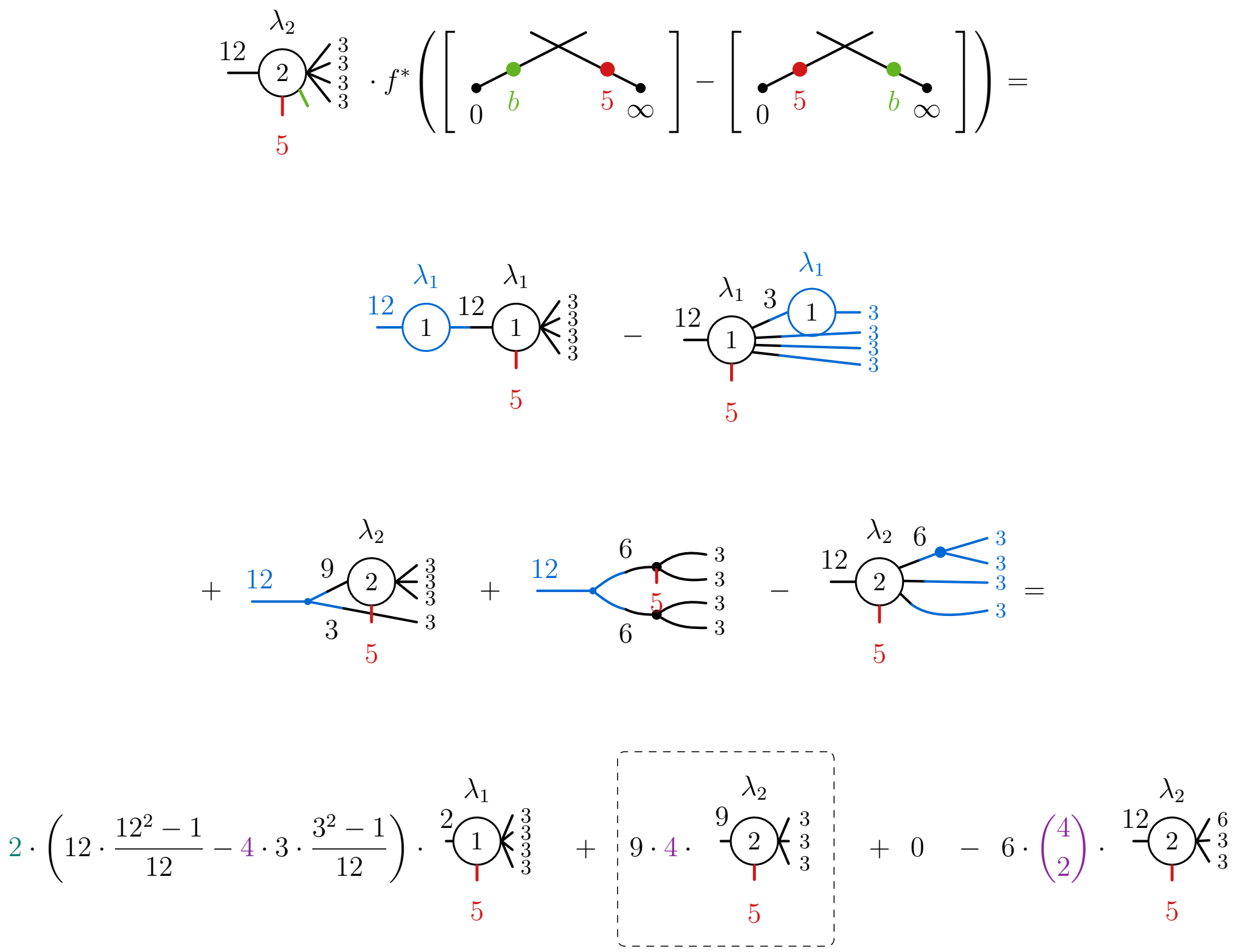}}}
    \caption{The WDVV relation from Example \ref{ex:wdvv2}. The boxed graph corresponds to an invariant for covers of strictly lower degree than $12$. We omit the flag corresponding to the point $b$ to avoid cluttering the graphs: it is understood that it is always on a vertex on the opposite side of the red $5$ flag.}
    \label{fig:wdvv2}
\end{figure}

\begin{ex}\label{ex:wdvv2}
    We consider the WDVV relation on the moduli space $\adm_{2}((12),(5)_{12},(3^4))$ with $\alpha=\lambda_2$, depicted in Figure \ref{fig:wdvv2}. In the second line we have the two divisor graphs with a genus $1$ $\Trz{1}$-graph, in the third line the divisor graphs with a genus $0$ $\Trz{1}$-graph. The fourth line contains the various multiplicities to obtain the relation
    \begin{equation*}
        135x_{4;1}(3^4) + 18x_{4;2}(3^3) - 18x_{4;2}(6,3^2)=0.
    \end{equation*}
    The middle term, coming from a disconnected graph with a trivial component, corresponds to an invariant for covers of degree strictly smaller than $12$.
\end{ex}

\subsection{The quantum $\star$-product and commutation} \label{subsec:star}
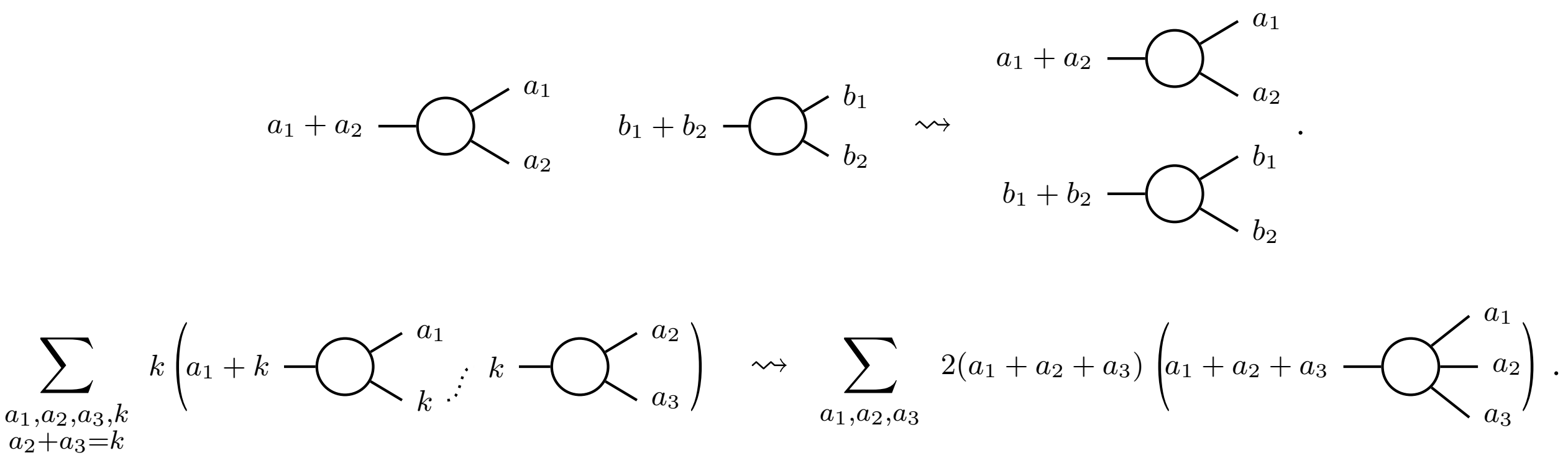
\begin{figure}
    \centering
 \begin{equation*}
    \begin{tikzpicture}[baseline=-0.5ex, thick, scale=0.6, font=\scriptsize, vertex/.style={draw, circle, minimum size=0.6cm, inner sep=0pt}]
        \node[vertex] (v) at (0,0) {};
        \draw (v.west) -- ++(-0.67, 0) node[left] {$a_1+a_2$};
        \draw (v.30) -- ++(0.67, 0.4) node[right] {$a_1$};
        \draw (v.-30) -- ++(0.67, -0.4) node[right] {$a_2$};
    \end{tikzpicture}
    \hspace{1em}
    \begin{tikzpicture}[baseline=-0.5ex, thick, scale=0.4, font=\scriptsize, vertex/.style={draw, circle, minimum size=0.6cm, inner sep=0pt}]
        \node[vertex] (v) at (0,0) {};
        \draw (v.west) -- ++(-0.67, 0) node[left] {$b_1+b_2$};
        \draw (v.30) -- ++(0.67, 0.4) node[right] {$b_1$};
        \draw (v.-30) -- ++(0.67, -0.4) node[right] {$b_2$};
    \end{tikzpicture}
    \hspace{0.5em} \leadsto \hspace{0.5em}
    \begin{tikzpicture}[baseline=-0.5ex, thick, scale=0.6, font=\scriptsize, vertex/.style={draw, circle, minimum size=0.6cm, inner sep=0pt}]
        \node[vertex] (vTop) at (0, 1.2) {};
        \draw (vTop.west) -- ++(-0.67, 0) node[left] {$a_1+a_2$};
        \draw (vTop.30) -- ++(0.67, 0.4) node[right] {$a_1$};
        \draw (vTop.-30) -- ++(0.67, -0.4) node[right] {$a_2$};
        \node[vertex] (vBot) at (0, -1.2) {};
        \draw (vBot.west) -- ++(-0.67, 0) node[left] {$b_1+b_2$};
        \draw (vBot.30) -- ++(0.67, 0.4) node[right] {$b_1$};
        \draw (vBot.-30) -- ++(0.67, -0.4) node[right] {$b_2$};
    \end{tikzpicture}.
\end{equation*}

\begin{equation*}
    \sum_{\substack{a_1,a_2,a_3,k\\a_2+a_3=k}} 
    \begin{tikzpicture}[
        baseline=-0.5ex, thick, scale=0.5,
        font=\scriptsize,
        vertex/.style={draw, circle, minimum size=0.6cm, inner sep=0pt}
    ]
        \node[vertex] (v1) at (0,0) {};
        \draw (v1.west) -- ++(-0.67, 0)
            node[left] {$k \ \ a_1+k$};
        \draw (v1.30) -- ++(0.67, 0.4)
            node[right] {$a_1$};
        \draw (v1.-30) -- ++(0.67, -0.4)
            node[right] (k1) {$k$};

        \node[vertex] (v2) at (5,0) {}; 
        \draw (v2.west) -- ++(-0.67, 0)
            node[left] (k2) {$k$};
        \draw (v2.30) -- ++(0.67, 0.4)
            node[right] {$a_2$};
        \draw (v2.-30) -- ++(0.67, -0.4)
            node[right] {$a_3$};

        \draw[dotted, thick]
            (k1) to[out=0, in=180] (k2);

        \node at (-3.5,0) {\scalebox{1.5}[3]{$($}};
        \node at (7.5,0) {\scalebox{1.5}[3]{$)$}};

    \end{tikzpicture}
    \hspace{0.5em} \leadsto \hspace{0.5em}
    \sum_{a_1,a_2,a_3}
    \begin{tikzpicture}[
        baseline=-0.5ex, thick, scale=0.6,
        font=\scriptsize,
        vertex/.style={draw, circle, minimum size=0.6cm, inner sep=0pt}
    ]
        \node[vertex] (v) at (0,0) {};
        \draw (v.west) -- ++(-0.67, 0)
            node[left] {$2(a_1+a_2+a_3) \ \  a_1+a_2+a_3$};
        \draw (v.45) -- ++(0.67, 0.53)
            node[right] {$a_1$};
        \draw (v.0) -- ++(0.67, 0)
            node[right] {$a_2$};
        \draw (v.-45) -- ++(0.67, -0.53)
            node[right] {$a_3$};

        \node at (-4.4,0) {\scalebox{1.5}[3]{$($}};
        \node at (2.1,0) {\scalebox{1.5}[3]{$)$}};

    \end{tikzpicture}.
\end{equation*}

    \caption{A schematic illustration of how the $\star$ product models the gluing of graphs: the top line corresponds to the $n=0$ term, the second line to $n=1$ in \eqref{eq:star_prod}. }
    \label{fig:starprod}
\end{figure}

To model the gluing of smooth curves to form nodal admissible covers of a nodal rational curve, we use the $\star$-product as in \cite{Intro2SFT} and \cite{BR2016q}. For any ring $R$, the $\star$-product on the $R$-algebra $R[\h][q_*][[p_*]]$ is given by
\begin{equation}\label{eq:star_prod}
    f \star g = f\exp\left(\sum_{k>0}\h k\overleftarrow{\frac{\partial}{\partial p_k}}\,\overrightarrow{\frac{\partial}{\partial q_k}}\right)g = \sum_{n\ge0}\frac{\h^n}{n!}\sum_{k_1,\dots,k_n>0}\frac{\partial^nf}{\partial p_1\cdots\partial p_n}\left(\prod_{i=1}^nk_i\right)\frac{\partial^ng}{\partial q_1\cdots\partial q_n}.
\end{equation}
Writing $f$ and $g$ in normal form, i.e. with all the $p_*$ variables to the left and $q_*$ variables to the right, the product $f\star g$ coincides with commuting all the $q_*$ variables of $f$ with the $p_*$ variables of $g$ according to the commutation rule
\begin{equation*}
    \left[q_j,p_k\right] = j\delta_{jk}.
\end{equation*}
To match the algebra in which the generating series $\Tr{\mu,d}^\bullet$ reside, we take $R=\mathbb{C}[\h^{-1},\e,v]$. The $\star$-product models the gluing of curves to form possibly disconnected nodal admissible covers. For example, consider
\begin{equation} \label{ex:starprod}
    f\star g = \left(\frac{1}{\h}\sum_{a_1,a_2\in\mathbb{N}}\frac{p_{a_1}p_{a_2}}{2!}q_{a_1+a_2}\right)\star\left(\frac{1}{\h}\sum_{b_1,b_2\in\mathbb{N}}\frac{p_{b_1}p_{b_2}}{2!}q_{b_1+b_2}\right).
\end{equation}
Both series in this product encode the invariants $x_{1;0}(a_1, a_2)$.
Since $g$ is linear in the variables $q_*$, only the $n=0$ and $n=1$ terms in \eqref{eq:star_prod} contribute to the $\star$-product \eqref{ex:starprod}. The $n=0$ term is  normal multiplication 
\begin{equation*}
    \frac{1}{\h^2}\sum_{a_1,a_2\in\mathbb{N}}\sum_{b_1,b_2\in\mathbb{N}}\frac{p_{a_1}p_{a_2}}{2!}\frac{p_{b_1}p_{b_2}}{2!}q_{a_1+a_2}q_{b_1+b_2},
\end{equation*}
and it corresponds to when no gluing takes place, giving rise to a disconnected cover with two components. The $n=1$ term
\begin{align*}
    & \;\h\sum_{k>0}k\frac{\partial}{\partial p_k}\left(\frac{1}{\h}\sum_{a_1,a_2\in\mathbb{N}}\frac{p_{a_1}p_{a_2}}{2!}q_{a_1+a_2}\right)\frac{\partial}{\partial q_k}\left(\frac{1}{\h}\sum_{b_1,b_2\in\mathbb{N}}\frac{p_{b_1}p_{b_2}}{2!}q_{b_1+b_2}\right) \\
    = & \; \frac{1}{\h}\sum_{k>0}\left(\sum_{a_1\in\mathbb{N}}p_{a_1}q_{a_1+k}\right)k\Biggl(\sum_{\substack{a_2,a_3\in\mathbb{N}\\a_2+a_3=k}}\frac{p_{a_2}p_{a_3}}{2!}\Biggr) = \frac{1}{\h}\sum_{a_1,a_2,a_3\in\mathbb{N}}2(a_1+a_2+a_3)\frac{p_{a_1}p_{a_2}p_{a_3}}{3!}q_{a_1+a_2+a_3}
\end{align*}
corresponds to when a point of ramification $k$ over $\infty$ on the first curve is glued to a point of ramification $k$ over $0$ on the second curve. See Figure \ref{fig:starprod} for an illustration.

We use the $\star$ product to encode WDVV relations. 

\begin{lem}\label{lem:TkWDVV}
    The commutation relation
    \begin{equation*}
        [\Trz{k},\Trz{1}] := \Trz{k}\star\Trz{1} - \Trz{1}\star\Trz{k} = 0
    \end{equation*}
    corresponds to the WDVV relations  involving the invariants $x_{k;g}(\ain)$ using the class $\alpha = \lambda_g$. In particular, for every $0\le g\le k$ and $\ain=(a_1,\dots,a_{k+2-g})$, the coefficient of 
        $\hbar^{g-1}\e^{2g}\frac{v^g}{g!}p_{\ain}q_{\lvert \ain\rvert}$
    in $[\Trz{k},\Trz{1}]$ corresponds to the WDVV relation \eqref{eq:WDVVcyc} on the space $\adm_g(\delta,(k+1),\ain)$.
\end{lem}
\begin{proof}
   The $n=0$ terms of the two $\star$ products \eqref{eq:star_prod} in the commutator cancel.  Since $\Trz{1},\Trz{k}$ are linear in the $q_*$ variables, only the $n=1$ terms survive: 
    \begin{equation*}
        [\Trz{k},\Trz{1}] = \h\sum_{\kappa>0}\kappa\left(\frac{\partial\Trz{k}}{\partial p_\kappa}\frac{\partial\Trz{1}}{\partial q_\kappa} - \frac{\partial\Trz{1}}{\partial p_\kappa}\frac{\partial\Trz{k}}{\partial q_\kappa}\right).
    \end{equation*}
    Let us write $\Trz{1}=\Trz{1}^{(0)}+\Trz{1}^{(1)}$, where
    \begin{equation*}
        \Trz{1}^{(0)} = \frac{1}{\h}\sum_{\ain\in\mathbb{N}^2}\frac{p_\ain}{2!}q_{\lvert\ain\rvert}, \quad \Trz{1}^{(1)}=\e^2v\sum_{a\in\mathbb{N}}\frac{a^2-1}{12}p_aq_a.
    \end{equation*}
    We have
    \begin{align*}
        & \h\sum_{\kappa>0}\kappa\frac{\partial\Trz{k}}{\partial p_\kappa}\frac{\partial\Trz{1}^{(0)}}{\partial q_\kappa} = \h\sum_{g=0}^k\h^{g-1}\e^{2g}\frac{v^g}{g!}\sum_{\kappa>0}\kappa \\
        & \cdot\Biggl(\sum_{a_3,\dots,a_{k+2-g}\in\mathbb{N}}x_{k;g}(\kappa,a_3,\dots,a_{k+2-g})\frac{p_{a_3}\cdots p_{a_{k+2-g}}}{(k-g)!}q_{\kappa+a_3+\cdots+a_{k+2-g}}\Biggr)\Biggl(\frac{1}{\h}\sum_{\substack{a_1,a_2\in\mathbb{N}\\a_1+a_2=\kappa}}\frac{p_{a_1}p_{a_2}}{2!}\Biggr) \\
        = & \sum_{g=0}^k\h^{g-1}\e^{2g}\frac{v^g}{g!}\frac12\sum_{\ain\in\mathbb{N}^{k+2-g}}(a_1+a_2)x_{k;g}(a_1+a_2,a_3,\dots,a_{k+2-g})\frac{p_{\ain}}{(k-g)!}q_{\lvert\ain\rvert} \\
        = & \sum_{g=0}^k\h^{g-1}\e^{2g}\frac{v^g}{g!}\sum_{\ain\in\mathbb{N}^{k+2-g}}\sum_{i<j}(a_i+a_j)x_{k;g}(\hat{a}_i,\hat{a}_j,a_1,\dots,a_{k+2-g},a_i+a_j)\frac{p_{\ain}}{(k+2-g)!}q_{\lvert\ain\rvert},
    \end{align*}
    where the hats denote omission. Next,
    \begin{align*}
        & \h\sum_{\kappa>0}\kappa\frac{\partial\Trz{k}}{\partial p_\kappa}\frac{\partial\Trz{1}^{(1)}}{\partial q_\kappa} = \h\sum_{g=0}^k\h^{g-1}\e^{2g}\frac{v^g}{g!}\sum_{\kappa>0}\kappa \\
        & \cdot\Biggl(\sum_{a_2,\dots,a_{k+1-g}\in\mathbb{N}}x_{k;g}(\kappa,a_2,\dots,a_{k+1-g})\frac{p_{a_2}\cdots p_{a_{k+1-g}}}{(k+1-g)!}q_{\kappa+a_2+\cdots+a_{k+1-g}}\Biggr)\e^2v\frac{\kappa^2-1}{12}p_\kappa \\
        = & \sum_{g=0}^k\h^g\e^{2g+2}\frac{v^{g+1}}{g!}\sum_{\ain\in\mathbb{N}^{k+1-g}}a_1\frac{a_1^2-1}{12}x_{k;g}(a_1,\dots,a_{k+1-g})\frac{p_{\ain}}{(k+1-g)!}q_{\lvert\ain\rvert} \\
        = & \sum_{g=1}^{k+1}\h^{g-1}\e^{2g}\frac{v^g}{g!}\sum_{\ain\in\mathbb{N}^{k+2-g}}g\sum_{i=1}^{k+2-g}a_i\frac{a_i^2-1}{12}x_{k;g-1}(a_1,\dots,a_{k+2-g})\frac{p_{\ain}}{(k+2-g)!}q_{\lvert\ain\rvert}.
    \end{align*}
    By a similar procedure,
    \begin{align*}
        & \h\sum_{\kappa>0}\kappa\frac{\partial\Trz{1}^{(0)}}{\partial p_\kappa}\frac{\partial\Trz{k}}{\partial q_\kappa} \\
        = & \; \sum_{g=0}^k\h^{g-1}\e^{2g}\frac{v^g}{g!}\sum_{\ain\in\mathbb{N}^{k+2-g}}\sum_{i=1}^{k+2-g}(\lvert\ain\rvert-a_i)x_{k;g}(\hat{a}_i,a_1,\dots,a_{k+2-g})\frac{p_\ain}{(k+2-g)!}q_{\lvert\ain\rvert},
    \end{align*}
    and 
    \begin{align*}
        & \h\sum_{\kappa>0}\kappa\frac{\partial\Trz{1}^{(1)}}{\partial p_\kappa}\frac{\partial\Trz{k}}{\partial q_\kappa} \\
        = & \; \sum_{g=1}^{k+1}\h^{g-1}\e^{2g}\frac{v^g}{g!}\sum_{\ain\in\mathbb{N}^{k+2-g}}g\,\lvert\ain\rvert\frac{\lvert\ain\rvert^2-1}{12}x_{k;g-1}(a_1\dots,a_{k+2-g})\frac{p_{\ain}}{(k+2-g)!}q_{\lvert\ain\rvert}.
    \end{align*}
    Therefore, for every $0\le g\le k$ and $\ain\in\mathbb{N}^{k+2-g}$, the coefficient of $\h^{g-1}\e^{2g}\frac{v^g}{g!}p_\ain q_{\lvert\ain\rvert}$ in the commutator $[\Trz{k},\Trz{1}]$ is 
    \begin{align*}
        & \sum_{i<j}x_{k;g}(\hat{a}_i,\hat{a}_j,a_1,\dots,a_{k+2-g},a_i+a_j)\cdot(a_i+a_j)\cdot x_{1;0}(a_i,a_j) \\
        & + g\sum_{i=1}^{k+2-g}x_{k,g-1}(a_1,\dots,a_{k+2-g})\cdot a_i\cdot x_{1;1}(a_i) \\
        & - \sum_{i=1}^{k+2-g}x_{1;0}(\lvert\ain\rvert-a_i,a_i)\cdot(\lvert\ain\rvert-a_i)\cdot x_{k;g}(\hat{a}_i,a_1,\dots,a_{k+2-g}) \\
        & - g\,x_{1;1}(\lvert\ain\rvert)\cdot \lvert\ain\rvert\cdot x_{k;g-1}(a_1,\dots,a_{k+2-g}) \,
    \end{align*}
   which one can recognize as \eqref{eq:WDVVcyc}, thus completing the proof.
\end{proof}

\begin{remark}\label{rem:Tk1k2_comm}
    It is a matter of careful book-keeping  to see that all commutators of the generating functions introduced correspond to WDVV relations on appropriate spaces of admissible covers (with four special branch points and possibly $\psi$ classes on two of them).  Hence:
    \begin{equation*}
        [\Tr{\mu_1,d_1}^\bullet,\Tr{\mu_2,d_2}^\bullet] = \Tr{\mu_1,d_1}^\bullet\star\Tr{\mu_2,d_2}^\bullet - \Tr{\mu_2,d_2}^\bullet\star\Tr{\mu_1,d_1}^\bullet=0.
    \end{equation*}
\end{remark}

\section{Reconstruction}

\subsection{Reconstructing the $\Trz{k}$}\label{sec:recTk}
This section is devoted to proving Proposition \ref{prop:Tk}. 
We start with some preliminary lemmas. First we derive the needed evaluation of $x_{k,0}(1^{k+1})$.  
\begin{lem}\label{lem:usingstring}
   For all $k \geq 0$ we have
   \[
   x_{k;0}(1^{k+1})= k!\,.
   \]
\end{lem}
\begin{proof} There is a unique isomorphism class of covers which are fully ramified over two base points and unramified over a third marked point, with a cyclic automorphism group of order $k+1$. But we are requiring the additional data of a marking of all inverse images of the unramified point, which contributes a factor of $(k+1)!$. All together we have $x_{k;0}(1^{k+1})= (k+1)!/(k+1) = k!$.
\end{proof}

\begin{remark}
    This evaluation is compatible with the unstable initial condition \eqref{eq:incond} and the string equation (see Lemma \ref{lem:string}).
    In fact, the unstable invariant is defined (for convenience) by requiring compatibility with the string equation.
\end{remark}
\begin{proof}[Proof of Proposition \ref{prop:Tk}]
    Fixing the value of $k$, the proof is by induction on the degree $\delta \geq k+1$. We describe carefully the base case $\delta = k+1$ as it contains already most of the features that are needed for the inductive step. 
    In degree $k+1$, the invariants $x_{k;g}(\nu)$ we need to compute, in decreasing genus, are  \[x_{k;k}(k+1), x_{k;k-1}(k,1), \ldots, x_{k;k-1}(\lceil{k/2}\rceil,\lfloor{k/2}\rfloor), \ldots, x_{k;k-j}(\nu_{j,i}), \ldots, x_{k;1}(2, 1^{k-1});\]
    for fixed genus $k-j$, as the index $i$ varies the $\nu_{j,i}$'s list all the partitions of $k+1$ of length $j+1$. We do not include the genus zero invariant $x_{k;0}(1^{k+1})$ as that is a known initial condition.
    
    We denote by $\rho_{k;k-j}(\nu)$ the $WDVV$ relation applied to the $1$ dimensional  class $\lambda_{k-j}$ in the space $\adm_{k-j}(k+1, k+1, \nu)$, where now  $\ell(\nu) = j+2$. There are no genus zero ($j=k$) relations, because in that case the length of the partition exceeds the degree of the cover.
    There is a natural bijection $\beta$ between invariants and relations by setting
    $$\beta(x_{k;g}(\nu))=\rho_{k;g+1}(\nu).$$ Unfortunately, the relation $\rho_{k;k+1}(k+1)$ is trivial as the corresponding graph is symmetric; but then again the genus zero invariant is already known so we still have the same number of equations as unknowns. We must show that the matrix of the corresponding linear system is invertible.

To achieve that, we observe the relations  more closely.
As described in Lemma \ref{lem:graphic_wdvv}, for a given relation $\rho_{k;g}(\nu)$, the terms with non-zero coefficients are described graphically as the remaining graph after one pulls out either a right-pointing genus zero tripod, or a two-valent genus one vertex carrying $\lambda_1$; since $\delta = k+1$, trying to pull out a tripod from the left produces a graph not of compact type, which is then killed by the  class $\lambda_g$. Tripods may only be pulled out from the right, corresponding to invariants $x_{k;g}(\nu')$, with $\nu'$ a partition of length one less than $\nu$. Pulling out a tripod from the right corresponds to being in the second term of the commutator, and therefore these invariants appear in the relation with negative coefficients.

A genus one vertex can instead be pulled out along any edge of the graph of $\rho_{k;g}(\nu)$, contributing a multiple of the invariant $x_{k;g-1}(\nu)$. Adding the signed coefficients obtained by pulling out the genus one vertex in all possible ways, one obtains that the coefficient 
of $x_{k;g-1}(\nu)$ is 
\[
\frac{g}{12}\cdot\left((k+1)^3 - \sum_{i=1}^\ell \nu_i^3\right),
\]
which is a positive number.
Now we can construct the (augmented) matrix of the linear system of relations, see Figure \ref{fig:T4} for an illustration in the case $k=4$. The matrix has a block decomposition by grouping rows (resp. columns) based on the genus of the relation (resp. invariant). There are square blocks corresponding to groups matched by the bijection $\beta$, which are diagonal matrices with positive coefficients. The only non-zero other coefficients (other than the last entry in the augmentation vector) appear in the (typically rectangular) blocks immediately to the left of the square blocks, and all the non-zero entries here are negative.
In order to show that the matrix of the linear system is invertible, we argue that its columns are linearly independent. Assume there is a linear combination of columns giving the zero vector. Since only the first two columns have non-zero first entries (with opposite signs), it follows that the first two coefficients of the linear combination are either both zero, or have the same sign. Arguing similarly for each successive column, one obtains that either all coefficients in the linear combination are zero, or they are all non-zero and in particular with the same sign. The latter case yields a contradiction as the resulting linear combination has a non-zero final entry, thus proving linear independence.
We have therefore shown that all degree $k+1$ invariants are reconstructed from the $WDVV$ relations from the initial condition $x_{k,0}(1^{k+1}) = k!$.

\begin{figure}
    \centering
 \resizebox{0.8 \textwidth}{!}{
\newcommand{\profile}[4]{%
  \begin{tikzpicture}[
    baseline=(current bounding box.center),
    x=1cm,y=1cm,
    line cap=round,
    every node/.style={font=\scriptsize}
  ]
    \path[use as bounding box] (-1.18,-0.64) rectangle (1.18,0.70);
    \node[circle,draw,thick,minimum size=6.2mm,inner sep=0pt] (v) at (0,0)
      {$#2$};
    \node[anchor=south,inner sep=0pt,font=\scriptsize] at (0,0.32)
      {$\lambda_{#1}$};
    \coordinate (leftend) at (-0.78,0);
    \draw[thick] (leftend) -- (-0.31,0);
    \node[anchor=east,inner sep=0pt] at (-0.79,0) {$5$};
    \foreach \w [count=\i from 1] in {#4}{
      \pgfmathsetmacro{\ya}{0.10*((#3+1)/2-\i)}
      \pgfmathsetmacro{\xa}{sqrt(0.31*0.31-\ya*\ya)}
      \pgfmathsetmacro{\ye}{0.18*((#3+1)/2-\i)}
      \draw[thick] (\xa,\ya) -- (0.78,\ye);
      \node[anchor=west,inner sep=1pt] at (0.78,\ye) {$\w$};
    }
  \end{tikzpicture}%
}

\newcommand{\limitprofile}{%
  \begin{tikzpicture}[
    baseline=(current bounding box.center),
    x=1cm,y=1cm,
    line cap=round,
    every node/.style={font=\scriptsize}
  ]
    \path[use as bounding box] (-1.18,-0.64) rectangle (1.18,0.70);
    \coordinate (v) at (0,0);
    \fill (v) circle (1.6pt);
    \coordinate (leftend) at (-0.72,0);
    \draw[thick] (leftend) -- (v);
    \node[anchor=east,inner sep=0pt] at (-0.74,0) {$5$};
    \foreach \yy in {-0.36,-0.18,0,0.18,0.36}
      \draw[thick] (v) .. controls (0.28,0.55*\yy) and (0.46,\yy) .. (0.70,\yy);
    \draw[densely dashed] (0.78,-0.47) -- (0.78,0.47);
    \draw[thick] (-1.02,0.61)
      .. controls (-1.20,0.38) and (-1.20,-0.38) .. (-1.02,-0.61);
    \draw[thick] (1.02,0.61)
      .. controls (1.20,0.38) and (1.20,-0.38) .. (1.02,-0.61);
  \end{tikzpicture}%
}

\begin{tikzpicture}[
  every node/.style={font=\normalsize},
  sign/.style={minimum width=11mm,minimum height=11mm,inner sep=0pt}
]

\matrix (m) [
  matrix of nodes,
  nodes={sign},
  column sep=6mm,
  row sep=5mm
]{
  & \profile{4}{4}{1}{5}
  & \profile{3}{3}{2}{4,1}
  & \profile{3}{3}{2}{3,2}
  & \profile{2}{2}{3}{3,1,1}
  & \profile{2}{2}{3}{2,2,1}
  & \profile{1}{1}{4}{2,1,1,1}
  & \limitprofile \\

  \profile{4}{4}{2}{4,1}
  & $-$ & $+$ & $0$ & $0$ & $0$ & $0$ & $0$ \\

  \profile{4}{4}{2}{3,2}
  & $-$ & $0$ & $+$ & $0$ & $0$ & $0$ & $0$ \\

  \profile{3}{3}{3}{3,1,1}
  & $0$ & $-$ & $-$ & $+$ & $0$ & $0$ & $0$ \\

  \profile{3}{3}{3}{2,2,1}
  & $0$ & $-$ & $-$ & $0$ & $+$ & $0$ & $0$ \\

  \profile{2}{2}{4}{2,1,1,1}
  & $0$ & $0$ & $0$ & $-$ & $-$ & $+$ & $0$ \\

  \profile{1}{1}{5}{1,1,1,1,1}
  & $0$ & $0$ & $0$ & $0$ & $0$ & $-$ & $\ast$ \\
};

\coordinate (gridNW) at ($(m-2-2.north west)+(-2mm,2mm)$);
\coordinate (gridNE) at ($(m-2-8.north east)+(2mm,2mm)$);
\coordinate (gridSW) at ($(m-7-2.south west)+(-2mm,-2mm)$);
\coordinate (gridSE) at ($(m-7-8.south east)+(2mm,-2mm)$);
\coordinate (tableNW) at ($(gridNW)+(-5mm,0)$);
\coordinate (tableNE) at ($(gridNE)+(5mm,0)$);
\coordinate (tableSW) at ($(gridSW)+(-5mm,0)$);
\coordinate (tableSE) at ($(gridSE)+(5mm,0)$);

\draw[thick]
  ($(tableNW)+(3mm,0)$) -- (tableNW)
  -- (tableNW |- tableSW) -- ($(tableNW |- tableSW)+(3mm,0)$);
\draw[thick]
  ($(tableNE)+(-3mm,0)$) -- (tableNE)
  -- (tableNE |- tableSE) -- ($(tableNE |- tableSE)+(-3mm,0)$);

\coordinate (xA) at ($(m-2-2.east)!0.5!(m-2-3.west)$);
\coordinate (xB) at ($(m-2-4.east)!0.5!(m-2-5.west)$);
\coordinate (xC) at ($(m-2-6.east)!0.5!(m-2-7.west)$);
\coordinate (xD) at ($(m-2-7.east)!0.5!(m-2-8.west)$);
\coordinate (yA) at ($(m-3-2.south)!0.5!(m-4-2.north)$);
\coordinate (yB) at ($(m-5-2.south)!0.5!(m-6-2.north)$);
\coordinate (yC) at ($(m-6-2.south)!0.5!(m-7-2.north)$);

\foreach \x in {xA,xB,xC}
  \draw[thick] (\x |- tableNW) -- (\x |- tableSW);

\fill[black!22]
  ($(xD |- tableNW)+(-1.1mm,0)$) rectangle
  ($(xD |- tableSW)+(1.1mm,0)$);
\draw[thin] ($(xD |- tableNW)+(-1.1mm,0)$) -- ($(xD |- tableSW)+(-1.1mm,0)$);
\draw[thin] ($(xD |- tableNW)+( 1.1mm,0)$) -- ($(xD |- tableSW)+( 1.1mm,0)$);

\foreach \y in {yA,yB,yC}{
  \draw[thick] (gridNW |- \y) -- ($(xD |- \y)+(-1.1mm,0)$);
  \draw[thick] ($(xD |- \y)+(1.1mm,0)$) -- (gridNE |- \y);
}

\end{tikzpicture}} 
    \caption{The signs of the coefficients of the augmented matrix of the linear system reconstructing the degree $k+1 = 5$ invariants of $\Trz{4}$. The columns correspond to the invariants and the rows to the relations, both depicted graphically. The matrix is ``one column away" from being lower triangular.}
    \label{fig:T4}
\end{figure}
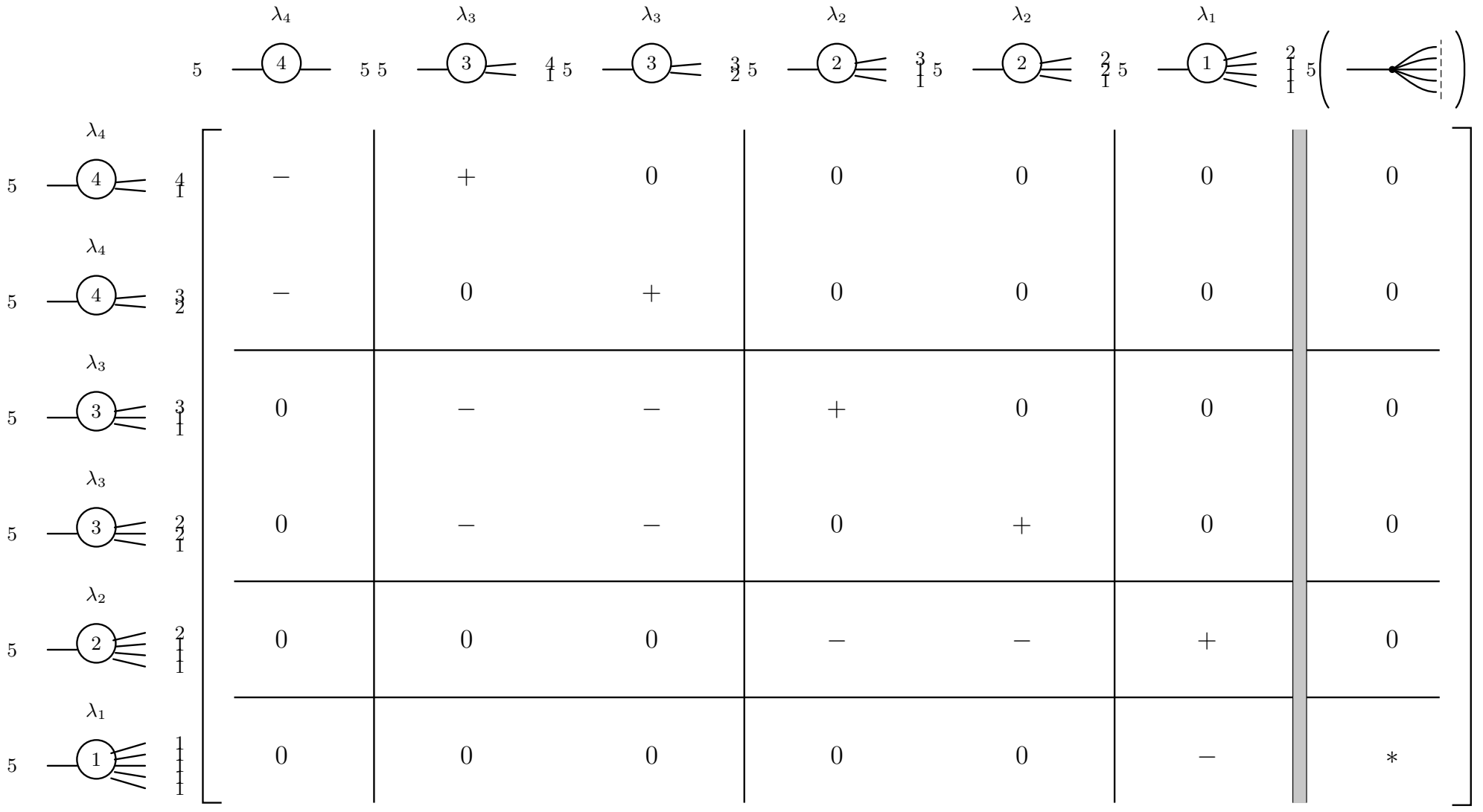
For the inductive step, assume $\delta>k+1$. One may still order invariants and relations in decreasing order of genus, and define an injective map $\beta$ assigning to each invariant the relation obtained by increasing $g$ by one. Just as before the genus $k+1$ relation $\rho_{k;k+1}(\delta)$ is trivial. The significant differences with respect to the base case are the following:
\begin{enumerate}
    \item the genus zero invariants are no longer initial conditions, but unknowns;
    \item with $\delta\geq k+2$, there are now genus zero relations, which are not in the image of the function $\beta$;
    \item since the ramification profile over the second branch point is $(k+1)$, it is now possible to pull out a tripod from the left, with a disconnected, compact type graph on the right. In that case, there are two possibilities: if there are at least two nontrivial connected components on the disconnected graph then we are in situation $(1)$ from Lemma \ref{lem:disceval}, and that contribution vanishes. If the only nontrivial component is the one hosting the point of ramification order $k+1$, then there are some trivial components and we are in situation $(2)$ of Lemma \ref{lem:disceval}. The degree of the cover for the non-trivial component is strictly less than $\delta$, and therefore by induction the invariant is known.
\end{enumerate}
We can then add exactly one genus zero relation to all the (non-trivial) relations in the image of $\beta$ and obtain an augmented matrix for the linear system which has a chance at a unique solution. By $(3)$ above all coefficients coming from different behavior than the base case are known and end up contributing to the augmentation vector. The matrix of principal parts has a similar block shape as in the base case and the argument for the linear independence of its columns carries over verbatim. This concludes the proof of Proposition \ref{prop:Tk}.
\end{proof}
This result can be used to effectively compute the invariants $x_{k;g}(\ain)$;  in Table \ref{tab:xkg} we collect  the coefficients for the generating series $\Trz{2},\Trz{3},\Trz{4}$ explicitly. In Appendix \ref{app:python}, we show how to use \texttt{pyBoLaNO} \cite{pybolano}, a \texttt{Python} symbolic package that is useful in implementing the $\star$-product \eqref{eq:star_prod}, to do these computations. It is worth noting that to perform these computations we use that $x_{k;g}(a_1,\dots,a_{k+1-g})$ is a symmetric polynomial in the $a_i$ of degree $2g$, which we prove in Corollary \ref{cor:polynomiality}.

\begin{table}[bt]
    \centering
    \renewcommand{\arraystretch}{1.8}
    \setlength{\tabcolsep}{6pt}
    \resizebox{\textwidth}{!}{%
    \begin{tabular}{c|ccccc}
        {\large $x_{k;g}(\ain)$} & $g=0$ & $1$ & $2$ & $3$ & $4$ \\
        \hline

        $k=2$
        &
        $2$
        &
        $\displaystyle \frac{h_2-3}{6}$
        &
        $\displaystyle \frac{(h_2-1)(h_2-4)}{60}$
        &
        &
        \\[2mm]

        $3$
        &
        $6$
        &
        $\displaystyle \frac{h_2-6}{2}$
        &
        $\displaystyle \frac{(h_2-3)(h_2-7)}{20}$
        &
        $\displaystyle \frac{3(h_2-1)(h_2-4)(h_2-9)}{560}$
        &
        \\[2mm]

        $4$
        &
        $24$
        &
        $\displaystyle 2(h_2-10)$
        &
        $\displaystyle
        \frac{3h_2^2-e_1e_3-50h_2+195}{15}$
        &
        $\displaystyle
        \frac{3h_2^3+(e_1e_2)^2-70h_2^2+455h_2-820}{140}$
        &
        $\displaystyle
        \frac{(h_2-1)(h_2-4)(h_2-9)(h_2-16)}{420}$
    \end{tabular}%
    }\vspace{.5cm}
    
    \caption{
        The coefficients $x_{k;g}$ of  $\Trz{k}$ for $2\leq k\leq 4$.
        In each entry,  $\ell(\ain) = k+1-g$. We write
        $h_2$
        for the complete homogeneous symmetric polynomial of degree two,$h_2
            = \sum_{i = 1}^{k+1-g} a_i^2+\sum_{i<j}a_i a_j$,
        and $e_i$
        for the $i$-th elementary symmetric polynomial in the same variables.
    }
    \label{tab:xkg}
\end{table}

\subsection{Reconstructing the $\Tr{\mu}^\bullet$} 

We turn to the more general invariants with no $\psi$-class, namely the $x_{\mu;g}(\ain)$, and show they can be expressed in terms of the building blocks $x_{k;g}(\ain)$. 
\begin{prop}\label{prop:TmuTkinv}
    Let $\mu$ be a partition. Any invariant of type $x_{\mu;g}(\ain)$ may be computed in terms of  invariants $x_{k,\tilde{g}}(\tilde{\ain})$, with $k\leq |\mu|$, $\tilde{g} \leq g$ and $|\tilde{\ain}|\leq |\ain|$.
\end{prop}
\begin{proof} Consider the data $(\delta = |\ain|, g, \mu)$
    We prove the result by induction on the following partial order on this data:
    \begin{align} \label{eq:pord}
    (\delta_1, g_1, \mu_1)>  (\delta_2, g_2, \mu_2)  \iff &\delta_1 >\delta_2\ \ \ \mbox{or} \nonumber \\
    & \delta_1 =\delta_2 \ \ \ \mbox{and}\ \ \  g_1>g_2 \ \ \ \mbox{or}\nonumber\\
    & \delta_1 =\delta_2 \ \ \ \mbox{and}\ \ \  g_1=g_2 \ \ \ \mbox{and}\ \ \ \ell(\mu_1)>\ell(\mu_2).
    \end{align}
     For any $(g,\ain)$, the statement is trivially true if  $\ell(\mu)=1$. Let $\mu$ be a partition of length at least 2, and consider a (possibly nonvanishing) invariant $x_{\mu;g}(\ain)$. Choose nonempty partitions $\mu_1$ and $\mu_2$ such that $\mu_1\sqcup\mu_2=\mu$. Consider the moduli space 
     \begin{equation*}
         \adm_g(\delta, \mu_1+1, \mu_2+1, \ain),
     \end{equation*}
of dimension $g+1$. We analyze the WDVV relation (using notation of \eqref{eq:WDVVcyc}):
\begin{equation} \label{eq:wdvvtmu}
      \lambda_g\cdot f^\ast D(\{1,4\}|\{2,3\})= \lambda_g \cdot f^\ast D(\{1,2\}|\{3,4\}).  
\end{equation}
There is one term on the left-hand side of \eqref{eq:wdvvtmu} corresponding to covers that, over the component containing the second and third special branch points,  consist of a disjoint union of trivial covers. For all trivial covers of degree greater than one, one of the ramification points maps to the node, the other must be a part of  either $\mu_1$ of $\mu_2$ and correspondingly maps to the appropriate branch point. This term contributes a non-zero multiple of the invariant $x_{\mu,g}(\ain)$.

It is easy to see that all other terms in \eqref{eq:wdvvtmu} give rise to products of invariants corresponding to data strictly lower than $(\delta, g, \mu)$ in \eqref{eq:pord}, thus establishing the inductive step. 
    \end{proof}

We have proved reconstruction for all connected invariants; we may use Corollary \ref{cor:expT} to obtain all disconnected ones.
This reconstruction result is efficiently expressed in terms of generating series.     
\begin{prop}\label{prop:TmuTk}
    Let $\mu=(\mu_1,\dots,\mu_\ell)$ be a partition, so that $\mu_1\ge\mu_2\ge\cdots\ge\mu_\ell$. The generating series $\Tr{\mu}^\bullet$ can be expressed in terms of $\star$-products of $\Trz{1},\dots,\Trz{\lvert\mu\rvert}$ as
    \begin{equation*}
        \Tr\mu^\bullet = \sum_{g=0}^{\lvert\mu\rvert-\mu_1}\h^g\e^{2g}\frac{v^g}{g!}\bigl(\Tr\mu^\bullet\bigr)^{(g)},
    \end{equation*}
    where $\bigl(\Tr\mu^\bullet\bigr)^{(g)}$ is a linear combination of terms $\Trz{k_1}\star\cdots\star\Trz{k_m}$, where $k_1+\dots+k_m=\lvert\mu\rvert-g$ and $m\le\max\{\ell-g,1\}$. Moreover, the leading term is
    \begin{equation*}
        (\Tr{\mu}^\bullet)^{(0)} = \frac{1}{\lvert\operatorname{Aut}(\mu)\rvert}\Trz{\mu_1}\star\Trz{\mu_2}\star\cdots\star\Trz{\mu_\ell} + \sum_{\substack{k_1,\dots,k_m\ge1\\k_1+\cdots+k_m=\lvert\mu\rvert\\m<\ell}}c_{k_1,\dots,k_m}\Trz{k_1}\star\cdots\star\Trz{k_m},
    \end{equation*}
    where $c_{k_1,\dots,k_m}$ are known constants.
\end{prop}

\begin{proof}[Proof of Proposition \ref{prop:TmuTk}]
    We rely on induction on the following partial order on  partitions:
    \begin{equation}\label{eq:partial_order}
        \nu_1 > \nu_2 \iff \lvert\nu_1\rvert>\lvert\nu_2\rvert \quad \text{or} \quad \lvert\nu_1\rvert=\lvert\nu_2\rvert \hspace{0.5em} \text{and} \hspace{0.5em} \ell(\nu_1)>\ell(\nu_2).
    \end{equation}
     The statement is true if  $\ell(\mu)=1$, since $\Tr{(k)}^\bullet=\Trz{k}$.
    Let $\mu$ be a partition of length at least 2, and choose nonempty partitions $\mu_1$ and $\mu_2$ such that $\mu_1\sqcup\mu_2=\mu$. 
Analogously to Definition \ref{def:Smud}, let $S_{\mu_1, \mu_2}$ denote the set whose elements $\sigma$ consist of
    \begin{itemize}
        \item an integer $1\le m\le\ell(\mu_1) + \ell(\mu_2)$;
        \item for every $j=1,2$ and every $i=1,\dots,m$, a \emph{possibly empty} partition $\mu_{j,i}$, such that $\coprod_{i=1}^m\mu_{j,i} = \mu_j$;
        \item for every $i=1,\dots,m$, an integer $0\le h_i\le\lvert\mu_{1,i}\rvert+\lvert\mu_{2,i}\rvert$;
        \item for every $i=1,\dots,m$, a vector of positive integers $\ain_i\in\mathbb{N}^{\lvert\mu_{1,i}\rvert+\lvert\mu_{2,i}\rvert+1-h_i}$.
    \end{itemize}
    The discrete data $\sigma = (\mu_{1,i},\mu_{2,i},h_i,\ain_i)_{i=1}^m$,  corresponds to an admissible cover space
 \begin{equation*}
  \adm_\sigma^\bullet := \dadm_{\vec{h}}\bigl(\vec{\delta},(\vec{\mu}_1+1),(\vec{\mu}_2+1),\vec{\ain}\bigr).
    \end{equation*}
    for disconnected covers of geometric genus $h = \sum h_i$; the length of the profiles $\vec{\ain}$ is chosen so the space has dimension $h+1$.

 For every $\sigma\in S_{\mu_1,\mu_2}$, we multiply the WDVV relation \eqref{eq:wdvvtmu} by the monomial   
 \begin{align}
        \label{eq:auxiliary_1}  \frac{1}{m!}\h^{h-m}\e^{2h}\frac{v^h}{h!}\binom{h}{h_1,\dots,h_m}\frac{1}{\prod_{j,i}\lvert\operatorname{Aut}(\mu_{j,i})\rvert} \prod_{i=1}^m\frac{p_{\ain_i}}{(\lvert\mu_{1,i}\rvert+\lvert\mu_{2,i}\rvert+1-h_i)!}q_{\lvert\ain_i\rvert}
    \end{align}

    Fix a geometric genus $0\le h\le\lvert\mu\rvert$, and let $S_{\mu_1,\mu_2}^{(h)}\subset S_{\mu_1,\mu_2}$ denote the subset of those $\sigma=(\mu_{1,i},\mu_{2,i},h_i,\ain_i)_{i=1}^m$ such that $\sum_{i}h_i=h$. The divisor on the right-hand side of \eqref{eq:wdvvtmu} equals
     \begin{equation*}
        f^\ast D(\{1,2\}|\{3,4\})=\sum_{\substack{g_1+g_2=h\\ |I_{g_i}|= g_i\\I_{g_1}\cap I_{g_2} =\emptyset}}br^\ast D(\{1,2\} \sqcup I_{g_1}\,\vert\,\{3,4\}\sqcup I_{g_2}) = \sum_{\substack{g_1,g_2\ge0\\g_1+g_2=h}}
        br^\ast \left[
        \begin{tikzpicture}[baseline=-3ex, scale=0.7,  transform shape]
            \tkzDefPoint(0,0){O} 
            \def\len{3} 
            \tkzDefPoint(35:{1/6*\len}){A}
            \tkzDefPoint(35:{-3/4*\len}){B}
            \tkzDefPoint(-35:{3/4*\len}){C}
            \tkzDefPoint(-35:{-1/6*\len}){D}
        
            \tkzDrawSegment[thick](B,A)
            \tkzDrawSegment[thick](D,C)
        
            \tkzDefPoint(35:{-5/8*\len}){zero}
            \draw[thick] ($(zero) + (0,2.5pt)$) -- ($(zero) - (0,2.5pt)$);
            \tkzLabelPoint[below](zero){$0$}

            \tkzDefPoint(35:{-15/32*\len}){P1}
            \draw[thick] ($(P1) + (0,2.5pt)$) -- ($(P1) - (0,2.5pt)$);
            \tkzLabelPoint[below=1.5pt](P1){$p_1$}

            \tkzDefPoint(35:{-9/24*\len}){p1}
            \draw[thick] ($(p1) + (0,2.5pt)$) -- ($(p1) - (0,2.5pt)$);
            
            \tkzDefPoint(35:{-3/24*\len}){p2}
            \draw[thick] ($(p2) + (0,2.5pt)$) -- ($(p2) - (0,2.5pt)$);

            \path (p1) -- (p2) node[midway, below, xshift=2pt, yshift = 12pt] {$\iddots$};

            \draw[decorate,decoration={brace,amplitude=3pt,mirror}] 
            ($(p1)-(-0.15,0.25)$) -- ($(p2)-(-0.15,0.25)$) node[midway,below right,font=\scriptsize] {$g_1$};

            \tkzDefPoint(-35:{15/32*\len}){P2}
            \draw[thick] ($(P2) + (0,2.5pt)$) -- ($(P2) - (0,2.5pt)$);
            \tkzLabelPoint[below=1pt](P2){$p_2$}

            \tkzDefPoint(-35:{5/8*\len}){infty}
            \draw[thick] ($(infty) + (0,2.5pt)$) -- ($(infty) - (0,2.5pt)$);
            \tkzLabelPoint[below=1pt](infty){$\infty$}
            
            \tkzDefPoint(-35:{8/24*\len}){p3}
            \draw[thick] ($(p3) + (0,2.5pt)$) -- ($(p3) - (0,2.5pt)$);

            \tkzDefPoint(-35:{3/24*\len}){p4}
            \draw[thick] ($(p4) + (0,2.5pt)$) -- ($(p4) - (0,2.5pt)$);

            \path (p3) -- (p4) node[midway, below, xshift=-1pt, yshift = 12pt] {$\ddots$};

            \draw[decorate,decoration={brace,amplitude=3.5pt}] 
            ($(p3)-(0.1,0.25)$) -- ($(p4)-(0.1,0.25)$) node[midway,below=2.5pt,font=\scriptsize] {$g_2$};
        \end{tikzpicture}
        \right].
    \end{equation*}
    For every $\sigma\in S_{\mu_1,\mu_2}^{(h)}$, the right hand side of \eqref{eq:wdvvtmu} is computed by summing over all possible ways to glue a smooth admissible cover of geometric genus $g_1$ lying over the twig with $0$ and $p_1$ (with ramification profile over $p_1$ described by $\mu_1$) to a smooth admissible cover of genus $g_2$ lying over the twig with $p_2$ and $\infty$ (with ramification profile over $p_2$ described by $\mu_2$) to obtain a nodal admissible cover with discrete data and ramification profiles specified by $\sigma$. Summing over $0\leq h\leq |\mu|$, and all possible $\sigma=(\mu_{1,i},\mu_{2,i},h_i,\ain_i)_{i=1}^m\in S_{\mu_1,\mu_2}^{(h)}$ we obtain
    \begin{gather}
        \nonumber \sum_{\sigma_1\in S_{\mu_1}}\sum_{\sigma_2\in S_{\mu_2}}\sum_{n\ge0}\frac{\h^n}{n!}\sum_{k_1,\dots,k_n>0} \\
        \nonumber \frac{\partial^n}{\partial p_{k_1}\cdots\partial p_{k_n}}\left(\frac{1}{m_1!}\h^{g_1-m_1}\e^{2g_1}\frac{v^{g_1}}{g_1!}\binom{g_1}{g_{1,1},\dots,g_{1,m_1}}\prod_{i_1=1}^{m_1}\frac{p_{\ain_{1,i_1}}}{(\lvert\mu_{1,i_1}\rvert+1-g_{1,i_1})!}q_{\lvert\ain_{1,i_1}\rvert}\right) \\
        \nonumber \prod_{i_1=1}^{m_1}\frac{x_{\mu_{1,i_1};g_{1,i_1}}(\ain_{1,i_1})}{\lvert\operatorname{Aut}(\mu_{1,i_1})\rvert} \, \prod_{j=1}^nk_n \, \prod_{i_2=1}^{m_2}\frac{x_{\mu_{2,i_2};g_{2,i_2}}(\ain_{2,i_2})}{\lvert\operatorname{Aut}(\mu_{2,i_2})\rvert} \\
        \nonumber \frac{\partial^n}{\partial q_{k_1}\cdots\partial    q_{k_n}}\left(\frac{1}{m_2!}\h^{g_2-m_2}\e^{2g_2}\frac{v^{g_2}}{g_2!}\binom{g_2}{g_{2,1},\dots,g_{2,m_2}}\prod_{i_2=1}^{m_2}\frac{p_{\ain_{2,i_2}}}{(\lvert\mu_{2,i_2}\rvert+1-g_{2,i_2})!}q_{\lvert\ain_{2,i_2}\rvert}\right) \\
        \label{eq:T1T2} = \Tr{\mu_1}^\bullet \star \Tr{\mu_2}^\bullet.
    \end{gather}

 We turn to the left-hand side of \eqref{eq:wdvvtmu}. The divisor here is:

    \begin{equation}\label{eq:Dsum}
    f^\ast D(\{1,4\}|\{2,3\})=\sum_{\substack{g_1+g_2=h\\ |I_{g_i}|= g_i\\I_{g_1}\cap I_{g_2} =\emptyset}}br^\ast D(\{1,4\} \sqcup I_{g_1}\,\vert\,\{2,3\}\sqcup I_{g_2}) = 
        br^\ast \sum_{g=0}^h 
        \left[
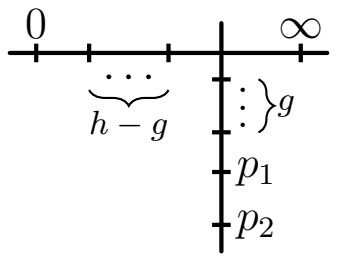
        \begin{tikzpicture}[baseline=-4ex, scale=0.7,  transform shape]]
            \tkzDefPoint(0,0){O} 
            \def\len{3} 
            \tkzDefPoint(0:{1/3*\len}){A}
            \tkzDefPoint(0:{-2/3*\len}){B}
            \tkzDefPoint(90:{3/32*\len}){C}
            \tkzDefPoint(90:{-5/8*\len}){D}
            \tkzDrawSegment[thick](B,A)
            \tkzDrawSegment[thick](D,C)
            \tkzDefPoint(0:{-7/12*\len}){zero}
            \draw[thick] ($(zero) + (0,2.5pt)$) -- ($(zero) - (0,2.5pt)$);
            \tkzLabelPoint[above](zero){$0$}
            \tkzDefPoint(0:{1/4*\len}){infty}
            \draw[thick] ($(infty) + (0,2.5pt)$) -- ($(infty) - (0,2.5pt)$);
            \tkzLabelPoint[above](infty){$\infty$}
            \tkzDefPoint(0:{-5/12*\len}){x1}
            \draw[thick] ($(x1) + (0,2.5pt)$) -- ($(x1) - (0,2.5pt)$);
            \tkzDefPoint(0:{-2/12*\len}){x2}
            \draw[thick] ($(x2) + (0,2.5pt)$) -- ($(x2) - (0,2.5pt)$);
            \path (x1) -- (x2) node[midway, below, xshift=1pt] {$\cdots$};
            \draw[decorate,decoration={brace,amplitude=3pt,mirror}] 
            ($(x1)-(0,0.35)$) -- ($(x2)-(0,0.35)$) node[midway,below=1.5pt,font=\scriptsize] {$h-g$};
            \tkzDefPoint(90:{-9/24*\len}){p1}
            \draw[thick] ($(p1) + (2.5pt,0)$) -- ($(p1) - (2.5pt,0)$);
            \tkzLabelPoint[right](p1){$p_1$}
            \tkzDefPoint(90:{-13/24*\len}){p2}
            \draw[thick] ($(p2) + (2.5pt,0)$) -- ($(p2) - (2.5pt,0)$);
            \tkzLabelPoint[right](p2){$p_2$}
            \tkzDefPoint(90:{-6/24*\len}){x3}
            \draw[thick] ($(x3) + (2.5pt,0)$) -- ($(x3) - (2.5pt,0)$);
            \tkzDefPoint(90:{-2/24*\len}){x4}
            \draw[thick] ($(x4) + (2.5pt,0)$) -- ($(x4) - (2.5pt,0)$);
            \path (x3) -- (x4) node[midway, right, yshift=3pt] {$\vdots$};
            \draw[decorate,decoration={brace,amplitude=3pt,mirror}] 
            ($(x3)+(0.35,0)$) -- ($(x4)+(0.35,0)$)      node[midway,right=1pt,font=\scriptsize] {$g$};
        \end{tikzpicture}
        \right].
    \end{equation}

We make some important observations about the invariants arising from pulling back the vertical twig of the boundary divisor  in \eqref{eq:Dsum}, see Figure \ref{fig:essentialcover} for an illustration:
\begin{itemize}
    \item the discrete data on the vertical twig does not depend on the total genus $h$ or on the profiles $\ain$. In particular, as we vary $\sigma$ along all of $S_{\mu_1, \mu_2}$, up to adding trivial connected components of degree one, by the Riemann Hurwitz formula there are only a finite number of possible options for the admissible covers  ``above" that twig. We call these finite number of admissible covers (ignoring the additional trivial degree one components that might be needed to get to the correct total degree) the {\it essential covers} of the vertical twig.
    \item because we are integrating the class $\lambda_h$, for a nonvanishing invariant above the vertical twig the number of non-special branch points on the twig must equal the geometric genus of the cover (and of course the dimension of the space of covers of the twig);
    \item because of the compact type condition needed for invariants not to vanish, the number of points above the node corresponding to the attaching points of the components of the essential part of the cover give a partition $\nu \vdash \delta'\leq \delta$ of length less than or equal to $\ell(\mu_1)+\ell(\mu_2)$. 
    \item essential covers of the vertical twig producing the same partition $\nu$ attach to the same covers of the horizontal twig with identical gluing factors. 
    \item any invariant from an essential cover can be promoted (by completing with ones the labeling of the inverse image over $p_3$ and rename those labels as $\ain$-type labels, to an invariant of type $x_{\vec{\mu}',g'}(\ain')$ which is inductively known. 
\end{itemize}
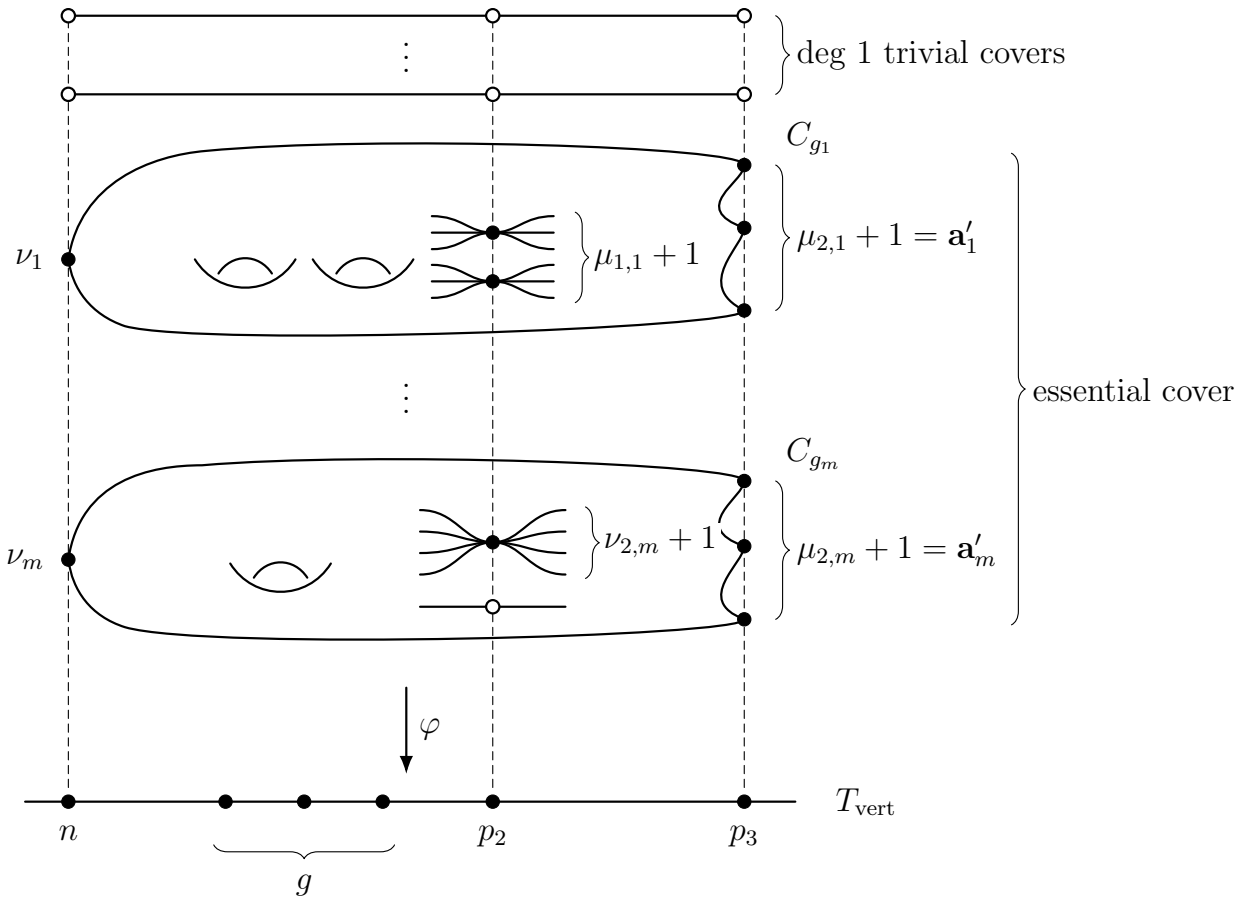
\begin{figure}
    \centering
  \resizebox{0.5 \textwidth}{!}{\begin{tikzpicture}[
  x=1cm,y=1cm,
  line cap=round,line join=round,
  dot/.style={circle,draw=black,fill=black,inner sep=1.7pt},
  open dot/.style={circle,draw=black,fill=white,thick,inner sep=1.7pt},
  guide/.style={draw=black,densely dashed,thin},
  strand/.style={draw=black,thick},
  lab/.style={draw=none,fill=white,inner sep=1.5pt},
  brace/.style={draw=black,decorate,decoration={brace,amplitude=5pt}},
  >=Latex
]

\def\xzero{0}
\def\xtwo{5.4}
\def\xinf{8.6}

\newcommand{\ram}[5]{%
  \ifnum#1=1
    \draw[strand] ({#2-#4/2},#3) -- (#2,#3) -- ({#2+#4/2},#3);
    \node[open dot] at (#2,#3) {};
  \else
    \node[dot] at (#2,#3) {};
    \foreach \k in {1,...,#1}{
      \pgfmathsetmacro{\ry}{#3-#5/2+(\k-1)*#5/(#1-1)}
      \draw[strand]
        ({#2-#4/2},\ry) .. controls ({#2-#4/4},\ry) and
        ({#2-#4/5},#3) .. (#2,#3)
        .. controls ({#2+#4/5},#3) and ({#2+#4/4},\ry) ..
        ({#2+#4/2},\ry);
    }
  \fi
}

\newcommand{\torus}[2]{%
  \draw[thick] ({#1-.64},#2)
    .. controls ({#1-.34},{#2-.48}) and ({#1+.34},{#2-.48}) ..
    ({#1+.64},#2);
  \draw[thick] ({#1-.34},{#2-.18})
    .. controls ({#1-.17},{#2+.08}) and ({#1+.17},{#2+.08}) ..
    ({#1+.34},{#2-.18});
}

\draw[thick] (-0.55,0) -- (9.25,0);
\foreach \x in {\xzero,\xtwo,\xinf} \node[dot] at (\x,0) {};
\node[lab,below=7pt] at (\xzero,0) {$n$};
\node[lab,below=7pt] at (\xtwo,0) {$p_2$};
\node[lab,below=7pt] at (\xinf,0) {$p_3$};
\node[lab,right=13pt] at (9.25,0) {$T_{\mathrm{vert}}$};
\node[dot] at (2,0) {};
\node[dot] at (3,0) {};
\node[dot] at (4,0) {};

\draw[brace,decoration={brace,amplitude=5pt,mirror}] (1.85,-0.55) -- (4.15,-0.55);
\node[lab,below=9pt] at (3,-0.55) {$g$};

\draw[thick] (0,10) -- (8.6,10);
\draw[thick] (0,9) -- (8.6,9);

\draw[brace] (9.00,10) -- (9.00,9);
\node[lab,right=6pt] at (9.00,9.5) {deg 1 trivial covers};

\draw[brace] (12.00,8.25) -- (12.00,2.25);
\node[lab,right=6pt] at (12.00,5.25) {essential cover};

\draw[guide] (\xzero,0.18) -- (\xzero,10);
\draw[guide] (\xtwo,0.18) -- (\xtwo,10);
\draw[guide] (\xinf,0.18) -- (\xinf,10);

\foreach \x in {\xzero,\xtwo,\xinf}{
  \node[open dot] at (\x,9) {};
  \node[open dot] at (\x,10) {};
}

\draw[->,thick] (4.30,1.45) -- (4.30,0.35);
\node[lab,right=3pt] at (4.30,0.92) {$\varphi$};

\draw[thick]
  (\xzero,3.08)
  .. controls (0.15,4.10) and (1.05,4.28) .. (1.70,4.28)
  .. controls (3.80,4.46) and (\xinf,4.30) .. (\xinf,4.08)
  .. controls (\xinf,3.88) and (7.88,3.50) .. (\xinf,3.25)
  .. controls (\xinf,3.04) and (7.94,2.57) .. (\xinf,2.32)
  .. controls (\xinf,2.10) and (1.75,1.92) .. (0.72,2.22)
  .. controls (0.28,2.36) and (0.05,2.69) .. (\xzero,3.08);
\node[dot] at (\xzero,3.08) {};
\node[lab,left=7pt] at (\xzero,3.08) {$\nu_m$};
\node[dot] at (\xinf,3.25) {};

\ram{4}{\xtwo}{3.30}{1.85}{0.82}
\ram{1}{\xtwo}{2.48}{1.85}{0}
\draw[brace] (6.58,3.75) -- (6.58,2.84);
\node[lab,right=5pt] at (6.58,3.30) {$ \nu_{2,m}+1$};
\torus{2.70}{3.03}
\node[dot] at (\xinf,4.08) {};
\node[dot] at (\xinf,2.32) {};
\draw[brace] (9.00,4.08) -- (9.00,2.32);
\node[lab,above right=3pt ] at (9.00,4.08) {$C_{g_m}$};
\node[lab,right=6pt] at (9.00,3.20) {$ \mu_{2,m}+1=\ain'_m$};

\node[draw=none] at (4.30,5.25) {$\vdots$};
\node[draw=none] at (4.30,9.6) {$\vdots$};
\draw[thick]
  (\xzero,6.90)
  .. controls (0.15,7.88) and (1.05,8.22) .. (1.75,8.28)
  .. controls (3.85,8.48) and (\xinf,8.34) .. (\xinf,8.10)
  .. controls (\xinf,7.88) and (7.88,7.56) .. (\xinf,7.30)
  .. controls (\xinf,7.08) and (7.94,6.58) .. (\xinf,6.25)
  .. controls (\xinf,6.02) and (1.75,5.78) .. (0.72,6.05)
  .. controls (0.28,6.20) and (0.05,6.53) .. (\xzero,6.90);
\node[dot] at (\xzero,6.90) {};
\node[lab,left=7pt] at (\xzero,6.90) {$\nu_1$};

\ram{3}{\xtwo}{6.62}{1.55}{0.42}
\ram{3}{\xtwo}{7.24}{1.55}{0.42}
\draw[brace] (6.45,7.50) -- (6.45,6.35);
\node[lab,right=5pt] at (6.45,6.92) {$ \mu_{1,1}+1$};
\torus{2.25}{6.90}
\torus{3.75}{6.90}

\node[dot] at (\xinf,8.10) {};
\node[dot] at (\xinf,7.30) {};
\node[dot] at (\xinf,6.25) {};
\draw[brace] (9.00,8.10) -- (9.00,6.25);
\node[lab,above right=3pt] at (9.00,8.10) {$C_{g_1}$};
\node[lab,right=6pt] at (9.00,7.18) {$ \mu_{2,1}+1= \ain'_1$};

\end{tikzpicture}} 
    \caption{A depiction of a cover of the vertical twig in \eqref{eq:Dsum}. By considering the profiles over $p_3$ as $\ain$-type conditions (after possibly labeling some unramified points), the invariant corresponding to the essential cover is an inductively known invariant.}
    \label{fig:essentialcover}
\end{figure}

These observations suggest to define by $c_{\mu_1, \mu_2}^{g,\nu}$ the sum of the contributions of all invariants coming from essential covers of the vertical twigs of geometric genus $g$ and with attaching partition $\nu$; we absorb into these coefficients also some of the combinatorial factors in a way that we make precise in Lemma \ref{lem:cmus} (we postpone this technical calculation to not break the flow of this proof).
Pending Lemma \ref{lem:cmus}, the left hand side of \eqref{eq:wdvvtmu} evaluates to
    \begin{align}\label{eq:cnununu}
        \frac{\lvert\operatorname{Aut}(\mu)\rvert}{\lvert\operatorname{Aut}(\mu_1)\rvert\lvert\operatorname{Aut}(\mu_2)\rvert}\Tr{\mu}^\bullet & + \sum_{\substack{\nu\vdash\lvert\mu\rvert\\\ell(\nu)<\ell(\mu)}}c_{\mu_1,\mu_2}^{0,\nu}\Tr{\nu}^\bullet + \sum_{g=1}^{\lvert\mu\rvert-\mu_1}\h^g\e^{2g}\frac{v^g}{g!}\hspace{-1.5em}\sum_{\substack{\nu\vdash\lvert\mu\rvert-g\\\ell(\nu)\le\max\{\ell(\mu)-g,1\}}}\hspace{-1.5em}c_{\mu_1,\mu_2}^{g,\nu}\Tr{\nu}^\bullet.
    \end{align}
    Equating this to \eqref{eq:T1T2} yields:
    \begin{align*}
        \Tr{\mu}^\bullet = & \; \frac{\lvert\operatorname{Aut}(\mu_1)\rvert\lvert\operatorname{Aut}(\mu_2)\rvert}{\lvert\operatorname{Aut}(\mu)\rvert}\Tr{\mu_1}^\bullet\star\Tr{\mu_2}^\bullet - \frac{\lvert\operatorname{Aut}(\mu_1)\rvert\lvert\operatorname{Aut}(\mu_2)\rvert}{\lvert\operatorname{Aut}(\mu)\rvert}\sum_{\substack{\nu\vdash\lvert\mu\rvert\\\ell(\nu)<\ell(\mu)}}c_{\mu_1,\mu_2}^{0,\nu}\Tr{\nu}^\bullet \\
        & - \frac{\lvert\operatorname{Aut}(\mu_1)\rvert\lvert\operatorname{Aut}(\mu_2)\rvert}{\lvert\operatorname{Aut}(\mu)\rvert}\sum_{g=1}^{\lvert\mu\rvert-\mu_1}\h^g\e^{2g}\frac{v^g}{g!}\hspace{-1.5em}\sum_{\substack{\nu\vdash\lvert\mu\rvert-g\\\ell(\nu)\le\max\{\ell(\mu)-g,1\}}}\hspace{-1.5em}c_{\mu_1,\mu_2}^{g,\nu}\Tr{\nu}^\bullet.
    \end{align*}
   All the generating series on the right-hand side are indexed by partitions that are strictly smaller than $\mu$ in the partial order \eqref{eq:partial_order}. This completes the proof.
\end{proof}

\begin{lem}\label{lem:cmus}
With notation as in the proof of Proposition \ref{prop:TmuTk}, we have:
      \begin{align*}
          c_{\mu_1,\mu_2}^{g,\nu}=& \sum_{\substack{\rho_{1,1},\dots,\rho_{1,l}\\ \coprod_{k}\rho_{1,k}=\mu_1}}\sum_{\substack{\rho_{2,1},\dots,\rho_{2,l}\\ \coprod_{k}\rho_{2,k}=\mu_2}} \sum_{\substack{g_1,\dots,g_l\ge0\\g_1+\cdots+g_l=g}}\binom{g}{2g_1+\nu_1-\lvert\rho_{1,1}\rvert-\lvert\rho_{2,1}\rvert,\dots,2g_l+\nu_l-\lvert\rho_{1,l}\rvert-\lvert\rho_{2,l}\rvert} \\
         & \hspace{5em} \prod_{k=1}^l \frac{\nu_k+1}{\lvert\operatorname{Aut}(\rho_{1,k})\rvert\lvert\operatorname{Aut}(\rho_{2,k})\rvert}\int_{\adm_{g_k}((\nu_k+1),\,(\rho_{1,k}+1),\,(\rho_{2,k}+1))}\lambda_{g_k}.
    \end{align*}
    Moreover, $c_{\mu_1,\mu_2}^{g,\nu}=0$ whenever any of the following statements are true:
    \begin{itemize}
        \item $\ell(\nu)>\max\{\ell(\mu)-g,1\}$;
        \item $\ell(\nu)=\ell(\mu)$ and $\nu\neq\mu$;
        \item $g>\lvert\mu\rvert-\mu_1$ (i.e. $\lvert\nu\rvert<\mu_1$).
    \end{itemize}
    When $\nu=\mu$ we have $c_{\mu_1,\mu_2}^{0,\mu} = \frac{\lvert\operatorname{Aut}(\mu)\rvert}{\lvert\operatorname{Aut}(\eta_1)\rvert\lvert\operatorname{Aut}(\eta_2)\rvert}$.
\end{lem}
\begin{proof}
We rewrite the right hand side of \eqref{eq:wdvvtmu} multiplied by the monomial \eqref{eq:auxiliary_1} reindexing the genera of the covers of the two components so that the horizontal component of \eqref{eq:Dsum} is covered by curves of genus $h$:
    \begin{align} \label{eq:righthandside}
          & \sum_{g=0}^{\lvert\mu\rvert}\h^g\e^{2g}\frac{v^g}{g!}\sum_{h=0}^{\lvert\mu\rvert-g}\hspace{-1em}\sum_{\substack{\sigma\in S_{\mu_1,\mu_2}^{(g+h)}\\\sigma=(\mu_{1,i},\mu_{2,i},k_i,\ain_i)_{i=1}^m}}\hspace{-1em}\frac{1}{m!}\h^{h-m}\e^{2h}\frac{v^h}{h!}\frac{g!h!}{(g+h)!}\binom{g+h}{k_1,\dots,k_m}\frac{1}{\prod_{j,i}\lvert\operatorname{Aut}(\mu_{j,i})\rvert} \\
        \nonumber & \hspace{7em} \int_{\adm_\sigma^\bullet}\lambda_{g+h}\,br^\ast  
        \left[
        \begin{tikzpicture}[baseline=-4ex, scale=0.7,  transform shape]
            \tkzDefPoint(0,0){O} 
            \def\len{3} 
            \tkzDefPoint(0:{1/3*\len}){A}
            \tkzDefPoint(0:{-2/3*\len}){B}
            \tkzDefPoint(90:{3/32*\len}){C}
            \tkzDefPoint(90:{-5/8*\len}){D}
            \tkzDrawSegment[thick](B,A)
            \tkzDrawSegment[thick](D,C)
            \tkzDefPoint(0:{-7/12*\len}){zero}
            \draw[thick] ($(zero) + (0,2.5pt)$) -- ($(zero) - (0,2.5pt)$);
            \tkzLabelPoint[above](zero){$0$}
            \tkzDefPoint(0:{1/4*\len}){infty}
            \draw[thick] ($(infty) + (0,2.5pt)$) -- ($(infty) - (0,2.5pt)$);
            \tkzLabelPoint[above](infty){$\infty$}
            \tkzDefPoint(0:{-5/12*\len}){x1}
            \draw[thick] ($(x1) + (0,2.5pt)$) -- ($(x1) - (0,2.5pt)$);
            \tkzDefPoint(0:{-2/12*\len}){x2}
            \draw[thick] ($(x2) + (0,2.5pt)$) -- ($(x2) - (0,2.5pt)$);
            \path (x1) -- (x2) node[midway, below, xshift=1pt] {$\cdots$};
            \draw[decorate,decoration={brace,amplitude=3pt,mirror}] 
            ($(x1)-(0,0.35)$) -- ($(x2)-(0,0.35)$) node[midway,below=1.5pt,font=\scriptsize] {$h$};
            \tkzDefPoint(90:{-9/24*\len}){p1}
            \draw[thick] ($(p1) + (2.5pt,0)$) -- ($(p1) - (2.5pt,0)$);
            \tkzLabelPoint[right](p1){$p_1$}
            \tkzDefPoint(90:{-13/24*\len}){p2}
            \draw[thick] ($(p2) + (2.5pt,0)$) -- ($(p2) - (2.5pt,0)$);
            \tkzLabelPoint[right](p2){$p_2$}
            \tkzDefPoint(90:{-6/24*\len}){x3}
            \draw[thick] ($(x3) + (2.5pt,0)$) -- ($(x3) - (2.5pt,0)$);
            \tkzDefPoint(90:{-2/24*\len}){x4}
            \draw[thick] ($(x4) + (2.5pt,0)$) -- ($(x4) - (2.5pt,0)$);
            \path (x3) -- (x4) node[midway, right, yshift=3pt] {$\vdots$};
            \draw[decorate,decoration={brace,amplitude=3pt,mirror}] 
            ($(x3)+(0.35,0)$) -- ($(x4)+(0.35,0)$)      node[midway,right=1pt,font=\scriptsize] {$g$};
        \end{tikzpicture}
        \right] \prod_{i=1}^m\frac{p_{\ain_i}}{(\lvert\mu_{1,i}\rvert+\lvert\mu_{2,i}\rvert+1-k_i)!}q_{\lvert\ain_i\rvert}.
    \end{align}
    For every $\sigma=(\mu_{1,i},\mu_{2,i},k_i,\ain_i)_{i=1}^m\in S_{\mu_1,\mu_2}^{(g+h)}$, over the horizontal twig, we have  invariants on a space of admissible covers of geometric genus $h$ and $m$ connected components. The $i$-th component is described on the left hand side of Figure \ref{fig:pieces}.

\begin{figure}
    \centering
   \begin{tikzpicture}[thick, font=\small, scale=0.7,  transform shape]]
            \node[draw, circle, minimum size=1.2cm, label={above left:$\lambda_{h_i}$}] (h) at (0,0) {$h_i$};
            \draw (h.west) -- ++(-1.0, 0) node[left] {$\delta_i$};
            \draw (h.105) -- ++(-0.175, 0.35) coordinate (t1_2);
            \draw (h.75) -- ++(0.175, 0.35) coordinate (t2_2);
            \draw[dotted, thick] (-0.15, 0.9) to[bend left=10] (0.22, 0.9);
            \draw[decorate, decoration={brace, amplitude=4pt, raise=2pt}] 
            (t1_2) -- (t2_2) node[midway, above=6pt] {$h_i$};
            \draw (h.20) -- ++(1.0, 0.4) node[right] {$a_{i,1}$};
            \draw (h.-20) -- ++(1.0, -0.4) node[right] {$a_{i,\lvert\eta_{1,i}\rvert+\lvert\eta_{2,i}\rvert+1-k_i}$};
            \draw[loosely dotted, very thick] (1.5, 0.35) to[bend left=25] (1.5, -0.4);
            \draw (h.-120) -- ++(-0.6, -1.0) node[below left] {$\nu_{i,1}^\prime+1$};
            \draw (h.-90) -- ++(0, -1.1) node[below] {$\nu_{i,\ell_i}^\prime+1$};
            \draw (h.-60) -- ++(0.6, -1.0) node[below right] {$1^{\delta_i-\lvert\nu^\prime_i\rvert-\ell_i}$};
            \draw[loosely dotted, very thick] (-0.7, -1.45) to[bend right=15] (-0.05, -1.55);
        \end{tikzpicture} \ \ \ \ \ \ \ \ 
        \begin{tikzpicture}[thick, font=\small, scale=0.7,  transform shape]]
            \node[draw, circle, minimum size=1.2cm, label={right:$\lambda_{g_{i,s}}$}] (g) at (0,0) {$g_{i,s}$};
            \draw (g.north) -- ++(0, 1.0) node[above] {$\nu^\prime_{i,s}+1$};
            \draw (g.-115) -- ++(-0.5, -0.71) node[below left, xshift=10mm] {$(\rho_{1,i,s}+1)$};
            \draw (g.-65) -- ++(0.5, -0.71) node[below right, xshift=-2mm] {$(\rho_{2,i,s}+1)$};
            \draw (g.155) -- ++(-0.3, 0.175) coordinate (l1);
            \draw (g.205) -- ++(-0.3, -0.175) coordinate (l2);
            \draw[dotted, thick] (-0.8, 0.2) to[bend right=15] (-0.8, -0.2);
            \draw[decorate, decoration={brace, amplitude=4pt, raise=4pt}] 
            (l2) -- (l1) node[midway, left=8pt] {$2g_{i,s}+\nu^\prime_{i,s}-\lvert\rho_{1,i,s}\rvert-\lvert\rho_{2,i,s}\rvert$};
        \end{tikzpicture}
    \caption{The connected components of the admissible covers mapping to the  horizontal (left) and vertical (right) twig of the divisor $f^\ast D(\{1,4\}|\{2,3\})$. }
    \label{fig:pieces}
\end{figure}
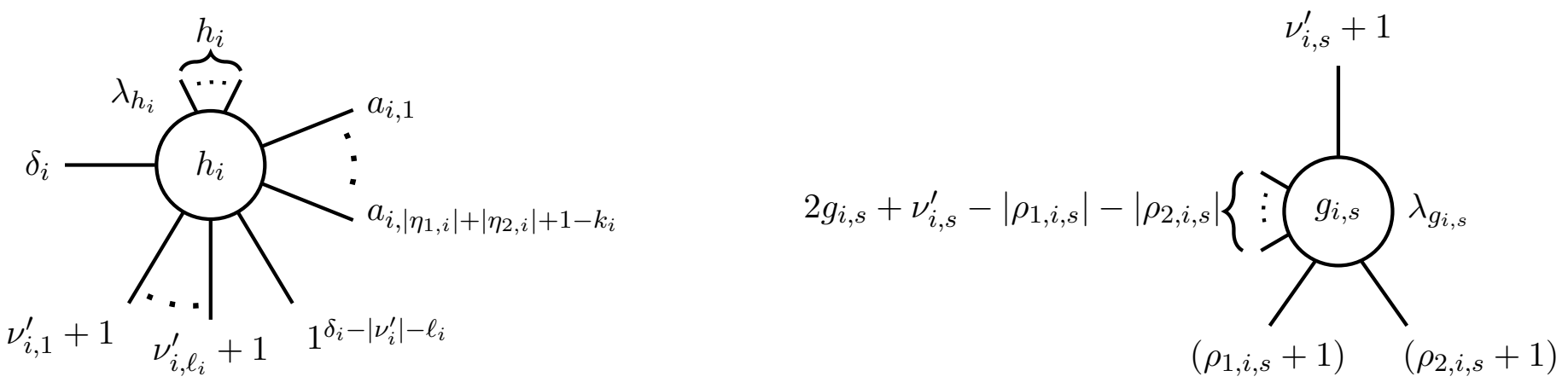
    The non-trivial ramification profiles over the node are given by  a partition $\nu_i^\prime$  of some length $\ell_i$. Since there are $h_i$ additional simple branch points, the size of $\nu_i^\prime$ is  $h_i-k_i+\lvert\eta_{1,i}\rvert+\lvert\eta_{2,i}\rvert$ by the Riemann-Hurwitz formula. The invariant from this admissible cover space is:
    \begin{equation*}
        \prod_{i=1}^m\frac{x_{\nu^\prime_i;h_i}(\ain_i)}{\lvert\operatorname{Aut}(\nu_i^\prime)\rvert}.
    \end{equation*}
    For every $i=1,\dots,m$ and every $s=1,\dots,\ell_i$, there is a curve glued over the node to the point labeled by $\nu_{i,s}^\prime+1$ mapping to the vertical twig, with cover depicted on the right hand side of Figure \ref{fig:pieces}.
    The ramification profile over $p_j$ is
    \begin{equation*}
        (\rho_{j,i,s}+1) = \bigl(\rho_{j,i,s}^{(1)}+1,\dots,\rho_{j,i,s}^{(\ell_{j,i,s})}+1,1^{\nu^\prime_{i,s}+1-\lvert\rho_{j,i,s}\rvert-\ell_{j,i,s}}\bigr), \quad \ell_{j,i,s}=\ell(\rho_{j,i,s}),
    \end{equation*}
    for some possibly empty partition $\rho_{j,i,s}$, and there are $2g_{i,s}+\nu^\prime_{i,s}-\lvert\rho_{1,i,s}\rvert-\lvert\rho_{2,i,s}\rvert$ additional simple branch points. To match the data specified by $\sigma=(\mu_{1,i},\mu_{2,i},k_i,\ain_i)_{i=1}^m$, it is necessary that $\sum_{s}g_{i,s}=k_i-h_i$ and $\coprod_{s}\rho_{j,i,s}=\mu_{j,i}$. The invariant associated to the vertical twig  is
    \begin{equation*}       \prod_{i=1}^m\prod_{s=1}^{\ell_i}\frac{1}{\lvert\operatorname{Aut}(\rho_{1,i,s})\rvert\lvert\operatorname{Aut}(\rho_{2,i,s})\rvert}\int_{\adm_{g_{i,s}}\bigl((\nu^\prime_{i,s}+1),(\rho_{1,i,s}+1),(\rho_{2,i,s}+1)\bigr)}\lambda_{g_{i,s}}.
    \end{equation*}
    
    With these evaluations, the gluing factor $\prod_{i,s}(\nu_{i,s}^\prime+1)$ and the ways to distribute the additional simple branch points on the various connected components, \eqref{eq:righthandside} is equal to
    \begin{gather*}
        \sum_{\substack{h_1,\dots,h_m\ge0\\h_1+\cdots+h_m=h}} \, \sum_{\substack{\nu^\prime_1,\dots,\nu^\prime_m\\\nu^\prime_i\vdash h_i-k_i+\lvert\mu_{1,i}\rvert+\lvert\mu_{2,i}\rvert}}\prod_{i=1}^m\binom{k_i}{h_i,k_i-h_i}\frac{x_{\nu^\prime_i;h_i}(\ain_i)}{\lvert\operatorname{Aut}(\nu^\prime_i)\rvert} \cdot \prod_{i=1}^m\prod_{s=1}^{\ell_i}(\nu^\prime_{i,s}+1) \\
        \sum_{\substack{\rho_{1,1,1},\dots,\rho_{1,1,\ell_1}\\ \coprod_s\rho_{1,1,s}=\mu_{1,1}}} \cdots \sum_{\substack{\rho_{1,m,1},\dots,\rho_{1,m,\ell_m}\\ \coprod_s\rho_{1,m,s}=\mu_{1,m}}} \, \sum_{\substack{\rho_{2,1,1},\dots,\rho_{2,1,\ell_1}\\ \coprod_s\rho_{2,1,s}=\mu_{2,1}}} \cdots \sum_{\substack{\rho_{2,m,1},\dots,\rho_{2,m,\ell_m}\\ \coprod_s\rho_{2,m,s}=\mu_{2,m}}} \, \sum_{\substack{g_{1,1},\dots,g_{1,\ell_1}\ge0\\\sum_sg_{1,s}=k_1-h_1}} \cdots \sum_{\substack{g_{m,1},\dots,g_{m,\ell_m}\ge0\\\sum_sg_{m,s}=k_m-h_m}} \\
        \prod_{i=1}^m\binom{k_i-h_i}{2g_{i,1}+\nu^\prime_{i,1}-\lvert\rho_{1,i,1}\rvert-\lvert\rho_{2,i,1}\rvert,\dots,2g_{i,\ell_i}+\nu^\prime_{i,\ell_i}-\lvert\rho_{1,i,\ell_i}\rvert-\lvert\rho_{2,i,\ell_i}\rvert} \\ \prod_{i=1}^m\prod_{s=1}^{\ell_i}\frac{1}{\lvert\operatorname{Aut}(\rho_{1,i,s})\rvert\lvert\operatorname{Aut}(\rho_{2,i,s})\rvert}\int_{\adm_{g_{i,s}}\bigl((\nu^\prime_{i,s}+1),(\rho_{1,i,s}+1)_{\nu^\prime_{i,s}+1},(\rho_{2,i,s}+1)_{\nu^\prime_{i,s}+1}\bigr)}\lambda_{g_{i,s}}.
    \end{gather*}
    For all partitions $\nu^\prime_1,\dots,\nu^\prime_m$ appearing in the sum above, let $\nu:=\nu^\prime_1\sqcup\dots\sqcup\nu^\prime_m$, and denote the individual elements of $\nu$ by $\nu_1,\dots,\nu_l$, where $l:=\ell(\nu)$. The total ramification degree over the node is independent of $h$:
    \begin{equation*}
        \lvert\nu\rvert = \sum_{i=1}^m\lvert\nu^\prime_i\rvert = h-(g+h)+\lvert\mu_1\rvert+\lvert\mu_2\rvert = \lvert\mu\rvert - g.
    \end{equation*}
 When we sum over all $0\le h\le\lvert\mu\rvert-g$ and $\sigma\in S_{\mu,d}^{(g+h)}$, we obtain
    \begin{align}
        \nonumber & \sum_{g=0}^{\lvert\mu\rvert}\h^g\e^{2g}\frac{v^g}{g!} \sum_{\nu\vdash\lvert\mu\rvert-g}\Biggl(\sum_{m\ge1}\frac{1}{m!}\sum_{\substack{\nu^\prime_1,\dots,\nu^\prime_m\\ \coprod_i\nu_i^\prime=\nu}}\sum_{h_1,\dots,h_m\ge0}\h^{h-m}\e^{2h}\frac{v^h}{h!}\binom{h}{h_1,\dots,h_m} \\
        \nonumber & \hspace{7em}\sum_{\ain_1\in\mathbb{N}^{\lvert\nu^\prime_1\rvert+1-h_1}}\cdots\sum_{\ain_m\in\mathbb{N}^{\lvert\nu^\prime_m\rvert+1-h_m}} \prod_{i=1}^m\frac{x_{\nu^\prime_i;h_i}(\ain_i)}{\lvert\operatorname{Aut}(\nu^\prime_i)\rvert}\frac{p_{\ain_i}}{(\lvert\nu^\prime_i\rvert+1-h_i)!}{q_{\lvert\ain_i\rvert}}\Biggr) \\
        \nonumber & \sum_{\substack{\rho_{1,1},\dots,\rho_{1,l}\\ \coprod_{k}\rho_{1,k}=\mu_1}}\sum_{\substack{\rho_{2,1},\dots,\rho_{2,l}\\ \coprod_{k}\rho_{2,k}=\mu_2}} \sum_{\substack{g_1,\dots,g_l\ge0\\g_1+\cdots+g_l=g}}\binom{g}{2g_1+\nu_1-\lvert\rho_{1,1}\rvert-\lvert\rho_{2,1}\rvert,\dots,2g_l+\nu_l-\lvert\rho_{1,l}\rvert-\lvert\rho_{2,l}\rvert} \\
        \label{eq:int_num_cnu} & \hspace{5em} \prod_{k=1}^l \frac{\nu_k+1}{\lvert\operatorname{Aut}(\rho_{1,k})\rvert\lvert\operatorname{Aut}(\rho_{2,k})\rvert}\int_{\adm_{g_k}((\nu_k+1),\,(\rho_{1,k}+1),\,(\rho_{2,k}+1))}\lambda_{g_k} \\
        \nonumber & = \sum_{g=0}^{\lvert\mu\rvert}\h^g\e^{2g}\frac{v^g}{g!} \sum_{\nu\vdash\lvert\mu\rvert-g}\Tr{\nu}^\bullet\,c_{\mu_1,\mu_2}^{g,\nu},
    \end{align}
    as desired. The vanishing statement of the Lemma follows from the observation that the intersection number \eqref{eq:int_num_cnu} vanishes whenever there exists some $1\le k\le\ell$ such that $\lvert\rho_{1,k}\rvert+\ell(\rho_{1,k})>\nu_k+1$ or $\lvert\rho_{2,k}\rvert+\ell(\rho_{2,k})>\nu_k+1$ or $(\rho_{1,k},\rho_{2,k})=(\emptyset,\emptyset)$. For $\nu=\mu$,
    the intersection number vanishes unless $(\rho_{1,k},\rho_{2,k})=((\mu_k),\emptyset)$ or $(\rho_{1,k},\rho_{2,k})=(\emptyset,(\mu_k))$ for every $k$: $\prod_{k=1}^\ell(\mu_k+1)\frac{1}{\mu_k+1}=1$ by Lemma \ref{lem:disceval}. Therefore $c^{0,\mu}_{\mu_1,\mu_2}$ is the number of ways of picking partitions $\rho_{1,1}\sqcup\cdots\sqcup\rho_{1,\ell}=\mu_1$ and $\rho_{2,1}\sqcup\cdots\sqcup\rho_{2,\ell}=\mu_2$ satisfying such conditions, which is $\frac{\lvert\operatorname{Aut}(\mu)\rvert}{\lvert\operatorname{Aut}(\mu_1)\rvert\lvert\operatorname{Aut}(\mu_2)\rvert}$.
\end{proof}

\subsection{Reconstructing the $\Tr{\mu,d}^\bullet$} 

This section is devoted to showing that all the invariants that we are considering, namely $x_{\mu,d;g}(\ain)$ and their disconnected counterparts $x_{\vec\mu,d;\vec{g}}(\vec{\ain})$, can be expressed in terms of the much simpler building blocks $x_{k;g}(\ain)$. It is natural to express these results in terms of relations between the generating series $\Tr{\mu,d}^\bullet$ and $\Trz{k}$. 

\begin{prop}\label{prop:Tds}
    For every partition $\mu$ and every integer $d\ge1$,
    \begin{align}
        \label{eq:TmudTmuprop}\Tr{\mu,d}^\bullet = \frac{v}{d}\Trz{1}\star\Tr{\mu,d-1}^\bullet & - \sum_{g=1}^{d-1}\h^g\e^{2g}\frac{v^{d+g}}{(d+g)!}\frac{g}{d}\sum_{\nu\vdash\lvert\mu\rvert+d-g}\gamma_{\mu,d-1}^{g,\nu}\Tr{\nu}^\bullet \\
        \nonumber & - \h^d\e^{2d}\frac{v^{2d}}{(2d)!}\gamma_{\mu,d-1}^{d,\mu}\Tr{\mu}^\bullet,
    \end{align}
    where $\gamma_{\mu,d-1}^{g,\nu}$ is the coefficient of 
    \begin{equation*}
        \h^{g-\ell(\nu)}\e^{2g}\frac{v^{d+g-1}}{(d+g-1)!}\frac{1}{\lvert\operatorname{Aut}(\nu)\rvert}\prod_{i=1}^{\ell(\nu)}\frac{p_1^{\nu_i+1}}{\nu_i+1}q_{\nu_i+1}
    \end{equation*} 
    in
    \begin{equation*}
        \left(p_1^2\frac{\partial}{\partial p_2} + \frac{1}{\h}\frac{p_1^2}{2}q_2\right)\Tr{\mu,d-1}^\bullet.
    \end{equation*}
\end{prop}
\begin{proof}
    For every $0\le h\le\lvert\mu\rvert+d$, let $S_{\mu,d}^{(h)}$ denote the subset of $S_{\mu,d}$ from Definition \ref{def:Smud} consisting of elements $\sigma=(\mu_i,g_i,n_i,\ain_i)_{i=1}^m$ for which $\sum_ig_i=h$, i.e. the geometric genus of the associated admissible cover is $h$. Then
    
    \begin{align}
        \label{eq:general_Tmud}\Tr{\mu,d}^\bullet = & \; \sum_{h=0}^{\lvert\mu\rvert+d}\sum_{\sigma\in S_{\mu,d}^{(h)}}\frac{1}{m!}\h^{h-m}\e^{2h}\frac{v^{d+h}}{(d+h)!}\binom{d+h}{r_1,\dots,r_m} \\
        \nonumber & \hspace{5em} \frac{1}{\prod_i\lvert\operatorname{Aut}(\mu_i)\rvert}\left(\int_{\adm_\sigma^\bullet}\lambda_h\,b^*(\psi_\bullet^d)\right)\prod_{i=1}^m\frac{p_{\ain_i}}{n_i!}q_{\lvert\ain_i\rvert}.
    \end{align}
    In this proof we denote the branch point with ramification profile described by $\mu$, which is the one on which we take the $\psi$-class, by the symbol $\bullet$, and we denote the $r=d+h$ additional branch points with simple ramification by $p_1,\dots,p_r$. One can write the class $\psi_\bullet$ on $\overline{\mathcal{M}}_{0,\{0,\bullet,p_1,\dots,p_r,\infty\}}$ as a sum of divisors by separating $\bullet$ from $0$ and one of the simple branch points $p_i$, or by separating $\bullet$ from $0$ and $\infty$. We choose the following linear combination of these possibilities:
    \begin{equation}\label{eq:psi_choice}
        \psi_\bullet = \frac{1}{d}\sum_{i=1}^{d+h}
        \left[
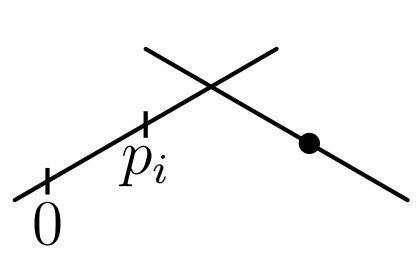
        \begin{tikzpicture}[baseline=-2ex]
            \tkzDefPoint(0,0){O} 
            \def\len{2} 
            \tkzDefPoint(30:{1/4*\len}){A}
            \tkzDefPoint(30:{-3/4*\len}){B}
            \tkzDefPoint(-30:{3/4*\len}){C}
            \tkzDefPoint(-30:{-1/4*\len}){D}
        
            \tkzDrawSegment[thick](B,A)
            \tkzDrawSegment[thick](D,C)
        
            \tkzDefPoint(30:{-5/8*\len}){zero}
            \draw[thick] ($(zero) + (0,2.5pt)$) -- ($(zero) - (0,2.5pt)$);
            \tkzLabelPoint[below](zero){$0$}

            \tkzDefPoint(30:{-1/4*\len}){pi}
            \draw[thick] ($(pi) + (0,2.5pt)$) -- ($(pi) - (0,2.5pt)$);
            \tkzLabelPoint[below](pi){$p_i$}

            \tkzDefPoint(-30:{3/8*\len}){p}
            \fill (p) circle (2pt);
        \end{tikzpicture}
        \right]
        - \frac{h}{d}
        \left[
        
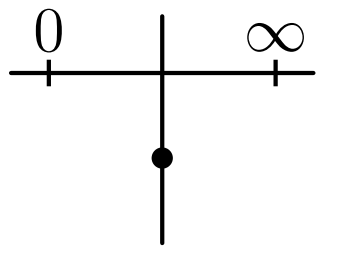
\begin{tikzpicture}[baseline=-2.5ex]
            \tkzDefPoint(0,0){O} 
            \def\len{3} 
            \tkzDefPoint(0:{1/3*\len}){A}
            \tkzDefPoint(0:{-1/3*\len}){B}
            \tkzDefPoint(90:{1/8*\len}){C}
            \tkzDefPoint(90:{-3/8*\len}){D}
        
            \tkzDrawSegment[thick](B,A)
            \tkzDrawSegment[thick](D,C)
        
            \tkzDefPoint(0:{-1/4*\len}){zero}
            \draw[thick] ($(zero) + (0,2.5pt)$) -- ($(zero) - (0,2.5pt)$);
            \tkzLabelPoint[above](zero){$0$}

            \tkzDefPoint(0:{1/4*\len}){infty}
            \draw[thick] ($(infty) + (0,2.5pt)$) -- ($(infty) - (0,2.5pt)$);
            \tkzLabelPoint[above](infty){$\infty$}

            \tkzDefPoint(90:{-3/16*\len}){psi}
            \fill (psi) circle (2pt);
        \end{tikzpicture}
        \right].
    \end{equation}
    Therefore $\psi_\bullet^d$ is a linear combination of terms of the form
    \begin{equation*}
       D(0;r_1 \,\vert\, \bullet,\infty;r_2)\psi_\bullet^{d-1}:= \sum_{\substack{|I_{r_i}| = r_i\\ I_{r_1}\cap I_{r_2} = \emptyset}} D(\{0\}\sqcup I_{r_1}\,\vert\,\{\bullet,\infty\} \sqcup I_{r_2})\psi_\bullet^{d-1} =
        \left[
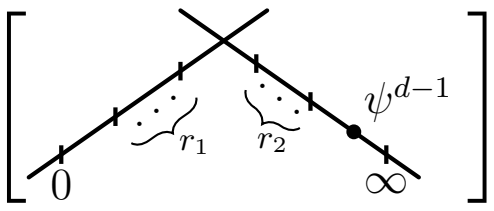
        \begin{tikzpicture}[baseline=-3ex, scale=0.7,  transform shape]
            \tkzDefPoint(0,0){O} 
            \def\len{3} 
            \tkzDefPoint(35:{1/6*\len}){A}
            \tkzDefPoint(35:{-3/4*\len}){B}
            \tkzDefPoint(-35:{3/4*\len}){C}
            \tkzDefPoint(-35:{-1/6*\len}){D}
        
            \tkzDrawSegment[thick](B,A)
            \tkzDrawSegment[thick](D,C)
        
            \tkzDefPoint(35:{-5/8*\len}){zero}
            \draw[thick] ($(zero) + (0,2.5pt)$) -- ($(zero) - (0,2.5pt)$);
            \tkzLabelPoint[below](zero){$0$}

            \tkzDefPoint(35:{-5/12*\len}){p1}
            \draw[thick] ($(p1) + (0,2.5pt)$) -- ($(p1) - (0,2.5pt)$);
            
            \tkzDefPoint(35:{-2/12*\len}){p2}
            \draw[thick] ($(p2) + (0,2.5pt)$) -- ($(p2) - (0,2.5pt)$);

            \path (p1) -- (p2) node[midway, below, xshift=2pt, yshift = 12pt] {$\iddots$};

            \draw[decorate,decoration={brace,amplitude=3pt,mirror}] 
            ($(p1)-(-0.15,0.25)$) -- ($(p2)-(-0.15,0.25)$) node[midway,below right,font=\scriptsize] {$r_1$};

            \tkzDefPoint(-35:{4/8*\len}){psi}
            \fill (psi) circle (2pt);
            \tkzLabelPoint[above right=-2pt](psi){$\psi^{d-1}$}

            \tkzDefPoint(-35:{5/8*\len}){infty}
            \draw[thick] ($(infty) + (0,2.5pt)$) -- ($(infty) - (0,2.5pt)$);
            \tkzLabelPoint[below=1pt](infty){$\infty$}
            
            \tkzDefPoint(-35:{8/24*\len}){p3}
            \draw[thick] ($(p3) + (0,2.5pt)$) -- ($(p3) - (0,2.5pt)$);

            \tkzDefPoint(-35:{3/24*\len}){p4}
            \draw[thick] ($(p4) + (0,2.5pt)$) -- ($(p4) - (0,2.5pt)$);

            \path (p3) -- (p4) node[midway, below, xshift=-1pt, yshift = 12pt] {$\ddots$};

            \draw[decorate,decoration={brace,amplitude=3.5pt}] 
            ($(p3)-(0.1,0.25)$) -- ($(p4)-(0.1,0.25)$) node[midway,below=2.5pt,font=\scriptsize] {$r_2$};
        \end{tikzpicture}
        \right]
    \end{equation*}
    and
    \begin{equation*}
        D(0,\infty;r_1 \,\vert\, \bullet;r_2)\psi_\bullet^{d-1}:= \sum_{\substack{|I_{r_i}| = r_i\\ I_{r_1}\cap I_{r_2} = \emptyset}}D(\{0, \infty\}\sqcup I_{r_1}\,\vert\,\{\bullet\} \sqcup I_{r_2})\psi_\bullet^{d-1} =
        \left[
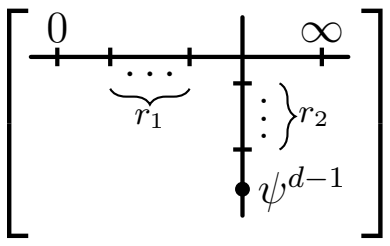
        \begin{tikzpicture}[baseline=-3ex, scale=0.7,  transform shape]
            \tkzDefPoint(0,0){O} 
            \def\len{3} 
            \tkzDefPoint(0:{1/3*\len}){A}
            \tkzDefPoint(0:{-2/3*\len}){B}
            \tkzDefPoint(90:{1/8*\len}){C}
            \tkzDefPoint(90:{-1/2*\len}){D}
        
            \tkzDrawSegment[thick](B,A)
            \tkzDrawSegment[thick](D,C)
        
            \tkzDefPoint(0:{-7/12*\len}){zero}
            \draw[thick] ($(zero) + (0,2.5pt)$) -- ($(zero) - (0,2.5pt)$);
            \tkzLabelPoint[above](zero){$0$}

            \tkzDefPoint(0:{1/4*\len}){infty}
            \draw[thick] ($(infty) + (0,2.5pt)$) -- ($(infty) - (0,2.5pt)$);
            \tkzLabelPoint[above](infty){$\infty$}

            \tkzDefPoint(0:{-5/12*\len}){p1}
            \draw[thick] ($(p1) + (0,2.5pt)$) -- ($(p1) - (0,2.5pt)$);

            \tkzDefPoint(0:{-2/12*\len}){p2}
            \draw[thick] ($(p2) + (0,2.5pt)$) -- ($(p2) - (0,2.5pt)$);

            \path (p1) -- (p2) node[midway, below, xshift=1pt, yshift=2pt] {$\cdots$};

            \draw[decorate,decoration={brace,amplitude=3pt,mirror}] 
            ($(p1)-(0,0.3)$) -- ($(p2)-(0,0.3)$) node[midway,below=1.5pt,font=\scriptsize] {$r_1$};

            \tkzDefPoint(90:{-5/12*\len}){psi}
            \fill (psi) circle (2pt);
            \tkzLabelPoint[right](psi){$\psi^{d-1}$}

            \tkzDefPoint(90:{-7/24*\len}){p3}
            \draw[thick] ($(p3) + (2.5pt,0)$) -- ($(p3) - (2.5pt,0)$);

            \tkzDefPoint(90:{-1/12*\len}){p4}
            \draw[thick] ($(p4) + (2.5pt,0)$) -- ($(p4) - (2.5pt,0)$);

            \path (p3) -- (p4) node[midway, right, yshift=3pt] {$\vdots$};

            \draw[decorate,decoration={brace,amplitude=3pt,mirror}] 
            ($(p3)+(0.35,0)$) -- ($(p4)+(0.35,0)$)      node[midway,right=1pt,font=\scriptsize] {$r_2$};
        \end{tikzpicture}
        \right].
    \end{equation*}
All divisors in the above summations are isomorphic, therefore the corresponding invariants are $\frac{(r_1+r_2)!}{r_1!r_2!}$ times the invariant on a single such irreducible divisor. For every $\sigma\in S_{\mu,d}^{(h)}$, the degeneration formula ensures that when pulling back these classes along the branch morphism $br:\adm_\sigma^\bullet\rightarrow\overline{\mathcal{M}}_{\{0,\bullet,p_1,\dots,p_r,\infty\}}$ and intersecting with $\lambda_h$, one obtains a product of two intersection numbers on admissible cover spaces determined by the data on each twig, where the curves are glued together over the node to form a nodal admissible cover with data corresponding to $\sigma$. The class $\lambda_h$ distributes itself over each twig according to the (geometric) genus of the associated admissible cover space, while the class $\psi_\bullet^{d-1}$ only remains on the twig with $\bullet$. We have
    \begin{align*}
        \int_{\adm_\sigma^\bullet}\lambda_h\,br^*\bigl(D(0;r_1 \,\vert\, \bullet,\infty;r_2)\psi_\bullet^{d-1}\bigr) & = 0 \quad \text{whenever } \; r_1\ge3, \\
        \int_{\adm_\sigma^\bullet}\lambda_h\,br^*\bigl(D(0,\infty;0\,\vert\, \bullet;d+h)\psi_\bullet^{d-1}\bigr) & = 0.
    \end{align*}
    The first line is due to dimension reasons on the twig containing $0$. The second line is because, by the Riemann-Hurwitz formula, the total ramification degree over the node is zero (i.e. the ramification type is $(1,1,1,\dots)$), so over the twig containing $\bullet$ there would be at least one curve of positive genus admitting a degree $1$ map to $\mathbb{P}^1$. Therefore, with the choice made in \eqref{eq:psi_choice}, we have 
    \begin{align}
        \label{eq:psi_bullet_1} \psi_\bullet^d = & \; \frac{1}{d}D(0;1\,\vert\,\bullet,\infty;d+h-1)\psi_\bullet^{d-1} + \frac{2}{d}D(0;2\,\vert\,\bullet,\infty;d+h-2)\psi_\bullet^{d-1} \\
        \label{eq:psi_bullet_2} & - \sum_{g=1}^{h-1}\frac{g}{d}D(0,\infty;h-g\,\vert\,\bullet,d+g)\psi_\bullet^{d-1},
    \end{align}
    modulo terms that vanish when pulled back along $br$ and intersected with $\lambda_h$. First, we consider the part of $\Tr{\mu,d}^\bullet$ in equation \eqref{eq:general_Tmud} corresponding to the $\frac{2}{d}D(0;2\,\vert\,\bullet,\infty;d+h-2)$ term in \eqref{eq:psi_bullet_1}, which involves intersection numbers
    \begin{equation*}
        \frac{1}{\prod_i\lvert\operatorname{Aut}(\mu_i)\rvert}\int_{\adm_\sigma^\bullet}\lambda_h\,br^*
        \left[
        \begin{tikzpicture}[baseline=-3ex, scale=0.7,  transform shape]
            \tkzDefPoint(0,0){O} 
            \def\len{3} 
            \tkzDefPoint(35:{1/6*\len}){A}
            \tkzDefPoint(35:{-3/4*\len}){B}
            \tkzDefPoint(-35:{3/4*\len}){C}
            \tkzDefPoint(-35:{-1/6*\len}){D}
        
            \tkzDrawSegment[thick](B,A)
            \tkzDrawSegment[thick](D,C)
        
            \tkzDefPoint(35:{-5/8*\len}){zero}
            \draw[thick] ($(zero) + (0,2.5pt)$) -- ($(zero) - (0,2.5pt)$);
            \tkzLabelPoint[below](zero){$0$}

            \tkzDefPoint(35:{-5.5/12*\len}){p1}
            \draw[thick] ($(p1) + (0,2.5pt)$) -- ($(p1) - (0,2.5pt)$);
            
            \tkzDefPoint(35:{-3/12*\len}){p2}
            \draw[thick] ($(p2) + (0,2.5pt)$) -- ($(p2) - (0,2.5pt)$);

            \tkzDefPoint(-35:{4/8*\len}){psi}
            \fill (psi) circle (2pt);
            \tkzLabelPoint[above right=-2pt](psi){$\psi^{d-1}$}

            \tkzDefPoint(-35:{5/8*\len}){infty}
            \draw[thick] ($(infty) + (0,2.5pt)$) -- ($(infty) - (0,2.5pt)$);
            \tkzLabelPoint[below=1pt](infty){$\infty$}
            
            \tkzDefPoint(-35:{8/24*\len}){p3}
            \draw[thick] ($(p3) + (0,2.5pt)$) -- ($(p3) - (0,2.5pt)$);

            \tkzDefPoint(-35:{3/24*\len}){p4}
            \draw[thick] ($(p4) + (0,2.5pt)$) -- ($(p4) - (0,2.5pt)$);

            \path (p3) -- (p4) node[midway, below, xshift=-1pt, yshift = 12pt] {$\ddots$};

            \draw[decorate,decoration={brace,amplitude=3.5pt}] 
            ($(p3)-(0.1,0.25)$) -- ($(p4)-(0.1,0.25)$) node[midway,below=3pt,font=\scriptsize] {$d+h-2$};
        \end{tikzpicture}
        \right], \quad \sigma\in S_{\mu,d}^{(h)}, \; 0\le h\le \lvert\mu\rvert+d.
    \end{equation*}
    Since the dimension of the admissible cover space associated with the twig on the left is $1$, the corresponding intersection number must be $\int_{\adm_1(\delta,\delta)}\lambda_1=\frac{\delta^2-1}{12}$ for some $\delta>0$. So for $\sigma=(\mu_i,g_i,n_i,\ain_i)_{i=1}^m$, the intersection number above can be obtained by gluing a two-valent genus 1 graph to a graph of type $(\mu_i,g_i-1,n_i,\ain_i)$ for some $i$, as shown in Figure \ref{fig:glueciambella}.
    \begin{figure}
        \centering
      \begin{tikzpicture}[thick, font=\small, scale = 0.7 , transform shape]
            \node[draw, circle, minimum size=0.8cm] (v1) at (0,0) {$1$};
            \node[above left, xshift=-2mm, yshift=-2mm] at (v1.north) {$\lambda_1$};
            \draw (v1.105) -- ++(-0.115, 0.25);
            \draw (v1.75) -- ++(0.115, 0.25);
            \draw (v1.180) -- ++(-1, 0) node[left] {$\delta_i$};
            \draw (v1.0) -- ++(1, 0) node[right] (d1) {$\delta_i$};
            \begin{scope}[shift={(5.5, 0)}]
                \node[draw, circle, minimum size=1.2cm] (g1) at (0,0) {$g_i-1$};
                \node[above left] at (g1.130) {$\lambda_{g_i-1}$};
                \draw (g1.180) -- ++(-1.2, 0) node[left] (d2) {$\delta_i$};
                \draw (g1.105) -- ++(-0.175, 0.35) coordinate (t1_1);
                \draw (g1.75) -- ++(0.175, 0.35) coordinate (t2_1);
                \draw[dotted, thick] (-0.15, 1) to[bend left=10] (0.22, 1);
                \draw[decorate, decoration={brace, amplitude=4pt, raise=2pt}] 
                (t1_1) -- (t2_1) node[midway, above=6pt] {$r_i-2$};
                \draw (g1.25) -- ++(1, 0.5) node[right] {$a_{i,1}$};
                \draw (g1.-25) -- ++(1, -0.5) node[right] {$a_{i,n_i}$};
                \draw[loosely dotted, very thick] (1.5, 0.5) to[bend left=15] (1.5, -0.45);
                \draw (g1.-120) -- ++(-0.85, -1.05) node[below, xshift=-2mm] {$\mu_{i,1}+1$};
                \draw (g1.-90) -- ++(0.15, -1.2) node[below] {$\mu_{i,\ell_i}+1$};
                \draw (g1.-65) -- ++(0.7, -1.05) node[below, xshift=8mm] {$1^{\delta_i-\lvert\mu_i\rvert-\ell_i}$};
                \draw[loosely dotted, very thick] (-0.85, -1.65) to[bend right=15] (-0.1, -1.8);
                \node at (-0.375, -1.3) {$\psi^{d-1}$};
            \end{scope}
            \draw[dashed, thick] (d1) -- (d2);
            \node at (8.5, 0) {\scalebox{1.5}{$\leadsto$}};
            \begin{scope}[shift={(11.5, 0)}]
                \node[draw, circle, minimum size=1.2cm] (g2) at (0,0) {$g_i$};
                \node[above left] at (g2.130) {$\lambda_{g_i}$};
                \draw (g2.180) -- ++(-1.2, 0) node[left] {$\delta_i$};
                \draw (g2.105) -- ++(-0.175, 0.35) coordinate (t1_2);
                \draw (g2.75) -- ++(0.175, 0.35) coordinate (t2_2);
                \draw[dotted, thick] (-0.15, 0.9) to[bend left=10] (0.22, 0.9);
                \draw[decorate, decoration={brace, amplitude=4pt, raise=2pt}] 
                (t1_2) -- (t2_2) node[midway, above=6pt] {$r_i$};
                \draw (g2.25) -- ++(1, 0.5) node[right] {$a_{i,1}$};
                \draw (g2.-25) -- ++(1, -0.5) node[right] {$a_{i,n_i}$};
                \draw[loosely dotted, very thick] (1.5, 0.5) to[bend left=15] (1.5, -0.45);
                \draw (g2.-120) -- ++(-0.75, -1.05) node[below, xshift=-3mm] {$\mu_{i,1}+1$};
                \draw (g2.-90) -- ++(0, -1.2) node[below, xshift=-1mm] {$\mu_{i,\ell_i}+1$};
                \draw (g2.-65) -- ++(0.6, -1.05) node[below, xshift=7mm] {$1^{\delta_i-\lvert\mu_i\rvert-\ell_i}$};
                \draw[loosely dotted, very thick] (-0.75, -1.4) to[bend right=15] (-0.15, -1.55);
                \node at (-0.34, -1.1) {$\psi^d$};
            \end{scope}
        \end{tikzpicture}
        \caption{Attaching a genus one vertex.}
        \label{fig:glueciambella}
    \end{figure}
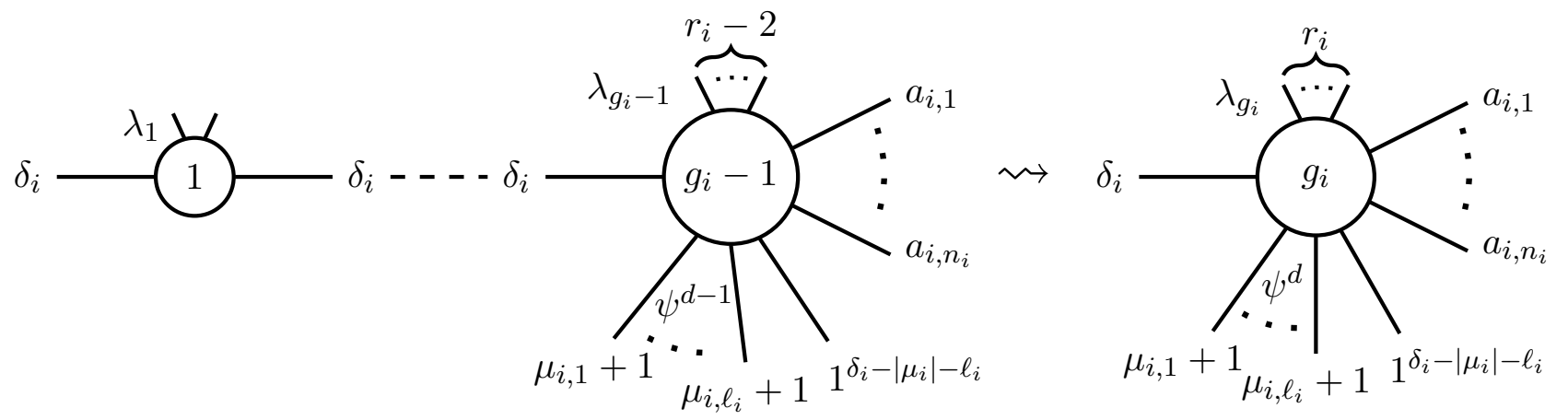

    There is an additional factor of $\binom{r_i}{2}$ coming from the number of ways to distribute two simple branch points on the left twig. Alternatively, if $(\mu_i,g_i,n_i,\ain_i)=(\emptyset,1,1,\delta_i)$ for some $i$, then we are gluing the genus one vertex above to a trivial cover, and the corresponding intersection number above is obtained by multiplying $\frac{\delta_i^2-1}{12}$ by $\int_{\adm_{\sigma^\prime}^\bullet}\lambda_h\psi_2^{d-1}$, where $\sigma^\prime = (\mu_j,g_j,n_j,\ain_j)_{j=1,j\neq i}^m$, together with a factor of $\binom{d+h}{2}$ from the simple branch points. This  implies that the part of $\Tr{\mu,d}^\bullet$ in \eqref{eq:general_Tmud} coming from  $\frac{2}{d}D(0;2\,\vert\,\bullet,\infty;d+h-2)$ is
    \begin{equation*}
        \frac{v}{d}\left(\e^2v\sum_{a\in\mathbb{N}}\frac{a^2-1}{12}p_aq_a\right) \star \Tr{\mu,d-1}^\bullet.
    \end{equation*}
    A similar argument with a tripod instead of a genus one vertex shows that the part associated to $\frac{1}{d}D(0;1\,\vert\,\bullet,\infty;d+h-1)$ is
    \begin{equation*}
        \frac{v}{d}\left(\frac{1}{\h}\sum_{\ain\in\mathbb{N}^2}\frac{p_\ain}{2!}q_{\lvert\ain\rvert}\right) \star \Tr{\mu,d-1}^\bullet.
    \end{equation*}
    Therefore the part of $\Tr{\mu,d}^\bullet$ in \eqref{eq:general_Tmud} corresponding to  \eqref{eq:psi_bullet_1}  is
    \begin{equation*}
        \frac{v}{d}\Trz{1}\star\Tr{\mu,d-1}^\bullet.
    \end{equation*}
    Consider the part of $\Tr{\mu,d}^\bullet$ associated to  \eqref{eq:psi_bullet_2}: 
    \begin{align*}
        \sum_{h=0}^{\lvert\mu\rvert+d}\sum_{\sigma\in S_{\mu,d}^{(h)}}&\frac{1}{m!}\h^{h-m}\e^{2h}\frac{v^{d+h}}{(d+h)!}\binom{d+h}{r_1,\dots,r_m}\frac{1}{\prod_i\lvert\operatorname{Aut}(\mu_i)\rvert} \\
        & \sum_{g=1}^{h-1}\frac{g}{d}\int_{\adm_\sigma^\bullet}\lambda_h\,br^*\bigl(D(0,\infty;h-g\,\vert\,\bullet,d+g)\psi_\bullet^{d-1}\bigr)\prod_{i=1}^m\frac{p_{\ain_i}}{n_i!}q_{\lvert\ain_i\rvert}.
    \end{align*}
    By changing the order of summation with respect to $h$ and $g$, and re-indexing $h$, one obtains
    \begin{align}
        \label{eq:Tcorrection} \sum_{g=1}^{\lvert\mu\rvert+d-1}\h^g\e^{2g}\frac{v^{d+g}}{(d+g)!}\frac{g}{d} & \sum_{h=0}^{\lvert\mu\rvert+d-g}\sum_{\sigma\in S_{\mu,d}^{(g+h)}}\frac{1}{m!}\h^{h-m}\e^{2h}\frac{v^h}{h!}\frac{(d+g)!h!}{(d+g+h)!}\binom{d+g+h}{r_1,\dots,r_m} \\
        \nonumber & \frac{1}{\prod_i\lvert\operatorname{Aut}(\mu_i)\rvert}\int_{\adm_\sigma^\bullet}\lambda_{g+h} \, b^*\bigl(D(0,\infty;h\,\vert\,\bullet,d+g)\psi_\bullet^{d-1}\bigr)\prod_{i=1}^m\frac{p_{\ain_i}}{n_i!}q_{\lvert\ain_i\rvert}.
    \end{align}
    For every $\sigma = (\mu_i,k_i,n_i,\ain_i)_{i=1}^m\in S_{\mu,d}^{(g+h)}$, consider the intersection number
    \begin{equation}\label{eq:admsigmapicture}
        \frac{1}{\prod_i\lvert\operatorname{Aut}(\mu_i)\rvert}\int_{\adm_\sigma^\bullet}\lambda_{g+h} \, b^*
        \left[
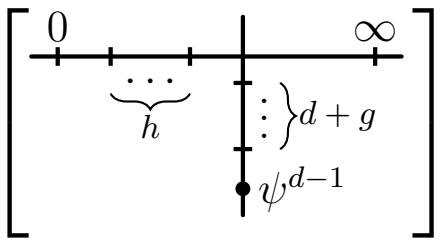
        \begin{tikzpicture}[baseline=-3ex , scale = 0.7 , transform shape]
            \tkzDefPoint(0,0){O} 
            \def\len{3} 
            \tkzDefPoint(0:{1/2*\len}){A}
            \tkzDefPoint(0:{-2/3*\len}){B}
            \tkzDefPoint(90:{1/8*\len}){C}
            \tkzDefPoint(90:{-1/2*\len}){D}
            
            \tkzDrawSegment[thick](B,A)
            \tkzDrawSegment[thick](D,C)
        
            \tkzDefPoint(0:{-7/12*\len}){zero}
            \draw[thick] ($(zero) + (0,2.5pt)$) -- ($(zero) - (0,2.5pt)$);
            \tkzLabelPoint[above](zero){$0$}

            \tkzDefPoint(0:{5/12*\len}){infty}
            \draw[thick] ($(infty) + (0,2.5pt)$) -- ($(infty) - (0,2.5pt)$);
            \tkzLabelPoint[above](infty){$\infty$}

            \tkzDefPoint(0:{-5/12*\len}){p1}
            \draw[thick] ($(p1) + (0,2.5pt)$) -- ($(p1) - (0,2.5pt)$);

            \tkzDefPoint(0:{-2/12*\len}){p2}
            \draw[thick] ($(p2) + (0,2.5pt)$) -- ($(p2) - (0,2.5pt)$);

            \path (p1) -- (p2) node[midway, below, xshift=1pt] {$\cdots$};

            \draw[decorate,decoration={brace,amplitude=3pt,mirror}] 
            ($(p1)-(0,0.35)$) -- ($(p2)-(0,0.35)$) node[midway,below=1.5pt,font=\scriptsize] {$h$};

            \tkzDefPoint(90:{-5/12*\len}){psi}
            \fill (psi) circle (2pt);
            \tkzLabelPoint[right](psi){$\psi^{d-1}$}

            \tkzDefPoint(90:{-7/24*\len}){p3}
            \draw[thick] ($(p3) + (2.5pt,0)$) -- ($(p3) - (2.5pt,0)$);

            \tkzDefPoint(90:{-1/12*\len}){p4}
            \draw[thick] ($(p4) + (2.5pt,0)$) -- ($(p4) - (2.5pt,0)$);

            \path (p3) -- (p4) node[midway, right, yshift=3pt] {$\vdots$};

            \draw[decorate,decoration={brace,amplitude=3pt,mirror}] 
            ($(p3)+(0.35,0)$) -- ($(p4)+(0.35,0)$)      node[midway,right=1pt,font=\scriptsize] {$d+g$};
        \end{tikzpicture}
        \right].
    \end{equation}
    The situation is rather analogous to \eqref{eq:righthandside}.
    Over the horizontal twig, the invariant  is
    \begin{equation*}
        \prod_{i=1}^m\frac{x_{\nu^\prime_i;h_i}(\ain_i)}{\lvert\operatorname{Aut}(\nu_i^\prime)\rvert},
    \end{equation*}
where $\nu_i^\prime$ is a partition of some length $\ell_i$, and the size of $\nu_i^\prime$ is $h_i-1+n_i$ by the Riemann-Hurwitz formula. For every $i=1,\dots,m$ and every $j=1,\dots,\ell_i$, there is a curve glued over the node to the point labeled by $\nu_{i,j}^\prime+1$ mapping to the vertical twig with degree $\nu_{i,j}^\prime+1$.

The ramification profile over $\bullet$ is
    \begin{equation*}
        (\rho_{i,j}+1) = \bigl(\rho_{i,j}^{(1)}+1,\dots,\rho_{i,j}^{(\ell_{i,j})}+1,1^{\nu^\prime_{i,j}+1-\lvert\rho_{i,j}\rvert-\ell_{i,j}}\bigr), \quad \ell_{i,j}=\ell(\rho_{i,j}),
    \end{equation*}
    for some possibly empty partition $\rho_{i,j}$, and there are $2g_{i,j}+\nu_{i,j}^\prime-\lvert\rho_{i,j}\rvert$ additional simple branch points. To match the data specified by $\sigma=(\mu_i,k_i,n_i,\mathbf{a}_i)_{i=1}^m$, it is necessary that $\sum_jg_{i,j}=k_i-h_i$ and $\coprod_{j}\rho_{i,j}=\mu_i$. Thus the invariant associated to the vertical twig is
    \begin{gather*}
        \frac{1}{\prod_{i,j}\lvert\operatorname{Aut}(\rho_{i,j})\rvert}\int_{\dadm_{\vec{g}}((\vec{\nu}^\prime+1),\,(\vec{\rho}+1))}\lambda_g\,\psi_\bullet^{d-1} \\
        = \prod_{i=1}^m\prod_{j=1}^{\ell_i}\frac{1}{(\nu^\prime_{i,j}+1)!}\frac{x_{\vec{\rho},d;\vec{g}}\bigl((1^{\nu^\prime_{1,1}+1}),\dots,(1^{\nu^\prime_{1,\ell_1}+1}),\dots\dots,(1^{\nu^\prime_{m,1}+1}),\dots,(1^{\nu^\prime_{m,\ell_m}+1})\bigr)}{\prod_{i,j}\lvert\operatorname{Aut}(\rho_{i,j})\rvert},
    \end{gather*}
    where $\vec{\rho}=(\rho_{i,j})_{1\le i\le m,1\le j\le\ell_i}$,  $\vec{g}=(g_{i,j})_{1\le i\le m,1\le j\le\ell_i}$ and $\vec{\nu}^\prime=(\nu_{i,j}^\prime)_{1\le i\le m,1\le j\le\ell_i}$, and we have used the string equation for disconnected invariants \eqref{eq:disc_string_1}, \eqref{eq:disc_string_2}. Including the gluing factors $\prod_{i,j}(\nu_{i,j}^\prime+1)$ and the ways to distribute the additional simple branch points on the various connected components, \eqref{eq:admsigmapicture} is equal to
    \begin{gather*}
        \sum_{\substack{h_1,\dots,h_m\ge0\\h_1+\cdots+h_m=h}} \, \sum_{\substack{\nu^\prime_1,\dots,\nu^\prime_m\\\nu^\prime_i\vdash h_i-1+n_i}}\prod_{i=1}^m\binom{r_i}{h_i,r_i-h_i}\frac{x_{\nu^\prime_i;h_i}(\ain_i)}{\lvert\operatorname{Aut}(\nu^\prime_i)\rvert} \cdot \prod_{i=1}^m\prod_{j=1}^{\ell_i}(\nu^\prime_{i,j}+1) \\
        \cdot \sum_{\substack{\rho_{1,1},\dots,\rho_{1,\ell_1}\\ \coprod_{j}\rho_{1,j}=\mu_1}} \cdots \sum_{\substack{\rho_{m,1},\dots,\rho_{m,\ell_m}\\ \coprod_j\rho_{m,j}=\mu_m}} \, \sum_{\substack{g_{1,1},\dots,g_{1,\ell_1}\ge0\\\sum_jg_{1,j}=k_1-h_1}} \cdots \sum_{\substack{g_{m,1},\dots,g_{m,\ell_m}\ge0\\\sum_jg_{m,j}=k_m-h_m}} \\
        \prod_{i=1}^m\binom{r_i-h_i}{2g_{i,1}+\nu^\prime_{i,1}-\lvert\rho_{i,1}\rvert,\dots,2g_{i,\ell_i}+\nu^\prime_{i,\ell_i}-\lvert\rho_{i,\ell_i}\rvert} \\
        \left(\prod_{i=1}^m\prod_{j=1}^{\ell_i}\frac{1}{(\nu^\prime_{i,j}+1)!}\right) \frac{x_{\vec{\rho},d;\vec{g}}\bigl((1^{\nu^\prime_{1,1}+1}),\dots,(1^{\nu^\prime_{1,\ell_1}+1}),\dots\dots,(1^{\nu^\prime_{m,1}+1}),\dots,(1^{\nu^\prime_{m,\ell_m}+1})\bigr)}{\prod_{i,j}\lvert\operatorname{Aut}(\rho_{i,j})\rvert},
    \end{gather*}
    For all partitions $\nu^\prime_1,\dots,\nu^\prime_m$ appearing in the sum above, let $\nu:=\nu^\prime_1\sqcup\dots\sqcup\nu^\prime_m$, and denote the individual elements of $\nu$ by $\nu_1,\dots,\nu_\ell$, where $\ell:=\ell(\nu)$. The total ramification degree over the node in \eqref{eq:admsigmapicture} is independent of $h$:
    \begin{equation*}
        \lvert\nu\rvert = \sum_{i=1}^m\lvert\nu^\prime_i\rvert = h-m+n = \lvert\mu\rvert + d-g,
    \end{equation*}
    Summing over all $0\le h\le\lvert\mu\rvert+d-g$ and $\sigma\in S_{\mu,d}^{(g+h)}$ in \eqref{eq:Tcorrection}, we obtain
    \begin{align*}
        & \sum_{g=1}^{\lvert\mu\rvert+d-1}\h^g\e^{2g}\frac{v^{d+g}}{(d+g)!}\frac{g}{d} \sum_{\nu\vdash\lvert\mu\rvert+d-g}\Biggl(\sum_{m\ge1}\frac{1}{m!}\sum_{\substack{\nu^\prime_1,\dots,\nu^\prime_m\\ \coprod_i\nu^\prime_i=\nu}}\sum_{h_1,\dots,h_m\ge0}\h^{h-m}\e^{2h}\frac{v^h}{h!}\binom{h}{h_1,\dots,h_m} \\
        & \hspace{13em}\sum_{\ain_1\in\mathbb{N}^{\lvert\nu^\prime_1\rvert+1-h_1}}\cdots\sum_{\ain_m\in\mathbb{N}^{\lvert\nu^\prime_m\rvert+1-h_m}} \prod_{i=1}^m\frac{x_{\nu^\prime_i;h_i}(\ain_i)}{\lvert\operatorname{Aut}(\nu^\prime_i)\rvert}\frac{p_{\ain_i}}{(\lvert\nu^\prime_i\rvert+1-h_i)!}{q_{\lvert\ain_i\rvert}}\Biggr) \\
        & \cdot \prod_{j=1}^{\ell}\frac{1}{\nu_j!}\sum_{\substack{\rho_1,\dots,\rho_\ell\\ \coprod_j\rho_j=\mu}}\sum_{\substack{g_1,\dots,g_\ell\ge0\\g_1+\cdots+g_\ell=g}}\binom{d+g}{2g_1+\nu_1-\lvert\rho_1\rvert,\dots,2g_\ell+\nu_\ell-\lvert\rho_\ell\rvert} \frac{x^{}_{\vec{\rho},d;\vec{g}}\left((1^{\nu_1+1}),\dots,(1^{\nu_\ell+1})\right)}{\prod_j\lvert\operatorname{Aut}(\rho_j)\rvert} \\
        & = \sum_{g=1}^{\lvert\mu\rvert+d-1}\h^g\e^{2g}\frac{v^{d+g}}{(d+g)!}\frac{g}{d} \sum_{\nu\vdash\lvert\mu\rvert+d-g} \Tr{\nu}^\bullet \, \gamma_{\mu,d-1}^{g,\nu},
    \end{align*}
    where $\gamma_{\mu,d-1}^{g,\nu}$ is the quantity in the second to last line. One sees that $\gamma_{\mu,d-1}^{g,\nu}$ is the coefficient of $\h^{g-\ell(\nu)}\e^{2g}\frac{v^{d+g-1}}{(d+g-1)!}\frac{1}{\lvert\operatorname{Aut}(\nu)\rvert}\prod_{i=1}^{\ell(\nu)}\frac{p_1^{\nu_i+1}}{\nu_i+1}q_{\nu_i+1}$ in $\left(p_1^2\frac{\partial}{\partial p_2} + \frac{1}{\h}\frac{p_1^2}{2}q_2\right)\Tr{\mu,d-1}^\bullet$ by Lemma \ref{lem:Tmu_string}. What is left to prove is that the $g\ge d+1$ terms vanish, and that the only contribution to the $g=d$ term is from $\gamma_{\mu,d-1}^{d,\mu}$. This is equivalent to the following two statements:
    \begin{align*}
        \lvert\nu\rvert < \lvert\mu\rvert & \implies \gamma_{\mu,d-1}^{g,\nu} = 0; \\
        \lvert\nu\rvert = \lvert\mu\rvert \text{ and } \nu \neq \mu & \implies \gamma_{\mu,d-1}^{g,\nu} =0.
    \end{align*}
    Both are derived from the fact that, for any choice of partitions $\coprod_{j=1}^\ell\rho_j=\mu$, the intersection number $x^{}_{\vec{\rho},d;\vec{g}}\left((1^{\nu_1+1}),\dots,(1^{\nu_\ell+1})\right)$ vanishes unless $\lvert\rho_j\rvert+\ell(\rho_j)\le\nu_j+1$ for every $j=1,\dots,\ell$.
\end{proof}
We write down the explicit formulas \eqref{eq:TmudTmuprop} for $d+\lvert\mu\rvert\le 5$ in Appendix \ref{appsec:TmudTmuprop}.

Combining Proposition \ref{prop:Tds} and Proposition \ref{prop:TmuTk}, we obtain that all invariants  $x_{\vec{\mu},d;\vec{g}}(\vec{\ain})$, can be expressed in terms of the building blocks $x_{k;g}(\ain)$.
\begin{thm}\label{thm:Tmud_in_Tk}
    For every partition $\mu=(\mu_1,\dots,\mu_\ell)$ and integer $d\ge0$ such that $(\ell,d)\neq(1,0)$, the generating function $\Tr{\mu,d}^\bullet$ can be expressed in terms of $\star$-products of $\Trz{1},\dots,\Trz{\lvert\mu\rvert+d-1}$ as
    \begin{equation*}
        \Tr{\mu,d}^\bullet = \frac{1}{\lvert\operatorname{Aut}(\mu)\rvert}\frac{v^d}{d!}\Trz{1}^{\star d}\star\Trz{\mu_1}\star\cdots\star\Trz{\mu_\ell} + \sum_{g=0}^{\lvert\mu\rvert+d-\mu_1}\h^g\e^{2g}\frac{v^{d+g}}{(d+g)!}(\Tr{\mu,d}^\bullet)^{(g)},
    \end{equation*}
    where $(\Tr{\mu,d}^\bullet)^{(g)}$ is a linear combination of terms $\Trz{k_1}\star\cdots\star\Trz{k_m}$ with $k_1+\cdots+k_m=\lvert\mu\rvert+d-g$ and $m\le d+\ell-1$.
\end{thm}
For all $d+\lvert\mu\rvert\le 5$, we write the explicit formulas for $\Tr{\mu,d}^\bullet$ in terms of $\Trz{1},\Trz{2},\Trz{3},\Trz{4}$ in Appendix \ref{appsec:Tmud_in_Tk}.
\begin{remark}
    It was mentioned in Remark \ref{rem:Tk1k2_comm} that the generating series for the building blocks $x_{k;g}(\ain)$ all $\star$-commute: $[\Trz{k_1},\Trz{k_2}]=\Trz{k_1}\star\Trz{k_2} - \Trz{k_2}\star\Trz{k_1}=0$. Therefore, an obvious consequence of Theorem \ref{thm:Tmud_in_Tk} is
    \begin{equation*}
        [\Tr{\mu_1,d_1}^\bullet,\Tr{\mu_2,d_2}^\bullet] = 0 \quad \forall\,\mu_1,\mu_2,d_1,d_2.
    \end{equation*}
    Geometrically, this is the same as pulling back the WDVV relation via the branch morphism to appropriate admissible cover spaces with four ramification profiles and intersecting with $\lambda_g$ and two (powers of) $\psi$-classes.
\end{remark}

\section{The KdV hierarchy and polynomiality}

\subsection{Admissible covers, stable maps, and the KdV hierarchy}
In \cite{bur15}, Buryak presented a construction of a Hamiltonian hierarchy of PDEs associated to any Cohomological Field Theory, called the double ramification (DR) hierarchy. Its construction is motivated by ideas in Symplectic Field Theory \cite{Intro2SFT}. Its approach differs from the one in \cite{dz01} and is based on integrating classes over the DR cycle on the moduli space of stable curves (although the DR and DZ constructions have been shown to be equivalent when the CohFT is semisimple \cite{bls24}). A description of the DR cycle is via the pushforward of the virtual fundamental class of the moduli space of stable maps to rubber $\mathbb{P}^1$, and Pixton's formula \cite{jppz17} gives an expression for the cycle as a sum over stable graphs. One of the authors' original questions that motivated this paper is whether there is an interesting relationship between the DR hierarchy and the analogous constructions made with admissible covers. \\

An important property of the DR cycle $\operatorname{DR}_g(a_1,\dots,a_n)$, which is crucial to the DR hierarchy construction of Buryak, is that it is polynomial in its entries \cite{pix23, spe24}. This is not true for the corresponding admissible covers cycle, thereby making it more difficult to carry out such a construction. This, among other considerations, is what motivated the authors to first consider the case where all admissible covers were fully ramified over a point. In such a restricted setting, all of the invariants that we consider turn out to be polynomial (see Theorem \ref{thm:disc_poly} later). \\

The first example of the DR hierarchy construction that was worked out is when the associated CohFT is trivial. In this case, the Hamiltonians of the DR hierarchy take the following form.
\begin{align*}
    \overline{g}_k & = \sum_{g\ge0,n\ge2}\frac{(-\e^2)^g}{n!}\sum_{\substack{a_1,\dots,a_n\in\mathbb{Z}\\a_1+\cdots+a_n=0}}\left(\int_{\operatorname{DR_g}(0,a_1,\dots,a_n)}\lambda_g\psi_1^k\right)p_{a_1}\cdots p_{a_n} \\
    & = \sum_{g=0}^k(-\e^2)^g\sum_{\substack{a_1,\dots,a_{k+2-g}\in\mathbb{Z}\\ a_1+\cdots+a_{k+2-g}=0}}\left(\int_{\operatorname{DR}_g(0,a_1,\dots,a_{k+2-g})}\lambda_g\psi_1^k\right)\frac{p_{a_1}\cdots p_{a_{k+2-g}}}{(k+2-g)!}, \quad k\ge1.
\end{align*}
Here the $p_*$ variables play a similar role to the one in this paper. To be precise, the variables $\{p_a:a>0\}$ agree, while the $\{p_{-a}:a>0\}$ correspond to the $\{q_a:a>0\}$ in previous sections. Buryak showed that these Hamiltonians correspond to the well-known KdV hierarchy of PDEs. This is  a version of the Witten-Kontsevich theorem \cite{wit90, kon92}, which relates intersections of $\psi$-classes on the moduli space of stable curves to the KdV hierarchy. Usually, the KdV equations are written in terms of a single function $u(x)$, rather than the variables $p_*$. They are related by the Fourier transformation
\begin{equation*}
    u(x) = \sum_{a\in\mathbb{Z}}p_ae^{iax}.
\end{equation*}
The variational derivative $\frac{\delta}{\delta u}=\sum_{n\ge0}(-\partial_x)^n\frac{\partial}{\partial u_n}$, where $u_n:=\partial^n_xu$, frequently appears. In Fourier coordinates it is given by $\frac{\delta}{\delta u}=\sum_{a\in\mathbb{Z}}e^{-iax}\frac{\partial}{\partial p_a}$. The variational derivative of the Hamiltonians above is
\begin{equation}\label{eq:hk_kdv}
    h_k := \frac{\delta\overline{g}_k}{\delta u} = \sum_{g=0}^k(-\e^2)^g\sum_{\ain\in\mathbb{Z}^{k+1-g}}\left(\int_{\operatorname{DR}_g(-\sum_ia_i,0,a_1,\dots,a_{k+1-g})}\hspace{-2em}\lambda_g\psi_2^k\right)\frac{p_{a_1}\cdots p_{a_{k+1-g}}}{(k+1-g)!}e^{i\sum_ia_ix}.
\end{equation}
The question of whether the DR hierarchy is endowed with a bihamiltonian structure is a difficult one. Recent results \cite{br_bihamiltonian} show that this is the case when the CohFT is semisimple and homogeneous. The existence of a bihamiltonian structure is useful, since it gives rise to powerful recursions for the corresponding Hamiltonians. In the case of the DR hierarchy associated to the trivial CohFT, this recursion corresponds to the Lenard-Magri \cite{magri_recursion} recursion for the KdV hierarchy, which can be written in terms of the objects $h_k$ as
\begin{equation}\label{eq:lenard_magri}
    \partial_xh_{k+1} = \frac{2}{2k+3}\left(u\partial_x+\frac12u_x+\frac{\e^2}{8}\partial_x^3\right)h_k, \quad k\ge1.
\end{equation}
Using this, we prove the following  lemma.
\begin{lem}\label{lem:DRsinh}
    For every $d\ge1$ and $g\ge0$,
    \begin{equation}\label{eq:DR_1d}
        \int_{\operatorname{DR}_g(-d,0,1^d)}\lambda_g\psi_2^{d+g-1} = \frac{(d+g)!}{2^gd!}[t^{2g}]\left(\frac{\sinh(t/2)}{t/2}\right)^{2d},
    \end{equation}
    where $1^d$ denotes the vector with $d$ entries equal to $1$, and $[t^{2g}]$ denotes the coefficient of $t^{2g}$ in the power series.
\end{lem}
\begin{proof}
    For $g=0$ it is immediate: $\int_{\overline{\mathcal{M}}_{0,d+2}}\psi_2^{d-1}=1$. Let $X_{d,g}$ denote the left-hand side of \eqref{eq:DR_1d}. This is the coefficient of $(-\e^2)^g\frac{p_1^d}{d!}e^{idx}$ in $h_{d+g-1}$. Setting $k=d+g-2$ in the Lenard-Magri recursion \eqref{eq:lenard_magri}, extracting the coefficient of $(-\e^2)^g\frac{p_1^d}{d!}e^{idx}$ and dividing by $id$ yields
    \begin{equation}\label{Xdg_recursion}
        X_{d,g} = \frac{1}{2d+2g-1}\left((2d-1)X_{d-1,g} +  \frac{d^2}{4}X_{d,g-1}\right),
    \end{equation}
    where $X_{0,g}$ is taken to be zero. To show that the right-hand side of \eqref{eq:DR_1d} satisfies this recursion, notice that $\frac{d^2}{dt^2}\sinh^{2d}(t/2)=\frac{d(2d-1)}{2}\sinh^{2d-2}(t)+d^2\sinh^{2d}(t)$. The coefficient of $t^{2d+2g-2}$ of this differential equation gives
    \begin{align*}
        \frac{(2d+2g)(2d+2g-1)}{2^{2d}}[t^{2g}]\left(\frac{\sinh(t/2)}{t/2}\right)^{2d} = & \; \frac{d(2d-1)}{2^{2d-1}}[t^{2g}]\left(\frac{\sinh(t/2)}{t/2}\right)^{2(d-1)} \\
        & + \frac{d^2}{2^{2d}}[t^{2g-2}]\left(\frac{\sinh(t/2)}{t/2}\right)^{2d},
    \end{align*}
    which gives the desired result after multiplying by $\frac{2^{2d-1}(d+g-1)!}{2^gd!(2d+2g-1)}$.
\end{proof}
Another way to write the intersection number \eqref{eq:DR_1d} is
\begin{align}
    \int_{\operatorname{DR}_g(-d,0,1^d)}\lambda_g\psi_2^{d+g-1} & = \frac{(2d-1)!!}{2^{2g}(2d+2g-1)!!}[t^{2g}]\prod_{n=1}^d\frac{1}{1-n^2t^2} \nonumber \\
    & = \frac{(2d-1)!!}{2^{2g}(2d+2g-1)!!}h_g(1^2,\dots,d^2), \label{eq:before}
\end{align}
where $h_g$ is the complete homogeneous symmetric polynomial of degree $g$ in $d$ variables. This is because the recursion \eqref{Xdg_recursion} is equivalent to the relation $(1-d^2t^2)\prod_{n=1}^d\frac{1}{1-n^2t^2}=\prod_{n=1}^{d-1}\frac{1}{1-n^2t^2}$. \newline

In \cite[Theorem 2.7]{nss2025}, a wall-crossing formula between Gromov-Witten cycles and Hurwitz cycles is proven:
\begin{equation}\label{eq:nss_wall}
    \mathsf{GW}_{g,\underline{\mu}} = \sum_{\Gamma}\frac{1}{\lvert\operatorname{Aut}(\Gamma)\rvert}\cdot\operatorname{gl}_{\Gamma*}\left(\mathsf{H}_{\underline{\mu},\underline{\eta}}\boxtimes\prod_{i=1}^kI_{g_i,n_i,\eta^i}(-\tilde{\psi}_i)\right).
\end{equation}
We give a brief overview of the  notation from \cite{nss2025} which is within the scope of our discussion. Let $\mathcal{X}\rightarrow\overline{M}_{0,m}^\circ$ be the universal curve over the moduli space of stable connected curves of genus $0$ with $m$ markings, and let $\underline{\mu}=(\mu^1,\dots,\mu^m)$ be a vector of partitions of an integer $\delta\ge1$ which specifies the ramification profiles of stable maps at the universal markings of $\mathcal{X}$. The Gromov-Witten cycle in \eqref{eq:nss_wall} is
\begin{equation*}
    \mathsf{GW}_{g,\underline{\mu}} = \frac{1}{\prod_i\lvert\operatorname{Aut}(\mu^i)\rvert}sr_*[\overline{M}_{g,n}(\mathcal{X},\underline{\mu})]^\text{vir} \in H^*(\overline{M}_{g,n+\sum_i\ell(\mu^i)}),
\end{equation*}
where $\overline{M}_{g,n}(\mathcal{X},\underline{\mu})$ is the $\overline{M}^\circ_{0,m}$-relative moduli space of possibly disconnected stable maps to $\mathcal{X}$ and $sr:\overline{M}_{g,n}(\mathcal{X},\underline{\mu})\rightarrow\overline{M}_{g,n+\sum_i\ell(\mu^i)}$ is the morphism that remembers the source curve. For another vector of partitions $\underline{\eta}=(\eta^1,\dots,\eta^k)$ of $\delta$, which specifies ramification profiles of admissible covers at additional $k$ moving markings, the Hurwitz cycle above is
\begin{equation*}
    \mathsf{H}_{\underline{\mu},\underline{\eta}} = \frac{1}{\prod_i\lvert\operatorname{Aut}(\mu^i)\rvert\prod_j\lvert\operatorname{Aut}(\eta^j)\rvert}sr_*[\overline{H}(\mathcal{X},\underline{\mu},\underline{\eta})]^\text{vir} \in H^*(\overline{M}_{g_0,\sum_j\ell(\mu^i)+\sum_j\ell(\eta^j)}),
\end{equation*}
where $\overline{H}(\mathcal{X},\underline{\mu},\underline{\eta})$ is the $\overline{M}^\circ_{g_0,m+k}$-relative moduli space of possibly disconnected admissible covers of $\mathcal{X}$, and $g_0$ is determined by the Riemann-Hurwitz formula. The sum in \eqref{eq:nss_wall} is over star shaped graphs $\Gamma$. Every graph $\Gamma$ has a root vertex with genus label $g_0$ corresponding to the unique admissible cover component of the maps, connected via an edge to $k$ other vertices with genus labels $g_1,\dots,g_k$, corresponding to contracted components. The $i$-th edge is labeled by a partition $\eta^i$ of $\delta$, and the $i$-th vertex has $n_i$ half-edges (called leafs) attached to it. The combinatorial data associated to graphs $\Gamma$ in the sum are subject to various conditions arising from the choice of $(g,n,\underline{\mu})$, see \cite[Section 2.5]{nss2025}. The $I$-functions $I_{g_i,n_i,\eta^i}(z)\in H^*(\overline{M}^\circ_{g_i,n_i+\ell(\eta^i)})[z^{\pm}]$ are given by an expression in $\lambda$ and $\psi$-classes arising from the virtual localization formula on $\mathbb{P}^1$ \cite{Graber-Pandharipande}. When $\delta\ge2$ and $\mathsf{GW}_{g,\underline{\mu}}$ is a cycle of dimension $g$, a consequence of \eqref{eq:nss_wall} is the  numerical equation
\begin{equation}\label{eq:GWHwithlag}   \lambda_g\mathsf{GW}_{g,\underline{\mu}} = \lambda_g\mathsf{H}_{\underline{\mu}}.
\end{equation}
Equation \eqref{eq:GWHwithlag} arises because intersecting with $\lambda_g$ forces $\eta^i=(\delta)$, and by dimension counting 
only the trivial graph in the sum contributes. A similar observation is made in \cite[Lemma 7.1]{br2025}. When we choose $\underline{\mu}$ to be $((\delta),(k+1),\ain)$ for $k\ge1$ and $\ain\in\mathbb{N}^{k+1-g}$, equation \eqref{eq:GWHwithlag} is equivalent to
\begin{equation}\label{eq:nss_us}
    \frac{1}{(\delta-k-1)!} \int_{[M_g(\mathcal{X},\underline{\mu})]^\text{vir}}sr^*(\lambda_g) = \frac{x_{k;g}(\ain)}{g!}.
\end{equation}
The $(\delta-k-1)!$ factor arises because on the right-hand side we are not marking the ramification points corresponding to the trivial part of $(k+1)=(k+1,1^{\delta-k-1})$. The $g!$ is because on the left-hand side we take the additional simple branch points to be unordered. In our context, the wall-crossing formula relates the building blocks $x_{k;g}(\ain)$ to integrals on the double ramification cycle. This is made precise in the proposition below.
\begin{prop}\label{prop:Adm2DR}
    For every $k\ge1$ and $0\le g\le k$ and $\ain=(a_1,\dots,a_{k+1-g})\in\mathbb{N}^{k+1-g}$,
    \begin{equation}\label{eq:admdr}
        \frac{x_{k;g}(\ain)}{k!g!} = \sum_{h=0}^g\frac{(k+1)!}{2^h(k+1-h)!}c_{k,h} \int_{\operatorname{DR}_{g-h}(-\sum_ia_i,0,a_1,\dots,a_{k+1-g})}\lambda_{g-h}\psi_2^{k-h},
    \end{equation}
    where $c_{k,h}$ is the coefficient of $t^{2h}$ in $\left(\frac{\operatorname{arcsinh}(t/2)}{t/2}\right)^{2(k+1-h)}$.
\end{prop}
\begin{proof}
    We prove the equivalent statement
    \begin{equation}\label{eq:reverse}
        \int_{\operatorname{DR}_g(-\sum_ia_i,0,a_1,\dots,a_{k+1-g})}\lambda_g\psi_2^k = \sum_{h=0}^g\frac{(k+1)!}{2^h(k+1-h)!}\tilde{c}_{k,h}\frac{x_{k-h;g-h}(\ain)}{(k-h)!(g-h)!},
    \end{equation}
    where $\tilde{c}_{k,h}$ is the coefficient of $t^{2h}$ in $\left(\frac{\sinh(t/2)}{t/2}\right)^{2(k+1-h)}$. In the notation from the previous discussion, the left-hand side of \eqref{eq:reverse} is
    \begin{equation*}
        \int_{[\overline{M}_{g,1}(\mathcal{X};(\delta),\ain)]^\text{vir}}sr^*(\lambda_g\psi_2^k),
    \end{equation*}
    where $\delta=\sum_{i=1}^{k+1-g}a_i$. We apply the tooth and comb degeneration \cite{op_trilogy1} to the target rational curve, thereby sending the point associated to the $\psi$-class onto the tooth. This gives rise to
    \begin{align*}
        \sum_{h=0}^g & \int_{[\overline{M}_{g-h}(\mathcal{X},\,\underline{\mu}^{(h)}]^\text{vir}}sr^*(\lambda_{g-h}) \cdot \frac{k+1-h}{(\delta-k-1+h)!} \\
        & \cdot \frac{1}{(k+1-h)!} \int_{[\overline{M}_{h,1}(\mathcal{X},(k+1-h),(1^{k+1-h})]^\text{vir}}sr^*(\lambda_h\psi_2^k) \\
        = \sum_{h=0}^g & \frac{1}{(k-h)!}\frac{1}{(\delta-k-1+h)!}\int_{[M_{g-h}(\mathcal{X},\underline{\mu}^{(h)})]^\text{vir}}sr^*(\lambda_{g-h}) \cdot \int_{\operatorname{DR}_h(-(k+1-h),0,1^{k+1-h})}\lambda_h\psi_2^k
    \end{align*}
    where $\underline{\mu}^{(h)}=((\delta),(k+1-h,1^{\delta-k-1+h}),\ain)$. The result follows from  \eqref{eq:nss_us} and Lemma \ref{lem:DRsinh}.
\end{proof}
Similarly to \eqref{eq:before}, one can write the coefficients $c_{k,h}$ appearing above as
\begin{align*}
    \frac{(k+1)!}{2^h(k+1-h)!}c_{k,h} & = \frac{(2k+1-2h)!!}{2^{2h}(2k+1)!!}[t^{2h}]\prod_{n=1}^k(1-n^2t^2) \\
    & = \frac{(2k+1-2h)!!}{2^{2h}(2k+1)!!}(-1)^he_h(1^2,\dots,k^2),
\end{align*}
where $e_h$ is the elementary symmetric polynomial of degree $h$ in $k$ variables.

The generating series $\Trz{k}$ for the invariants $x_{k;g}(\ain)$ are determined by the variational derivatives of the KdV hierarchy \eqref{eq:hk_kdv}, as a consequence of Proposition \ref{prop:Adm2DR}. To show this, we must tweak the definition of $\Trz{k}$ slightly:
\begin{equation*}
    \widetilde{\mathsf{T}}_k := \sum_{g=0}^k\e^{2g}\sum_{\ain\in\mathbb{Z}^{k+1-g}}\frac{x_{k;g}(\ain)}{k!g!}\frac{p_{\ain}}{(k+1-g)!}e^{i\lvert\ain\rvert x}.
\end{equation*}
Aside from dividing by $k!$ and setting $\h=v=1$, we have replaced $q_{\lvert\ain\rvert}$ with $e^{i\lvert\ain\rvert x}$, and we are summing over all integers instead of just positive integers. This makes sense due to the fact that the $x_{k;g}(\ain)$ are polynomial in the entries (see Corollary \ref{cor:polynomiality}).
\begin{cor}\label{cor:Tk_KdV}
    For every $k\ge1$,
    \begin{equation*}
        \widetilde{\mathsf{T}}_k = \sum_{g=0}^k\e^{2g}\frac{(k+1)!}{2^g(k+1-g)!}c_{k,g}h_{k-g}.
    \end{equation*}
\end{cor}

\subsection{Polynomiality}
We provide a proof of Theorem \ref{thm:poly}, i.e. we show that all the invariants $x_{\mu,d;g}(a_1,\dots,a_{\lvert\mu\rvert+d+1-g})$ are symmetric polynomials in their entries of degree $2g+d+\ell(\mu)-1$. We also give a more general statement involving the disconnected invariants $x_{\vec{\mu},d;\vec{g}}(\vec{\ain})$.

We start by discussing the building blocks $x_{k;g(\ain)}$. The DR cycle $\operatorname{DR}_g(a_1,\dots,a_n)$ is a polynomial of degree $2g$ in the entries $a_1,\dots,a_n$ \cite{pix23, spe24}. Therefore, a significant corollary of Proposition \ref{prop:Adm2DR} is the following.
\begin{cor}\label{cor:polynomiality}
    For every $k\ge1$ and $0\le g\le k$, the family of invariants $x_{k;g}(a_1,\dots,a_{k+1-g})$ is encoded as a symmetric polynomial in $a_{1},\dots,a_{k+1-g}$ of degree $2g$.
\end{cor}
\begin{remark}\label{rem:top_positive}
    The top degree monomials of $x_{k;g}(\ain)$ have positive coefficients. This follows from Proposition \ref{prop:Adm2DR}, the relationship between the DR integrals and the KdV hierarchy \cite{bur15}, and the positivity of the coefficients appearing the KdV hierarchy \cite[Theorem 1.4]{asKdV00}.
\end{remark}
Having established that the building blocks $x_{k;g}(\ain)$ are polynomial of degree $2g$, we unpack the statement of Theorem \ref{thm:Tmud_in_Tk} to deduce the result. The crucial step is showing that the $\star$-product preserves polynomiality when applied to the kind of generating series (i.e. the $\Trz{k}$) that we are considering. 
\begin{lem}\label{lem:starpol}
    For integers $n,m\ge1$, let $F\in\mathbb{Q}[x_1,\dots,x_n]$ and $G\in\mathbb{Q}[y_1,\dots,y_m]$ be symmetric polynomials of degrees $d$ and $e$, respectively. Let
    \begin{equation*}
        f(p_*,q_*) := \sumn{\ain}{n}F(a_1,\dots,a_n)\frac{p_{\ain}}{n!}q_{\lvert\ain\rvert}, \quad g(p_*,q_*) := \sumn{\mathbf{b}}{m}G(b_1,\dots,b_m)\frac{p_{\mathbf{b}}}{m!}q_{\lvert\mathbf{b}\rvert}.
    \end{equation*}
    Then
    \begin{equation*}
        f\odot g:= f\star g-fg = \h\sumn{\ain}{n+m-1}H(a_1,\dots,a_{n+m-1})\frac{p_{\ain}}{(n+m-1)!}q_{\lvert\ain\rvert},
    \end{equation*}
    where $H\in\mathbb{Q}[z_1,\dots,z_{n+m-1}]$ is a symmetric polynomial of degree $d+e+1$.
\end{lem}
\begin{proof}
    We have
    \begin{align*}
        \frac{1}{\h}(f\star g-fg) & = \sum_{k>0}k\frac{\partial f}{\partial p_k}\frac{\partial g}{\partial q_k} \\
        & = \sum_{k>0}k\sumn{\ain}{n-1}F(a_1,\dots,a_{n-1},k)\frac{p_{\ain}}{(n-1)!} q_{\lvert\ain\rvert+k}\sum_{\substack{\mathbf{b}\in\mathbb{N}^m\\\lvert\mathbf{b}\rvert=k}}G(b_1,\dots,b_m)\frac{\pb}{m!} \\
        & = \sumn{\ain}{n-1}\sumn{\mathbf{b}}{m}\lvert\mathbf{b}\rvert F\bigl(a_1,\dots,a_{n-1},\lvert\mathbf{b}\rvert\bigr)G(b_1,\dots,b_m)\frac{p_\ain p_{\mathbf{b}}}{(n-1)!\,m!}q_{\lvert\ain\rvert+\lvert\mathbf{b}\rvert} \\
        & = \sumn{\ain}{n+m-1}H(a_1,\dots,a_{n+m-1})\frac{p_\ain}{(n+m-1)!}q_{\lvert\ain\rvert}.
    \end{align*}
    The polynomial $H(a_1,\dots,a_{n+m-1})$ is given by
    \begin{align*}
        &= \frac{1}{(n-1)!\,m!}\sum_{\sigma\in S_{n+m-1}}\bigl(a_{\sigma(n)}+\cdots+a_{\sigma(n+m-1)}\bigr)F(a_{\sigma(1)},\dots,a_{\sigma(n-1)},a_{\sigma(n)}+\cdots+a_{\sigma(n+m-1)})\\
        & \hspace{10em} \cdot G(a_{\sigma(n)},\dots,a_{\sigma(n+m-1)}),
    \end{align*}
    where $S_{n+m-1}$ is the symmetric group on $n+m-1$ objects.
\end{proof}

The $\star$-product does not preserve polynomiality in general. It does in our context  because we have imposed full ramification over one of the branch points of the admissible cover, which ensures that the generating series $\Trz{k}$ are linear in the $q_*$ variables and that we only sum over vectors of positive integers $\ain\in\mathbb{N}^n$. If this assumption is removed, and the appropriate generating series are defined, one encounters situations such as the one below.
\begin{align*}
    & \Biggl(\frac{1}{\h}\sum_{\substack{\ain\in\mathbb{Z}^3\\\lvert\ain\rvert=0}}\frac{p_{\ain}}{3!}\Biggr)\star\Biggl(\frac{1}{\h}\sum_{\substack{\mathbf{b}\in\mathbb{Z}^3\\\lvert\mathbf{b}\rvert=0}}\frac{p_{\ain}}{3!}\Biggr) = \frac{1}{\h^2}\sum_{\substack{\ain\in\mathbb{Z}^3\\\lvert\ain\rvert=0}}\sum_{\substack{\mathbf{b}\in\mathbb{Z}^3\\\lvert\mathbf{b}\rvert=0}}\frac{p_{\ain}}{3!}\frac{p_{\mathbf{b}}}{3!} \\
    & + \frac{1}{\h}\sum_{\substack{\ain\in\mathbb{Z}^4\\\lvert\ain\rvert=0}}\frac{\lvert a_1+a_2\rvert + \lvert a_1+a_3\rvert + \lvert a_1+a_4\rvert + \lvert a_2+a_3\rvert + \lvert a_2+a_4\rvert + \lvert a_3+a_4\rvert}{2}\,\frac{p_\ain}{4!} \\
    & + \sum_{\substack{\ain\in\mathbb{Z}^2\\\lvert\ain\rvert=0}}\frac{\lvert a_1\rvert^3+\lvert a_2\rvert^3 - \lvert a_1\rvert -\lvert a_2\rvert}{24}\,\frac{p_\ain}{2!} =\\ 
    &\frac{1}{\h^2}\sum_{\substack{\ain\in\mathbb{Z}^3\\\lvert\ain\rvert=0}}\sum_{\substack{\mathbf{b}\in\mathbb{Z}^3\\\lvert\mathbf{b}\rvert=0}}\frac{p_{\ain}}{3!}\frac{p_{\mathbf{b}}}{3!}+
\frac{1}{\h}\sum_{\substack{\ain\in\mathbb{Z}^4\\\lvert\ain\rvert=0}}2\max\{|a_i|\}_{i=1}^4\,\frac{p_\ain}{4!}    +
\sum_{a\in \mathbb{N}} \frac{a(a^2-1)}{12}{p_ap_{-a}}
\end{align*}
where we have replaced the variables $q_a$ by $p_{-a}$. 
\begin{proof}[Proof of Theorem \ref{thm:poly}]
    By Proposition \ref{prop:TmuTk} and Theorem \ref{thm:Tmud_in_Tk}, for every partition $\mu=(\mu_1,\dots,\mu_\ell)$ and integer $d\ge0$,
    \begin{equation}\label{eq:Tmu_useful}
        \Tr{\mu,d}^\bullet = \frac{1}{\lvert\operatorname{Aut}(\mu)\rvert}\frac{v^d}{d!}\Trz{1}^{\star d}\star\Trz{\mu_1}\star\cdots\star\Trz{\mu_\ell} + \sum_{\substack{k_1,\dots,k_m\ge1\\k_1+\cdots+k_m\le\lvert\mu\rvert+d\\m<d+\ell}}\tilde{c}_{k_1,\dots,k_m}\Trz{k_1}\star\cdots\star\Trz{k_m},
    \end{equation}
    where $\tilde{c}_{k_1,\dots,k_m}$ are some constants multiplied by appropriate powers of $\h,\e^2$ and $v$. Since the generating series $\Trz{k}$ encode polynomial invariants $x_{k;g}(\ain)$, it is immediate from \eqref{eq:Tmu_useful} and Lemma \ref{lem:starpol} that all the invariants $x_{\mu,d;g}(\ain)$ are polynomial. For a fixed $0\le g\le\lvert\mu\rvert+d$, one extracts $x_{\mu,d;g}(\ain)$ by looking at the coefficient of $\h^{g-1}\e^{2g}\frac{v^{d+g}}{(d+g)!}$ in $\Tr{\mu,d}^\bullet$. Recalling the definition of $\odot$ from Lemma \ref{lem:starpol}, in the leading term of \eqref{eq:Tmu_useful}, one must look at
    \begin{equation*}
        \frac{1}{\lvert\operatorname{Aut}(\mu)\rvert}\frac{v^d}{d!}\Trz{1}^{\odot d}\odot\Trz{\mu_1}\odot\cdots\odot\Trz{\mu_\ell}.
    \end{equation*}
    This corresponds to choosing $g_1,\dots,g_d\in\{0,1\}$ (from the $\Trz{1}$ factors) and $g_{d+i}\le\mu_i$ (from the $\Trz{\mu_i}$ factors) for every $i=1,\dots,\ell$, such that $g_1+\cdots+g_{d+\ell}=g$. Since, by Corollary \ref{cor:polynomiality}, each choice of $g_i$ corresponds to a polynomial of degree $2g_i$, by Lemma \ref{lem:starpol} this will contribute a polynomial of degree $2g+d+\ell-1$ to the expression for $x_{\mu,d;g}(\ain)$. The top degree of this polynomial is nonvanishing by Remark \ref{rem:top_positive}. A similar argument applies to the remaining terms in $\Tr{\mu,d}^\bullet$, but they will contribute a polynomial of strictly lower degree since there are $<d+\ell$ factors in each term. This completes the proof.
\end{proof}
Below we give a more general statement, which includes the intersection numbers arising from disconnected covers too. 
\begin{thm}\label{thm:disc_poly}
    Fix a partition $\mu$, an integer $d\ge0$, a genus $g\ge0$, and the number of connected components $m\ge1$. Fix also $\eta_1,\dots,\eta_m$ (possibly empty) partitions such that $\coprod_{i=1}^m\eta_i=\mu$, genera $g_1,\dots,g_m\ge0$ summing to $g$, and $n_1,\dots,n_m\ge1$. Then, for every $i=1,\dots,m$, the intersection number $x_{\vec{\mu},d;\vec{g}}(\vec{\ain})$ is a polynomial in the $\ain_i\in\mathbb{N}^{n_i}$ of degree
    \begin{equation*}
        3g_i-2+n_i+\ell(\eta_i)-\lvert\eta_i\rvert.
    \end{equation*}
\end{thm}
\begin{proof}
    We may assume that all the conditions of Lemma \ref{lem:conditions} are satisfied, so that $x_{\vec{\mu},d;\vec{g}}(\vec{\ain})$ is not identically zero. As before, equation \eqref{eq:Tmu_useful} and Lemma \ref{lem:starpol} imply that all the invariants $x_{\vec{\mu},d;\vec{g}}(\vec{\ain})$ are polynomial. Also, the highest degree contribution for each polynomial comes from the leading term
    \begin{equation*}
        \Tr{\mu,d}^{\bullet} = \frac{1}{\lvert\operatorname{Aut}(\mu)\rvert}\frac{v^d}{d!}\Trz{1}^{\star d}\star\Trz{\mu_1}\star\cdots\star\Trz{\mu_\ell} + \cdots.
    \end{equation*}
    Let us analyze where the invariant $x_{\vec{\mu},d;\vec{g}}(\vec{\ain})$ appears in this leading term. For every $i$, write $\eta_i=(\mu_{i,1},\dots,\mu_{i,\ell_i})$ where $\ell_i=\ell(\eta_i)$ and the $\mu_{i,j}$ are elements of the partition $\mu$. Let $r_i=2g_i-1+n_i-\lvert\eta_i\rvert$ be the number of additional simple branch points on the $i$-th component. We claim that $x_{\vec{\mu},d;\vec{g}}(\vec{\ain})$ only appears in this piece of the leading term:
    \begin{equation*}
        \prod_{i=1}^m\sum_{h_i=0}^{g_i}[\e^{2g_i-2h_i}]\Trz{1}^{\odot(r_i-g_i)}\odot[\e^{2h_i}]\bigl(\Trz{\mu_{i,1}}\odot\cdots\odot\Trz{\mu_{i,\ell_i}}\bigr),
    \end{equation*}
    where $[\e^{2k}]$ denotes the coefficient of $\e^{2k}$ in the generating series. This is justified by the following argument. Firstly, one takes the standard product over $i=1,\dots,m$ since each $i$ corresponds to a different connected component. Now fix an $i\in\{1,\dots,m\}$. Since the ramification profile over $1$ of the $i$-th connected component is described by the partition $\eta_i=(\mu_{i,1},\dots,\mu_{i,\ell_i})$, one must ``effectively glue'' the generating series $\Trz{\mu_{i,1}},\dots,\Trz{\mu_{i,\ell_i}}$, thus obtaining $\Trz{\mu_{i,1}}\odot\cdots\odot\Trz{\mu_{i,\ell_i}}$. By extracting the coefficient of $\e^{2h_i}$, one obtains a connected curve of genus $h_i$ (formed by gluing $\ell_i$ curves $C_{i,1},\dots,C_{i,\ell_i}$ of genus $h_{i,j}\le\mu_{i,j}$) with ramification profile over $1$ described by $\eta_i$, and with $h_i$ additional simple branch points. In order to obtain a connected curve of genus $g_i$ with $r_i$ additional simple branch points, one must glue to it a curve of genus $g_i-h_i$ with $r_i-h_i$ additional simple branch points, coming from $\Trz{1}^{\star d}$. This latter curve is obtained by gluing $r_i-g_i$ curves together (hence the $\Trz{1}^{\odot(r_i-g_i)})$, of which $g_i-h_i$ of them are genus $1$ curves (which carry two additional simple branch points) and the remaining $r_i-2g_i+h_i$ are genus zero (which carry one additional simple branch point). By Lemma \ref{lem:starpol}, these choices give rise to a polynomial of degree
    \begin{align*}
        \underbrace{2(g_i-h_i) + r_i-g_i-1}_{\text{from }[\e^{2g_i-2h_i}]\Trz{1}^{\odot(r_i-g_i)}} + \underbrace{\sum_{j=1}^{\ell_i}2h_{i,j}+\ell_i-1}_{\text{from }[\e^{2h_i}](\Trz{\mu_{i,1}}\odot\cdots\odot\Trz{\mu_{i,\ell_i}})} + 1 & = r_i+g_i+\ell_i-1 \\
        & = 3g_i-2+n_i+\ell_i-\lvert\eta_i\rvert.
    \end{align*}
    The top degree of this polynomial is nonvanishing by Remark \ref{rem:top_positive}.
\end{proof}

\section{Some fun formulas}
We conclude this work with some formulas arising from the structure of invariants described in this paper that we find  pleasingly simple.

 This first Lemma is an immediate corollary of Proposition \ref{prop:Tds}.

\begin{lem}\label{cor:eps_zero}
    For every partition $\mu$ and integer $d\ge0$,
    \begin{equation}\label{eq:eps_zero}
        \left.\Tr{\mu,d}^\bullet\right\vert_{\e=0} = \left.\frac{v^d}{d!}\Trz{1}^{\star d}\star\Tr{\mu}^\bullet\right\vert_{\e=0}.
    \end{equation}
\end{lem}

Setting $\e=0$ means that we only consider invariants on spaces of admissible covers of genus zero,
\begin{equation*}
    x_{\mu,d;0}(\ain) = \int_{\adm_0(\delta,(\mu+1),\ain)}\psi_2^d.
\end{equation*}
Since the top power of a $\psi$-class is the class of a point in $\overline{\mathcal{M}}_{0,\{0,1,p_1,\dots,p_d,\infty\}}$, this intersection number is a genus zero triple Hurwitz number, with full ramification over one point. So \eqref{eq:eps_zero} is saying that all such Hurwitz numbers can be reconstructed from those Hurwitz numbers that have no additional simple branch points ($d=0$). Since $\mathbb{T}^\bullet_0=\sum_\mu\Tr{\mu}^\bullet$ and $\mathbb{T}^\circ_0=\sum_\mu\Tr{\mu}^\circ$ satisfy the exponential relation $\mathbb{T}^\bullet_0=\exp(\mathbb{T}^\circ_0)$, the above corollary can be written as
\begin{equation*}
    \left.\mathbb{T}^\bullet\right\vert_{\e=0} = \left.\exp_\star\left(v\Trz{1}\right)\star\exp\left(\mathbb{T}_0^\circ\right)\right\vert_{\e=0},
\end{equation*}
where $\mathbb{T}^\bullet=\sum_\mu\sum_d\Tr{\mu,d}^\bullet$ and the ``$\star$-exponential'' is
\begin{equation*}
    \exp_\star(v\Trz{1}) := \sum_{d\ge0}\frac{v^d}{d!}\Trz{1}^{\star d}.
\end{equation*}

In genus zero, it is straightforward to give closed formulas for the invariants  $x_{(k),d;0}(\ain)$. The following result is established using the formula for $H$ from Lemma \ref{lem:starpol}, and Lemma \ref{cor:eps_zero}.
\begin{lem}\label{lem:g0}
 For every $k\ge1$, $d\ge0$,
\begin{equation*}
    x_{(k),d;0}(a_1,\dots,a_{k+d+1}) \; = 
    \; (k+d)!(a_1+\cdots+a_{k+d+1})^d.
\end{equation*}   
\end{lem}

\begin{lem}
For any $k\geq 1$ and $\ain = (a_1, \ldots, a_k)$, we have:
    \begin{equation}
        x_{k;1}(\ain) = \frac{k!}{12}\left(h_2(\ain)- \binom{k+1}{2}\right) = \frac{k!}{24}\left[(\delta^2-k)+\sum_{i=1}^k(a_i^2-k)\right],
    \end{equation}
    with $h_2$ denoting the degree $2$ complete homogeneous symmetric polynomial in the $a_i$.
\end{lem}
\begin{proof}
    The invariants $x_{k;1}(\ain)$ is a quadratic symmetric polynomial, we can therefore express it in the monomial basis as:
    \begin{equation}
     x_{k;1}(\ain) = \alpha_{m_\emptyset}+ \alpha_{m_{(1)}}m_{(1)}(\ain) + \alpha_{m_{(2)}}m_{(2)}(\ain)+\alpha_{m_{(1,1)}}m_{(1,1)}(\ain),
    \end{equation}
 where the $\alpha$'s are coefficients that we wish to determine. We first do so by imposing that the invariants satisfy all WDVV relations of type $\rho_{k;1}(\ain)$, where we use the notation introduced in the proof of Proposition \ref{prop:Tk}, and $\ain$ has now length $k+1$.
 Using the genus zero formula $x_{k;0}(\ain)=k!$ from Lemma \ref{lem:g0}, we evaluate $\rho_{k;1}(\ain)$ to be
 \begin{equation} \label{eq:relfin}
   \left(\delta^3-\sum_{i=1}^{k+1} a_i^3\right)\frac{k!}{12} +\sum_{i=1}^{k+1}(\delta-a_i)x_{k;1}(\ain\smallsetminus a_i)-\sum_{i<j=1}^{k+1}(a_i+a_j)x_{k;1}(\ain\smallsetminus \{a_i,a_j\}\cup a_i+a_j)=0.
 \end{equation}
The left hand side of \eqref{eq:relfin} is a degree three symmetric polynomial in $a_1, \ldots, a_{k+1}$, that can be written in the monomial basis. Setting the coefficient of $m_{(1,1)}$ equal to zero determines $\alpha_{(1)}=0$. The coefficient of $m_{(2,1)}$ gives the equation:
\begin{equation}\label{eq:prima}
    \frac{k!}{4}+(k-1)(\alpha_{(2)}+\alpha_{(1,1)})-((k+2)\alpha_{(2)}+(k-1)\alpha_{(1,1)}) = 0
\end{equation}
 The coefficient of $m_{(1,1,1)}$
 \begin{equation}\label{eq:seconda}
    \frac{k!}{2}+3(k-2)(\alpha_{(1,1)})-3k\alpha_{(1,1)} = 0.
\end{equation}
The two equations \eqref{eq:prima}, \eqref{eq:seconda} determine $\alpha_{(2)} = \alpha_{(1,1)} = k!/12 $, thus showing that 
\begin{equation}\label{eq:daje}
x_{k;1}(\ain) = \alpha_{m_\emptyset}+(k!/12 )h_2(\ain) .
\end{equation}
We compute the constant term by imposing the geometrically trivial information $x_{k;1}(1^k)=0$: $h_2(1^k) = \binom{k+1}{2}$ must be a root of \eqref{eq:daje}, thus concluding the proof.
\end{proof}

We conclude with a closed formula for the invariants $x_{g;g}(a)$.
\begin{lem}\label{cor:xgg}
    For every $g\ge0$ and $a\in\mathbb{N}$:
    \begin{equation*}
        x_{g;g}(a) = \frac{(g!)^2}{2^{2g}(2g+1)!!}\prod_{n=1}^g(a^2-n^2).
    \end{equation*}
\end{lem}
\begin{proof}
    We have
    \begin{equation*}
        x_{g;g}(a) = (g!)^2\sum_{h=0}^g\frac{(g+1)!}{2^h(g+1-h)!}c_{g,h} \int_{\operatorname{DR}_{g-h}(-a,0,a)}\lambda_{g-h}\psi_2^{g-h},
    \end{equation*}
    where $c_{g,h}=[t^{2h}]\left(\frac{\operatorname{arcsinh}(t/2)}{t/2}\right)^{2(g+1-h)}$. By Lemma \ref{lem:DRsinh} $\int_{\operatorname{DR}_g(-1,0,1)}\lambda_g\psi_2^g = \frac{1}{2^{2g}(2g+1)!!}$. Since $\lambda_g\operatorname{DR}_g(-a,0,a)$ is a homogeneous polynomial of degree $2g$ in the variable $a$ \cite{pix23, spe24}, it is completely determined by its value at $a=1$. So
    \begin{equation}\label{eq:DR_!!}
        \int_{\operatorname{DR}_g(-a,0,a)}\lambda_g\psi_2^g = \frac{a^{2g}}{2^{2g}(2g+1)!!}.
    \end{equation}
    The central factorial numbers of the first kind $t(n,k)$ appear in the expansion of the inverse hyperbolic sine \cite{riordan} as
    \begin{equation*}
        \frac{1}{k!}\left(2\operatorname{arcsinh}(t/2)\right)^k = \sum_{n\ge k}t(n,k)\frac{t^n}{n!}.
    \end{equation*}
    Therefore
    \begin{equation*}
        c_{g,h} = \frac{(2g-2h+2)!}{(2g+2)!}t(2g+2,2g-2h+2).
    \end{equation*}
    Combining this with \eqref{eq:DR_!!}, we obtain
    \begin{align*}
        x_{g;g}(a) & = \frac{(g!)^2}{2^{2g}(2g+1)!!}\sum_{h=0}^gt(2g+2,2g-2h+2)a^{2g-2h} \\
        & = \frac{(g!)^2}{2^{2g}(2g+1)!!}\frac{1}{a^2}\sum_{h=0}^gt(2g+2,2h+2)a^{2h+2} = \frac{(g!)^2}{2^{2g}(2g+1)!!}\frac{1}{a^2}a^{[2g+2]} \\
        & = \frac{(g!)^2}{2^{2g}(2g+1)!!}\prod_{n=1}^g(a^2-n^2),
    \end{align*}
    where $a^{[n]}=a\left(a+\frac{n}{2}-1\right)\left(a+\frac{n}{2}-2\right)\cdots\left(a-\frac{n}{2}+1\right)$ is the central factorial.
\end{proof}

\appendix

\section{Programming the recursion}\label{app:python}
Below we present a \texttt{Python} program that can be used to calculate the $x_{k;g}(\ain)$ from the commutation relation $[\Trz{k},\Trz{1}]=0$. We do this step by step for $k=2$, then we show how to do the same for general $k$. The program relies on \texttt{pyBoLaNO}, a \texttt{Python} symbolic package based on \texttt{SymPy} that quickly normal-orders polynomials in annihilation and creation operators \cite{pybolano}. The annihilation operators $\hat{b}_i$ and the creation operators $\hat{b}_i^\dagger$ from quantum mechanics satisfy the canonical commutation relations
\begin{equation*}
    \bigl[\hat{b}_i, \hat{b}_j\bigr] = \bigl[\hat{b}^\dagger_i, \hat{b}^\dagger_j\bigr] = 0, \quad \bigl[\hat{b}_i,\hat{b}^\dagger_j\bigr] = \delta_{ij}.
\end{equation*}
Using these, any polynomial in the $\hat{b}_i,\hat{b}^\dagger_i$ can be normal-ordered so that all the $\hat{b}^\dagger_i$ are to the left and all the $\hat{b}_i$ are to the right, for example
\begin{equation*}
    \hat{b}_2\hat{b}^\dagger_2 = 1 + \hat{b}^\dagger_2\hat{b}_2.
\end{equation*}
The \texttt{pyBoLaNO} package does this via its \texttt{NO} function.
\begin{inputbox}{1}
import pybolano as bl

b_2, bd_2 = bl.ops(2) 
b_2*bd_2, bl.NO(b_2*bd_2)
\end{inputbox}
\begin{outputbox}{1}
    $b_2b_2^\dagger$, $1+b^\dagger_2b_2$
\end{outputbox}
The $p_i,q_i$ variables that we use in our encoding satisfy the commutation relations
    $\left[q_i,q_j\right] = \left[p_i, p_j\right] = 0$ and $\left[q_i,p_j\right] = i\delta_{ij}.$
Therefore we set $q_i=\sqrt{i}\,\hat{b}_i$ and $p_i=\sqrt{i}\,\hat{b}^\dagger_i$.
\begin{inputbox}{2}
import sympy as sm
    
q_2, p_2 = sm.sqrt(2)*bl.ops(2)[0], sm.sqrt(2)*bl.ops(2)[1]
bl.NO(q_2*p_2)
\end{inputbox}
\begin{outputbox}{2}
    $2 + 2b_2^\dagger b_2$
\end{outputbox}
Let us encode
\begin{equation*}
    \Trz{1} = \frac{1}{\h}\sum_{\ain\in\mathbb{N}^2}\frac{p_\ain}{2!}q_{\rvert\ain\lvert} + \e^2v\sum_{a\in\mathbb{N}}\frac{a^2-1}{12}p_aq_a
\end{equation*}
for covers of degree up to $5$. We set $\h=v=1$ since the $\e^2$ parameter alone suffices.
\begin{inputbox}{3}
D = 5 #consider covers of degree <= 5
b = [0] + [bl.ops(sm.Symbol(f"{k}"))[0] for k in range(1, D + 1)]
bd = [0] + [bl.ops(sm.Symbol(f"{k}"))[1] for k in range(1 , D + 1)]
q = [0] + [sm.sqrt(k) * b[k] for k in range(1, D + 1)]
p = [0] + [sm.sqrt(k) * bd[k] for k in range(1, D + 1)]
e = sm.symbols('e') #keeps track of genus
    
T_1_0 = 0 #genus 0 part of T_1
for a1 in range(1, D): 
    for a2 in range(1, D - a1 + 1): #sum over positive integers a1, a2 summing to at most D
        x_1_0 = 1
        T_1_0 += x_1_0*p[a1]*p[a2]*q[a1+a2]/2

T_1_1 = 0 #genus 1 part of T_1
for a in range(1, D + 1):
    x_1_1 = (a**2 - 1)/12
    T_1_1 += e**2*x_1_1*p[a]*q[a]

T_1 = T_1_0 + T_1_1
\end{inputbox}
\vspace{0.5em}
We now want to write
\begin{equation*}
    \Trz{2} = \frac{1}{\h}\sum_{\ain\in\mathbb{N}^3}x_{2;0}(\ain)\frac{p_{\ain}}{3!}q_{\lvert\ain\rvert} + \e^2v\sum_{\ain\in\mathbb{N}^2}x_{2;1}(\ain)\frac{p_\ain}{2!}q_{\lvert\ain\rvert} + \h\e^4\frac{v^2}{2!}\sum_{a\in\mathbb{N}}x_{2;2}(a)p_aq_a
\end{equation*}
in order to compute the three invariants $x_{2;0}(a_1,a_2,a_{3})$, $x_{2;1}(a_1,a_2)$, and $x_{2;2}(a)$. Again, we set $\h=v=1$. In Corollary \ref{cor:polynomiality}, we prove that $x_{k;g}(a_1,\dots,a_{k+1-g})$ is a symmetric polynomial in the $a_i$'s of degree $2g$, whose odd degree parts are zero. We know from Lemma \ref{lem:usingstring} that $x_{2;0}(1,1,1)=2$, so actually $x_{2;0}(a_1,a_2,a_3)=2$ for all $a_1,a_2,a_3$. Let us check that this agrees with $[\Trz{2},\Trz{1}]=0$ by using \texttt{pyBoLaNO}'s commutator function.
\begin{inputbox}{4}
T_2_0 = 0 #genus 0 part of T_2

for a1 in range(1, D - 1):
    for a2 in range(1, D - a1):
        for a3 in range(1, D - a1 - a2 + 1):
            x_2_0 = 2
            T_2_0 += x_2_0*p[a1]*p[a2]*p[a3]*q[a1+a2+a3]/(1*6)

bl.NO_commutator(T_2_0, T_1_0) #genus 0 part of [T_2, T_1]
\end{inputbox}
\begin{outputbox}{4}
    $0$
\end{outputbox}
By Corollary \ref{cor:polynomiality}, we can write the genus $1$ invariant as
\begin{equation*}
    x_{2;1}(a_1,a_2) = c_1(a_1^2 + a_2^2) + c_2a_1a_2 + c_3,
\end{equation*}
where $c_1,c_2,c_3$ are constants to be determined. We know that $x_{2;1}(1,1)=0$ because the ramification index over the internal point is $k+1=3$, so we use this as initial condition. In general, the polynomiality result of Corollary \ref{cor:polynomiality} ensures that if we set the maximum degree $D$ of the covers to be high enough ($D=5$ suffices in this $k=2$ example), the finitely many relations obtained from $[\Trz{k},\Trz{1}]$ are enough to determine all the coefficients of the $x_{k;g}(\ain)$.
\begin{inputbox}{5}
c1, c2, c3 = sm.symbols('c1 c2 c3')

T_2_1 = 0 #genus 1 part of T_2
for a1 in range(1, D):
    for a2 in range(1, D - a1 + 1):
        if (a1, a2) == (1, 1):
            x_2_1 = 0 #initial condition
        else:
            x_2_1 = c1*(a1**2 + a2**2) + c2*a1*a2 + c3
        T_2_1 += e**2*x_2_1*p[a1]*p[a2]*q[a1+a2]/(1*2)

X_1 = bl.NO_commutator(T_2_1, T_1_0) + bl.NO_commutator(T_2_0, T_1_1) #genus 1 part of [T_2, T_1]

X_1
\end{inputbox}
\begin{outputbox}{5}
    $-5\sqrt{3}c_1{b_1^\dagger}^3b_3 - 46\sqrt{2}c_1{b_1^\dagger}^2b_2^\dagger b_4 - 84\sqrt{5}c_1b_1^\dagger {b_2^\dagger}^2b_5 - 41\sqrt{15}c_1{b_1^\dagger}^2b_3^\dagger b_5 -2\sqrt{3}c_2{b_1^\dagger}^3b_3$ \newline
    $- 14\sqrt{2}c_2{b_1^\dagger}^2b_2^\dagger b_4 - 24\sqrt{5}c_2b_1^\dagger{b_2^\dagger}^2b_5 - 10\sqrt{15}c_2{b_1^\dagger}^2b_3^\dagger b_5 - \sqrt{3}c_3{b_1^\dagger}^3b_3 - 2\sqrt{2}{b_1^\dagger}^2b_2^\dagger b_4$ \newline
    $+18\sqrt{5}c_3b_1^\dagger{b_2^\dagger}^2b_5 - \sqrt{15}c_3{b_1^\dagger}^2b_3^\dagger b_5$
\end{outputbox}
The coefficient of each monomial in the $b^\dagger$ and $b$ variables gives a linear equation satisfied by $c_1,c_2,c_3$, which we can solve.
\begin{inputbox}{6}
eq1 = X_1.coeff(e**2 * bd[1]*bd[1]*bd[1] * b[3] * sm.sqrt(3))
eq2 = X_1.coeff(e**2 * bd[1]*bd[1]*bd[2] * b[4] * sm.sqrt(2))
eq3 = X_1.coeff(e**2 * bd[1]*bd[2]*bd[2] * b[5] * sm.sqrt(5))

sm.solve([eq1, eq2, eq3], [c1, c2, c3])
\end{inputbox}
\begin{outputbox}{6}
    \texttt{\{c1:1/6, c2:1/6, c3:-1/2\}}
\end{outputbox}
Therefore $x_{2;1}(a_1,a_2) = \frac{a_1^2+a_2^2}{6} + \frac{a_1a_2}{6} - \frac12$. We continue on to the genus $2$ invariant, which we write as
\begin{equation*}
    x_{2;2}(a) = (a-1)(a-2)(c_4a^2+c_5a+c_6).
\end{equation*}
We already know that $(a-1)$ and $(a-2)$ are factors of $x_{2;2}(a)$ because the invariant vanishes for $a=1$ and $a=2$ due to the same reason as before.
\begin{inputbox}{7}
T_2_1 = T_2_1.subs({c1:1/6, c2:1/6, c3:-1/2}) #substitute genus 1 invariant data
c4, c5, c6 = sm.symbols('c4 c5 c6')

T_2_2 = 0 #genus 2 part of T_2
for a in range(1, D+1):
    x_2_2 = (a - 1)*(a - 2)*(c4*a**2 + c5*a + c6)
    T_2_2 += e**4*x_2_2*p[a]*q[a]/(2*1)

X_2 = bl.NO_commutator(T_2_2, T_1_0) + bl.NO_commutator(T_2_1, T_1_1) #genus 2 part of [T_2, T_1]

eq4 = X_2.coeff(e**4 * bd[1]*bd[2] * b[3] * sm.sqrt(6))
eq5 = X_2.coeff(e**4 * bd[1]*bd[3] * b[4] * sm.sqrt(3))
eq6 = X_2.coeff(e**4 * bd[2]*bd[2] * b[4])

sm.solve([eq4, eq5, eq6], [c4, c5, c6])
\end{inputbox}
\begin{outputbox}{7}
    \texttt{\{c4: 1/12*c6 + 1/72, c5: -7/12*c6 + 5/72\}}
\end{outputbox}
The individual equations extracted from each coefficient of $[\Trz{2},\Trz{1}]$ are not guaranteed to be linearly independent, so we have to look at one more equation to find $c_6$.
\begin{inputbox}{8}
eq7 = X_2.coeff(e**4 * bd[1]*bd[4] * b[5] * sm.sqrt(5))

sm.solve([eq4, eq5, eq6, eq7], [c4, c5, c6])
\end{inputbox}
\begin{outputbox}{8}
    \texttt{\{c4:1/60, c5:1/20, c6:1/30\}}
\end{outputbox}
Therefore $x_{2;2}(a)=(a-1)(a-2)\left(\frac{a^2}{60}+\frac{a}{20}+\frac{1}{30}\right)=\frac{(a^2-1)(a^2-4)}{60}$. The general building block invariants are 
\begin{equation*}
    x_{k;g}(a_1,\dots,a_{k+1-g}) = \sum_{h=0}^g\sum_{\substack{\nu\vdash 2h\\\ell(\nu)\le k+1-g}}c_{g,\nu}\,m_\nu(a_1,\dots,a_{k+1-g}), \quad c_{g,\nu}\in\mathbb{Q},
\end{equation*}
where $m_\nu$ are the monomial symmetric polynomials. The following function writes $\Trz{k}$ up to some degree $D$ with unknown coefficients $c_{g,\nu}$.
\begin{inputbox}{9}
from sympy.utilities.iterables import partitions, multiset_permutations

eps = sm.Symbol(r'\varepsilon')
v = sm.Symbol('v')

def T_series(k, D):
    bd, b, p, q = [0], [0], [0], [0]
    for i in range(1, D+1):
        bd.append(bl.ops(sm.Symbol(f"{i}"))[1])
        p.append(sm.sqrt(i)*bl.ops(sm.Symbol(f"{i}"))[1])
        b.append(bl.ops(sm.Symbol(f"{i}"))[0])
        q.append(sm.sqrt(i)*bl.ops(sm.Symbol(f"{i}"))[0])
    Tk = 0
    for g in range(k+1): #sum over genus g from 0 to k
        ell = k + 1 - g
        unknown_coeffs = [] #build the basis of unknown symmetric polynomials for x_{k;g}
        if g > 0:
            for deg in range(0, 2*g+1, 2): #iterate over even degrees
                if deg == 0:
                    sym = sm.Symbol(f'c_{{{g},()}}')
                    unknown_coeffs.append(((), sym))
                else:
                    for mu_dict in partitions(deg):
                        if sum(mu_dict.values()) > ell: #disregard partitions of length >ell
                            continue
                        mu_tuple = tuple(part for part in sorted(mu_dict.keys(), reverse=True) for _ in range(mu_dict[part]))
                        mu_str = ",".join(map(str, mu_tuple))
                        sym = sm.Symbol(f'c_{{{g},({mu_str})}}')
                        unknown_coeffs.append((mu_tuple, sym))
#continues on next page...
\end{inputbox}
\begin{inputbox}{9}
        def evaluate_xkg_monomial(lam_tuple): #define evaluation of x_{k;g} on a specific partition lambda
            if g == 0:
                return sm.factorial(k) #x_{k;0} = k!
            val = 0
            for mu_tuple, sym in unknown_coeffs:
                if not mu_tuple: #constant term (degree 0)
                    val += sym
                    continue
                mu_padded = list(mu_tuple) + [0] * (ell - len(mu_tuple)) #pad the partition mu with zeros to match the length ell
                m_val = 0 
                for perm in multiset_permutations(mu_padded): #evaluate monomial symmetric polynomial m_mu on lam_tuple
                    term_val = 1
                    for a, power in zip(lam_tuple, perm):
                        term_val *= (a ** power)
                    m_val += term_val
                val += sym * m_val
            return val

        for d in range(1, D+1): #sum over all degrees d up to D
            for lam_dict in partitions(d, m=ell):
                if sum(lam_dict.values()) != ell: #sum over partitions of d of length exactly ell
                    continue
                lam_tuple = tuple(part for part in sorted(lam_dict.keys(), reverse=True) for _ in range(lam_dict[part])) #convert lambda dictionary to tuple
                xkg_val = evaluate_xkg_monomial(lam_tuple) #evaluate the polynomial for this specific combination of a_i
                p_prod = 1 #calculate the p terms and the symmetry factor (1 / prod(m_i!))
                m_fact_prod = 1
                for part, mult in lam_dict.items():
                    p_prod *= (p[part])**mult
                    m_fact_prod *= sm.factorial(mult)
                term = (eps**(2*g) * v**g / sm.factorial(g)) * xkg_val * (p_prod / m_fact_prod) * q[d]
                Tk += term
    return Tk

T_series(2,3)
\end{inputbox}
\begin{outputbox}{9}
    $\e^4v^2(c_{2,()} + c_{2,(2)} + c_{2,(4)})b_1^\dagger b_1 + \e^4v^2(c_{2,()} + 4c_{2,(2)} + 16c_{2,(4)})b_2^\dagger b_2$ \newline
    $+ \e^4v^2(c_{2,()} + 9c_{2,(2)} + 81c_{2,(4)})b_3^\dagger b_3 + \sqrt{2}\e^2v(c_{1,()} + c_{1,(1,1)} + 2c_{1,(2)}){b_1^\dagger}^2b_2$ \newline
    $2\sqrt{2}\e^2v(c_{1,()} + 2c_{1,(1,1)} + 5c_{1,(2)})b_2^\dagger b_1^\dagger b_3 + \frac{4{b_1^\dagger}^3b_3}{3}$
\end{outputbox}

\section{Explicit examples of generating functions}\label{app:ex}

We collect here some explicit examples of the generating functions we reconstruct.

\subsection{The series $\Trz{k}$}\label{appsec:Tk}
We collect here the series  $\Trz{k}$ for $k\leq 4$.

\begin{align*}
\Trz{1}
&= \frac{1}{2\hbar}
   \sum_{a_1,a_2\ge1}
   p_{a_1}p_{a_2}q_{a_1+a_2}
 + \frac{\varepsilon^2v}{12}
   \sum_{a\ge1}(a^2-1)\,p_aq_a,
\\[1em]
\Trz{2}
&= \frac{1}{3\hbar}
   \sum_{\mathbf a\in\mathbb N^3}
   p_{\mathbf a}q_{|\mathbf a|}
 + \frac{\varepsilon^2v}{12}
   \sum_{a_1,a_2\ge1}
   (a_1^2+a_2^2+a_1a_2-3)\,
   p_{a_1}p_{a_2}q_{a_1+a_2}
\\
&\quad
 + \frac{\hbar\varepsilon^4v^2}{120}
   \sum_{a\ge1}
   (a^2-1)(a^2-4)\,p_aq_a,
\\[1em]
\Trz{3}
&= \frac{1}{4\hbar}
   \sum_{\mathbf a\in\mathbb N^4}
   p_{\mathbf a}q_{|\mathbf a|}
 + \frac{\varepsilon^2v}{12}
   \sum_{\mathbf a\in\mathbb N^3}
   \left(
      \sum_{i=1}^3 a_i^2
      + \sum_{1\le i<j\le3}a_ia_j - 6
   \right)p_{\mathbf a}q_{|\mathbf a|}
\\
&\quad
 + \frac{\hbar\varepsilon^4v^2}{80}
   \sum_{a_1,a_2\ge1}
   \Bigl(
      a_1^4+a_2^4
      +2a_1a_2(a_1^2+a_2^2)+3a_1^2a_2^2
\\
&\hspace{9em}
      -10(a_1^2+a_2^2+a_1a_2)+21
   \Bigr)p_{a_1}p_{a_2}q_{a_1+a_2}
\\
&\quad
 + \frac{\hbar^2\varepsilon^6v^3}{1120}
   \sum_{a\ge1}
   (a^2-1)(a^2-4)(a^2-9)\,p_aq_a,
\\[1em]
\Trz{4}
&= \frac{1}{5\hbar}
   \sum_{\mathbf a\in\mathbb N^5}
   p_{\mathbf a}q_{|\mathbf a|}
 + \frac{\varepsilon^2v}{12}
   \sum_{\mathbf a\in\mathbb N^4}
   \left(
      \sum_{i=1}^4 a_i^2
      + \sum_{1\le i<j\le4}a_ia_j - 10
   \right)p_{\mathbf a}q_{|\mathbf a|}
\\
&\quad
 + \frac{\hbar\varepsilon^4v^2}{180}
   \sum_{\mathbf a\in\mathbb N^3}
   \Biggl(
      3\sum_{i=1}^3 a_i^4
      +6\sum_{\substack{1\le i,j\le3\\i\ne j}}a_i^3a_j
      +9\sum_{1\le i<j\le3}a_i^2a_j^2
\\
\end{align*}
\begin{align*}
&\hspace{9em}
      +11a_1a_2a_3(a_1+a_2+a_3)
      -50\left(
         \sum_{i=1}^3 a_i^2
         +\sum_{1\le i<j\le3}a_ia_j
      \right)
      +195
   \Biggr)p_{\mathbf a}q_{|\mathbf a|}
\\
&\quad
 + \frac{\hbar^2\varepsilon^6v^3}{1680}
   \sum_{a_1,a_2\ge1}
   \Bigl(
      3(a_1^6+a_2^6)+9a_1a_2(a_1^4+a_2^4)
\\
&\hspace{9em}
      +19a_1^2a_2^2(a_1^2+a_2^2)+23a_1^3a_2^3
      -70(a_1^4+a_2^4)
\\
&\hspace{9em}
      -140a_1a_2(a_1^2+a_2^2)-210a_1^2a_2^2
      +455(a_1^2+a_2^2+a_1a_2)-820
   \Bigr)p_{a_1}p_{a_2}q_{a_1+a_2}
\\
&\quad
 + \frac{\hbar^3\varepsilon^8v^4}{10080}
   \sum_{a\ge1}
   (a^2-1)(a^2-4)(a^2-9)(a^2-16)\,p_aq_a.
\end{align*}

\subsection{Examples from Proposition \ref{prop:TmuTk}}\label{appsec:TmuTk}

    For $\lvert\mu\rvert\le4$:
    \begin{align*}
        \Tr{(1,1)}^\bullet & = \frac{1}{2}\Trz{1}^{\star 2} -\frac{3}{2}\Trz{2} - \h\e^2v\frac14\Trz{1}, \\
        \Tr{(2,1)}^\bullet & = \Trz{2}\star\Trz{1} - 4\Trz{3} - \hbar\e^2v\,2 \Trz{2}, \\
        \Tr{(1,1,1)}^\bullet & = \frac{1}{6}\Trz{1}^{\star 3} -\frac{3}{2}\Trz{2} \star \Trz{1}+ \frac{10}{3}\Trz{3} - \h\e^2v\left(\frac14\Trz{1}^{\star2} - \frac52\Trz{2} \right) + \h^2\e^4\frac{v^2}{2!}\frac16\Trz{1}, \\
        \Tr{(3,1)}^\bullet = & \; \Trz{3}\star\Trz{1} - 5\Trz{4} - \h\e^2v\,5\Trz{3}, \\
        \Tr{(2,2)}^\bullet = & \; \frac12\Trz{2}^{\star2} - \frac52\Trz{4} - \h\e^2v\frac{10}{3}\Trz{3} - \h^2\e^4\frac{v^2}{2!}\Trz{2}, \\
        \Tr{(2,1,1)}^\bullet = & \; \frac12\Trz{2}\star\Trz{1}^{\star2} - 4\Trz{3}\star\Trz{1} - \frac32\Trz{2}^{\star2} + 15\Trz{4} - \h\e^2v\left(\frac{9}{4}\Trz{2}\star\Trz{1} - 25\Trz{3}\right) + \h^2\e^4\frac{v^2}{2!}8\Trz{2} \\
        \Tr{(1,1,1,1)}^\bullet = & \; \frac{1}{24}\Trz{1}^{\star4} - \frac34\Trz{2}\star\Trz{1}^{\star2} + \frac{10}{3}\Trz{3}\star\Trz{1} + \frac98\Trz{2}^{\star2} - \frac{35}{4}\Trz{4} \\
        & - \h\e^2v\left(\frac18\Trz{1}^{\star3} - \frac{23}{8}\Trz{2}\star\Trz{1} + \frac{35}{2}\Trz{3}\right) + \h^2\e^4\frac{v^2}{2!}\left(\frac{11}{48}\Trz{1}^{\star2} - \frac{63}{8}\Trz{2}\right) - \h^3\e^6\frac{v^3}{3!}\frac{3}{16}\Trz{1}.
    \end{align*}

\subsection{Examples from Proposition \ref{prop:Tds}}\label{appsec:TmudTmuprop}
\begingroup

\setlength{\jot}{2pt}
    For $d=1$ and $\lvert\mu\rvert\le4$:
{\small
\setlength{\jot}{3pt}
\begin{alignat*}{2}
\Tr{(1),1}^\bullet
&= v\Trz{1}\star\Tr{(1)}^\bullet
 - \h\e^2\frac{v^2}{2!}\frac12\Tr{(1)}^\bullet,
&\qquad
\Tr{(2),1}^\bullet
&= v\Trz{1}\star\Tr{(2)}^\bullet
 - \h\e^2\frac{v^2}{2!}2\Tr{(2)}^\bullet,
\\
\Tr{(1,1),1}^\bullet
&= v\Trz{1}\star\Tr{(1,1)}^\bullet
 - \h\e^2\frac{v^2}{2!}\Tr{(1,1)}^\bullet,
&
\Tr{(3),1}^\bullet
&= v\Trz{1}\star\Tr{(3)}^\bullet
 - \h\e^2\frac{v^2}{2!}5\Tr{(3)}^\bullet,
\\
\Tr{(2,1),1}^\bullet
&= v\Trz{1}\star\Tr{(2,1)}^\bullet
 - \h\e^2\frac{v^2}{2!}\frac52\Tr{(2,1)}^\bullet,
&
\Tr{(1,1,1),1}^\bullet
&= v\Trz{1}\star\Tr{(1,1,1)}^\bullet
 - \h\e^2\frac{v^2}{2!}\frac32\Tr{(1,1,1)}^\bullet,
\\
\Tr{(4),1}^\bullet
&= v\Trz{1}\star\Tr{(4)}^\bullet
 - \h\e^2\frac{v^2}{2!}10\Tr{(4)}^\bullet,
&
\Tr{(3,1),1}^\bullet
&= v\Trz{1}\star\Tr{(3,1)}^\bullet
 - \h\e^2\frac{v^2}{2!}\frac{11}{2}\Tr{(3,1)}^\bullet,
\\
\Tr{(2,2),1}^\bullet
&= v\Trz{1}\star\Tr{(2,2)}^\bullet
 - \h\e^2\frac{v^2}{2!}4\Tr{(2,2)}^\bullet,
&
\Tr{(2,1,1),1}^\bullet
&= v\Trz{1}\star\Tr{(2,1,1)}^\bullet
 - \h\e^2\frac{v^2}{2!}3\Tr{(2,1,1)}^\bullet,
\\
\Tr{(1,1,1,1),1}^\bullet
&= v\Trz{1}\star\Tr{(1,1,1,1)}^\bullet
 - \h\e^2\frac{v^2}{2!}2\Tr{(1,1,1,1)}^\bullet.
\end{alignat*}
}

    
    For $d = 2$ and $\lvert\mu\rvert\le 3$:
    \begin{align*}
        \Tr{(1),2}^\bullet = & \; \frac{v}{2}\Trz{1}\star\Tr{(1),1}^\bullet - \h\e^2\frac{v^3}{3!}\left(\frac52\Tr{(1,1)}^\bullet + \frac{27}{4}\Tr{(2)}^\bullet\right) - \h^2\e^4\frac{v^4}{4!}\frac34\Tr{(1)}^\bullet, \\
        \Tr{(2),2}^\bullet = & \; \frac{v}{2}\Trz{1}\star\Tr{(2),1}^\bullet - \h\e^2\frac{v^3}{3!}\left(\frac72\Tr{(2,1)}^\bullet + 24\Tr{(3)}^\bullet\right) - \h^2\e^4\frac{v^4}{4!}12\Tr{(2)}^\bullet, \\
        \Tr{(1,1),2}^\bullet = & \; \frac{v}{2}\Trz{1}\star\Tr{(1,1),1}^\bullet - \h\e^2\frac{v^3}{3!}\left(6\Tr{(1,1,1)}^\bullet + 9\Tr{(2,1)}^\bullet + 11\Tr{(3)}^\bullet\right) - \h^2\e^4\frac{v^4}{4!}3\Tr{(1,1)}^\bullet, \\
        \Tr{(3),2}^\bullet = & \; \frac{v}{2}\Trz{1}\star\Tr{(3),1}^\bullet - \h\e^2\frac{v^3}{3!}\left(8\Tr{(3,1)}^\bullet + \frac{125}{2}\Tr{(4)}^\bullet\right) - \h^2\e^4\frac{v^4}{4!}75\Tr{(3)}^\bullet, \\
        \Tr{(2,1),2}^\bullet = & \; \frac{v}{2}\Trz{1}\star\Tr{(2,1),1}^\bullet - \h\e^2\frac{v^3}{3!}\left(\frac{17}{2}\Tr{(2,1,1)}^\bullet + 27\Tr{(3,1)}^\bullet + \frac{63}{2}\Tr{(2,2)}^\bullet + \frac{225}{4}\Tr{(4)}^\bullet\right) - \h^2\e^4\frac{v^4}{4!}\frac{75}{4}\Tr{(2,1)}^\bullet, \\
        \Tr{(1,1,1),2}^\bullet = & \; \frac{v}{2}\Trz{1}\star\Tr{(1,1,1),1}^\bullet - \h\e^2\frac{v^3}{3!}\left(11\Tr{(1,1,1,1)}^\bullet + \frac{45}{4}\Tr{(2,1,1)}^\bullet + \frac{25}{2}\Tr{(3,1)}^\bullet\right) -\h^2\e^4\frac{v^4}{4!}\frac{27}{4}\Tr{(1,1,1)}^\bullet.
    \end{align*}
    For $d=3$ and $\lvert\mu\rvert\le 2$:
    \begin{align*}
        \Tr{(1),3}^\bullet = & \; \frac{v}{3}\Trz{1}\star\Tr{(1),2}^\bullet - \h\e^2\frac{v^4}{4!}\left(14\Tr{(1,1,1)}^\bullet + 32\Tr{(2,1)}^\bullet + \frac{304}{3}\Tr{(3)}^\bullet\right) - \h^2\e^4\frac{v^5}{5!}\left(\frac{43}{3}\Tr{(1,1)}^\bullet + \frac{147}{2}\Tr{(2)}^\bullet\right) \\
        & \; - \h^3\e^6\frac{v^6}{6!}\frac{15}{8}\Tr{(1)}^\bullet, \\
        \Tr{(2),3}^\bullet = & \; \frac{v}{3}\Trz{1}\star\Tr{(2),2}^\bullet - \h\e^2\frac{v^4}{4!}\left(\frac{32}{3}\Tr{(2,1,1)}^\bullet + \frac{208}{3}\Tr{(3,1)}^\bullet + 38\Tr{(2,2)}^\bullet + \frac{1025}{3}\Tr{(4)}^\bullet\right) \\
        & - \h^2\e^4\frac{v^5}{5!}\left(\frac{164}{3}\Tr{(2,1)}^\bullet + 672\Tr{(3)}^\bullet\right) - \h^3\e^6\frac{v^6}{6!}120\Tr{(2)}^\bullet, \\
        \Tr{(1,1),3}^\bullet = & \; \frac{v}{3}\Trz{1}\star\Tr{(1,1),2}^\bullet - \h\e^2\frac{v^4}{4!}\left(40\Tr{(1,1,1,1)}^\bullet + 69\Tr{(2,1,1)}^\bullet + \frac{448}{3}\Tr{(3,1)}^\bullet + 108\Tr{(2,2)}^\bullet + \frac{950}{3}\Tr{(4)}^\bullet\right) \\
        & - \h^2\e^4\frac{v^5}{5!}\left(54\Tr{(1,1,1)}^\bullet + 126\Tr{(2,1)}^\bullet + \frac{892}{3}\Tr{(3)}^\bullet\right) - \h^3\e^6\frac{v^6}{6!}15\Tr{(1,1)}^\bullet.
    \end{align*}
    For $d=4$ and $\mu=(1)$:
    \begin{align*}
        \Tr{(1),4}^\bullet = & \; \frac{v}{4}\Trz{1}\star\Tr{(1),3}^\bullet - \h\e^2\frac{v^5}{5!}\left(90\Tr{(1,1,1,1)}^\bullet + \frac{375}{2}\Tr{(2,1,1)}^\bullet + 545\Tr{(3,1)}^\bullet + 360\Tr{(2,2)}^\bullet + 1750\Tr{(4)}^\bullet\right) \\
        & - \h^2\e^4\frac{v^6}{6!}\left(\frac{759}{4}\Tr{(1,1,1)}^\bullet + \frac{5499}{8}\Tr{(2,1)}^\bullet + 4042\Tr{(3)}^\bullet\right) - \h^3\e^6\frac{v^7}{7!}\left(\frac{1593}{16}\Tr{(1,1)}^\bullet + \frac{31887}{32}\Tr{(2)}^\bullet\right) \\
        & - \h^4\e^8\frac{v^8}{8!}\frac{105}{16}\Tr{(1)}^\bullet.
    \end{align*}

\subsection{Examples from Theorem \ref{thm:Tmud_in_Tk}}\label{appsec:Tmud_in_Tk}
 {\small
    For $d=1$ and $\lvert\mu\rvert\le4$:
   {\small
    \begin{align*}
        \Tr{(1),1}^\bullet = & \; v\Trz{1}^{\star2} - \h\e^2\frac{v^2}{2!}\frac12\Trz{1}, \qquad \Tr{(2),1}^\bullet = v\Trz{2}\star\Trz{1} - \h\e^2\frac{v^2}{2!}2\Trz{2}, \\
        \Tr{(1,1),1}^\bullet = & \; v\left(\frac12\Trz{1}^{\star3} - \frac32\Trz{2}\star\Trz{1}\right) - \h\e^2\frac{v^2}{2!}\left(\Trz{1}^{\star2} - \frac32\Trz{2}\right) + \h^2\e^4\frac{v^3}{3!}\frac{3}{4}\Trz{1}, \\
        \Tr{(3),1}^\bullet = & \; v\Trz{3}\star\Trz{1} - \h\e^2\frac{v^2}{2!}5\Trz{3}, \\       
        \Tr{(2,1),1}^\bullet = & \; v\left(\Trz{2}\star\Trz{1}^{\star2} - 4\Trz{3}\star\Trz{1}\right) - \h\e^2\frac{v^2}{2!}\left(\frac{13}{2}\Trz{2}\star\Trz{1} - 10\Trz{3}\right) + \h^2\e^4\frac{v^3}{3!}15\Trz{2}, \\
        \Tr{(1,1,1),1}^\bullet = & \; v\left(\frac16\Trz{1}^{\star4} - \frac32\Trz{2}\star\Trz{1}^{\star2} + \frac{10}{3}\Trz{3}\star\Trz{1}\right) - \h\e^2\frac{v^2}{2!}\left(\frac34\Trz{1}^{\star3} - \frac{29}{4}\Trz{2}\star\Trz{1} 
        + 5\Trz{3}\right) \\
        & + \h^2\e^4\frac{v^3}{3!}\left(\frac{13}{8}\Trz{1}^{\star2} - \frac{45}{4}\Trz{2}\right) - \h^3\e^6\frac{v^4}{4!}\frac32\Trz{1}, \\
        \Tr{(4),1}^\bullet = & \; v\Trz{4}\star\Trz{1} - \h\e^2\frac{v^2}{2!}10\Trz{4}, \\
        \Tr{(3,1),1}^\bullet = & \; v\left(\Trz{3}\star\Trz{1}^{\star2} - 5\Trz{4}\star\Trz{1}\right) - \h\e^2\frac{v^2}{2!}\left(\frac{31}{2}\Trz{3}\star\Trz{1} - \frac{55}{2}\Trz{4}\right) + \h^2\e^4\frac{v^3}{3!}\frac{165}{2}\Trz{3}, \\
        \Tr{(2,2),1}^\bullet = & \; v\left(\frac12\Trz{2}^{\star2}\star\Trz{1} - \frac52\Trz{4}\star\Trz{1}\right) - \h\e^2\frac{v^2}{2!}\left(\frac{20}{3}\Trz{3}\star\Trz{1} + 2\Trz{2}^{\star2} - 10\Trz{4}\right) \\
        & - \h^2\e^4\frac{v^3}{3!}\left(3\Trz{2}\star\Trz{1} - 40\Trz{3}\right) + \h^3\e^6\frac{v^4}{4!}24\Trz{2},  \\
        \Tr{(2,1,1),1}^\bullet = & \; v\left(\frac12\Trz{2}\star\Trz{1}^{\star3} - 4\Trz{3}\star\Trz{1}^{\star2} - \frac32\Trz{2}^{\star2}\star\Trz{1} + 15\Trz{4}\star\Trz{1}\right) \\
        & - \h\e^2\frac{v^2}{2!}\left(6\Trz{2}\star\Trz{1}^{\star2} - 62\Trz{3}\star\Trz{1} - \frac92\Trz{2}^{\star2} + 45\Trz{4}\right) \\
        & + \h^2\e^4\frac{v^3}{3!}\left(\frac{177}{4}\Trz{2}\star\Trz{1} - 225\Trz{3}\right) - \h^3\e^6\frac{v^4}{4!}144\Trz{2}, \\
        \Tr{(1,1,1,1),1}^\bullet = & \; v\left(\frac{1}{24}\Trz{1}^{\star5} - \frac34\Trz{2}\star\Trz{1}^{\star3} + \frac{10}{3}\Trz{3}\star\Trz{1}^{\star2} + \frac98\Trz{2}^{\star2}\star\Trz{1} - \frac{35}{4}\Trz{4}\star\Trz{1}\right) \\
        & - \h\e^2\frac{v^2}{2!}\left(\frac13\Trz{1}^{\star4} - \frac{29}{4}\Trz{2}\star\Trz{1}^{\star2} + \frac{125}{3}\Trz{3}\star\Trz{1} + \frac94\Trz{2}^{\star2} - \frac{35}{2}\Trz{4}\right) \\
        & + \h^2\e^4\frac{v^3}{3!}\left(\frac{23}{16}\Trz{1}^{\star3} - \frac{327}{8}\Trz{2}\star\Trz{1} + 105\Trz{3}\right) - \h^3\e^6\frac{v^4}{4!}\left(\frac{7}{2}\Trz{1}^{\star2} - \frac{189}{2}\Trz{2}\right) \\
        & + \h^4\e^8\frac{v^5}{5!}\frac{15}{4}\Trz{1}.
    \end{align*}
    For $d=2$ and $\lvert\mu\rvert\le3$:
    \begin{align*}
        \Tr{(1),2}^\bullet = & \; \frac{v^2}{2!}\Trz{1}^{\star3} - \h\e^2\frac{v^3}{3!}\left(2\Trz{1}^{\star2} + 3\Trz{2}\right) + \h^2\e^4\frac{v^4}{4!}\frac74\Trz{1}, \\
        \end{align*}
\begin{align*}
        \Tr{(2),2}^\bullet = & \; \frac{v^2}{2!}\Trz{2}\star\Trz{1}^{\star2} - \h\e^2\frac{v^3}{3!}\left(\frac{13}{2}\Trz{2}\star\Trz{1} + 10\Trz{3}\right) + \h^2\e^4\frac{v^4}{4!}16\Trz{2}, \\
        \Tr{(1,1),2}^\bullet = & \; \frac{v^2}{2!}\left(\frac12\Trz{1}^{\star4} - \frac32\Trz{2}\star\Trz{1}^{\star2}\right) - \h\e^2\frac{v^3}{3!}\left(\frac52\Trz{1}^{\star3} - \frac94\Trz{2}\star\Trz{1} - 5\Trz{3}\right) \\
        & + \h^2\e^4\frac{v^4}{4!}\left(6\Trz{1}^{\star2} + \frac{33}{2}\Trz{2}\right) - \h^3\e^6\frac{v^5}{5!}\frac{25}{4}\Trz{1}, \\
        \Tr{(3),2}^\bullet = & \; \frac{v^2}{2!}\Trz{3}\star\Trz{1}^{\star2} - \h\e^2\frac{v^3}{3!}\left(\frac{31}{2}\Trz{3}\star\Trz{1} + \frac{45}{2}\Trz{2}\right) + \h^2\e^4\frac{v^4}{4!}85\Trz{3}, \\
        \Tr{(2,1),2}^\bullet = & \; \frac{v^2}{2!}\left(\Trz{2}\star\Trz{1}^{\star3} - \Trz{3}\star\Trz{1}^{\star2}\right) - \h\e^2\frac{v^3}{3!}\left(14\Trz{2}\star\Trz{1}^{\star2} - 22\Trz{3}\star\Trz{1} + 3\Trz{2}^{\star2} - 30\Trz{4}\right) \\
        & + \h^2\e^4\frac{v^4}{4!}\left(\frac{351}{4}\Trz{2}\star\Trz{1} + 185\Trz{2}\right) - \h^3\e^6\frac{v^5}{5!}\frac{355}{2}\Trz{2}, \\
        \Tr{(1,1,1),2}^\bullet = & \; \frac{v^2}{2!}\left(\frac16\Trz{1}^{\star5} - \frac32\Trz{2}\star\Trz{1}^{\star3} + \frac{10}{3}\Trz{3}\star\Trz{1}^{\star2}\right) \\
        & - \h\e^2\frac{v^3}{3!}\left(\frac{19}{12}\Trz{1}^{\star4} - \frac{27}{2}\Trz{2}\star\Trz{1}^{\star2} + \frac{35}{3}\Trz{3}\star\Trz{1} - \frac92\Trz{2}^{\star2} + 10\Trz{4}\right) \\
        & + \h^2\e^4\frac{v^4}{4!}\left(\frac{61}{8}\Trz{1}^{\star3} - \frac{301}{8}\Trz{2}\star\Trz{1} - \frac{255}{2}\Trz{3}\right) - \h^3\e^6\frac{v^5}{5!}\left(\frac{985}{48}\Trz{1}^{\star2} + \frac{945}{8}\Trz{2}\right) \\
        & + \h^4\e^8\frac{v^6}{6!}\frac{195}{8}\Trz{1}.
    \end{align*}
    For $d=3$ and $\lvert\mu\rvert \le 2$:
    \begin{align*}
        \Tr{(1),3}^\bullet = & \; \frac{v^3}{3!}\Trz{1}^{\star4} - \h\e^2\frac{v^4}{4!}\bigl(5\Trz{1}^{\star3} + 15\Trz{2}\star\Trz{1} + 20\Trz{3}\bigr) + \h^2\e^4\frac{v^5}{5!}\left(\frac{53}{4}\Trz{1}^{\star2} + 93\Trz{2}\right) \\
        & - \h^3\e^6\frac{v^6}{6!}\frac{123}{8}\Trz{1}, \\
        \Tr{(2),3}^\bullet = & \; \frac{v^3}{3!}\Trz{2}\star\Trz{1}^{\star3} - \h\e^2\frac{v^4}{4!}\left(14\Trz{2}\star\Trz{1}^{\star2} + 40\Trz{3}\star\Trz{1} + 3\Trz{2}^{\star2} + 60\Trz{4}\right) \\
        & + \h^2\e^4\frac{v^5}{5!}\left(92\Trz{2}\star\Trz{1} + 580\Trz{3}\right) - \h^3\e^6\frac{v^6}{6!}172\Trz{2}, \\
        \Tr{(1,1),3}^\bullet = & \; \frac{v^3}{3!}\left(\frac12\Trz{1}^{\star5} - \frac32\Trz{2}\star\Trz{1}^{\star3}\right) - \h\e^2\frac{v^4}{4!}\left(5\Trz{1}^{\star4} + \frac32\Trz{2}\star\Trz{1}^{\star2} - \frac92\Trz{2}^{\star2} - 15\Trz{4}\right) \\
        & + \h^2\e^4\frac{v^5}{5!}\left(26\Trz{1}^{\star3} + \frac{735}{4}\Trz{2}\star\Trz{1} + 435\Trz{3}\right) - \h^3\e^6\frac{v^6}{6!}\left(\frac{153}{2}\Trz{1}^{\star2} + \frac{2421}{2}\Trz{2}\right) \\
        & + \h^4\e^8\frac{v^7}{7!}\frac{399}{4}\Trz{1}.
    \end{align*}
    For $d=4$ and $\mu=(1)$:
    \begin{align*}
        \Tr{(1),4}^\bullet = & \; \frac{v^4}{4!}\Trz{1}^{\star5} - \h\e^2\frac{v^5}{5!}\left(10\Trz{1}^{\star4} + 45\Trz{2}\star\Trz{1}^{\star2} + 120\Trz{3}\star\Trz{1} + 150\Trz{4}\right) \\
        & + \h^2\e^4\frac{v^6}{6!}\left(\frac{223}{4}\Trz{1}^{\star3} + \frac{1431}{2}\Trz{2}\star\Trz{1} + 2950\Trz{3}\right) - \h^3\e^6\frac{v^7}{7!}\left(\frac{711}{4}\Trz{1}^{\star2} + \frac{14403}{3}\Trz{2}\right) \\
        & + \h^4\e^8\frac{v^8}{8!}\frac{4033}{16}\Trz{1}.
    \end{align*}}}


\bibliographystyle{alpha}
\bibliography{biblio}

@misc{nss2025,
    author={Nesterov, Denis and Schimpf, Maximilian and Schmitt, Johannes},
    title={Admissible covers and stable maps}, 
    note = {arxiv:2505.03487},
    year={2025}}

@article{pybolano,
    title = {py{B}o{L}a{NO}: A Python symbolic package for normal ordering involving bosonic ladder operators},
    author = {Lim, Hendry M. and Dwiputra, Donny and Ukhtary, M. Shoufie and Nugraha, Ahmad R. T.},
    journal = {Computer Physics Communications},
    volume = {313},
    pages = {109622},
    year = {2025}}

@article{Graber-Pandharipande,
    AUTHOR = {Graber, T. and Pandharipande, R.},
     TITLE = {Localization of virtual classes},
   JOURNAL = {Invent. Math.},
  FJOURNAL = {Inventiones Mathematicae},
    VOLUME = {135},
      YEAR = {1999},
    NUMBER = {2},
     PAGES = {487--518}}

@article{jppz17,
    AUTHOR = {Janda, Felix and Pandharipande, Rahul and Pixton, Aaron and Zvonkine, Dimitri},
     TITLE = {Double ramification cycles on the moduli spaces of curves},
   JOURNAL = {Publ. Math. Inst. Hautes \'Etudes Sci.},
  FJOURNAL = {Publications Math\'ematiques. Institut de Hautes \'Etudes
              Scientifiques},
    VOLUME = {125},
      YEAR = {2017},
     PAGES = {221--266}}

@article{bur15,
    AUTHOR = {Buryak, Alexandr},
     TITLE = {Double ramification cycles and integrable hierarchies},
   JOURNAL = {Comm. Math. Phys.},
  FJOURNAL = {Communications in Mathematical Physics},
    VOLUME = {336},
      YEAR = {2015},
    NUMBER = {3},
     PAGES = {1085--1107}}

@article{ssv:g0,
    AUTHOR = {Shadrin, Sergey and Shapiro, Michael Z. and Vainshtein, Alek},
     TITLE = {Chamber behavior of double {H}urwitz numbers in genus 0},
   JOURNAL = {Adv. Math.},
  FJOURNAL = {Advances in Mathematics},
    VOLUME = {217},
      YEAR = {2008},
    NUMBER = {1},
     PAGES = {79--96}}

@article{bssz15,
    AUTHOR = {Buryak, Alexandr and Shadrin, Sergey and Spitz, Loek and Zvonkine, Dimitri},
     TITLE = {Integrals of {$\psi$}-classes over double ramification cycles},
   JOURNAL = {Amer. J. Math.},
  FJOURNAL = {American Journal of Mathematics},
    VOLUME = {137},
      YEAR = {2015},
    NUMBER = {3},
     PAGES = {699--737}}

@article{Moyal_1949,
    AUTHOR = {Moyal, Jos\'e Enrique},
     TITLE = {Quantum mechanics as a statistical theory},
   JOURNAL = {Proc. Cambridge Philos. Soc.},
  FJOURNAL = {Proceedings of the Cambridge Philosophical Society},
    VOLUME = {45},
      YEAR = {1949},
     PAGES = {99--124}}

@article{cc:c3z3,
    AUTHOR = {Cadman, Charles and Cavalieri, Renzo},
     TITLE = {Gerby localization, {$Z_3$}-{H}odge integrals and the {GW} theory of {$[\mathbb{C}^3/\mathbb{Z}_3]$}},
   JOURNAL = {Amer. J. Math.},
  FJOURNAL = {American Journal of Mathematics},
    VOLUME = {131},
      YEAR = {2009},
    NUMBER = {4},
     PAGES = {1009--1046},
      ISSN = {0002-9327},
     CODEN = {AJMAAN},
   MRCLASS = {14N35},
  MRNUMBER = {2543921 (2010e:14051)},
MRREVIEWER = {Hsian-Hua Tseng},
       DOI = {10.1353/ajm.0.0063},
       URL = {http://dx.doi.org/10.1353/ajm.0.0063}}

@preamble{
   "\def\cprime{$'$} "}

@article{acv03,
    AUTHOR = {Abramovich, Dan and Corti, Alessio and Vistoli, Angelo},
     TITLE = {Twisted bundles and admissible covers},
      NOTE = {Special issue in honor of Steven L. Kleiman},
   JOURNAL = {Comm. Algebra},
  FJOURNAL = {Communications in Algebra},
    VOLUME = {31},
      YEAR = {2003},
    NUMBER = {8},
     PAGES = {3547--3618}}

@article{Jarvis_Kaufmann_Kimura,
    AUTHOR = {Jarvis, Tyler J. and Kaufmann, Ralph and Kimura, Takashi},
     TITLE = {Pointed admissible {$G$}-covers and {$G$}-equivariant
     cohomological field theories},
   JOURNAL = {Compos. Math.},
    VOLUME = {141},
      YEAR = {2005},
    NUMBER = {4},
     PAGES = {926--978}}

@unpublished{bgp:crc,
	Author = {Jim Bryan and Tom Graber and Rahul Pandharipande},
	Note = {Preprint:math.AG/0510335},
	Title = {The Orbifold Quantum Cohomology of $\mathbb{C}^2/\mathbb{Z}_3$ and {Hurwitz-Hodge} Integrals},
	Url = {Preprint:math.AG/0510335},
	Year = {2005}}

@article{elsv:hnaiomsoc,
	Author = {Torsten Ekedahl and Sergei Lando and Michael Shapiro and Alek Vainshtein},
	Journal = {Invent. Math.},
	Pages = {297--327},
	Title = {{Hurwitz} Numbers and intersections on moduli spaces of curves},
	Volume = {146},
	Year = {2001}}

@article{elsv:ohnahi,
	Author = {Torsten Ekedahl and Sergei Lando and Michael Shapiro and Alek Vainshtein},
	Journal = {C.R. Acad.Sci.Paris Ser.I Math.},
	Pages = {1175--1180},
	Title = {On {Hurwitz} Numbers and {Hodge} Integrals},
	Volume = {328},
	Year = {1999}}

@article{fp:hiagwt,
	Author = {Carel Faber and Rahul Pandharipande},
	Journal = {Invent. Math.},
	Pages = {173--199},
	Title = {Hodge integrals and {G}romov-{W}itten theory},
	Volume = {139(1)},
	Year = {2000}}

@article{gjv:ttgodhn,
    AUTHOR = {Goulden, Ian P. and Jackson, David M. and Vakil, Ravi},
     TITLE = {Towards the geometry of double {H}urwitz numbers},
   JOURNAL = {Adv. Math.},
  FJOURNAL = {Advances in Mathematics},
    VOLUME = {198},
      YEAR = {2005},
    NUMBER = {1},
     PAGES = {43--92}}

@article{gv05,
    AUTHOR = {Graber, T. and Vakil, R.},
     TITLE = {Relative virtual localization and vanishing of tautological classes on moduli spaces of curves},
   JOURNAL = {Duke Math. J.},
  FJOURNAL = {Duke Mathematical Journal},
    VOLUME = {130},
      YEAR = {2005},
    NUMBER = {1},
     PAGES = {1--37}}

@article{hm:kd,
	Author = {Joe Harris and David Mumford},
	Journal = {Invent. Math.},
	Pages = {23--88},
	Title = {On the {Kodaira} dimension of the moduli space of curves},
	Volume = {67},
	Year = {1982}}

@unpublished{k:pc,
	Author = {Joachim Kock},
	Note = {Notes. http://mat.uab.es/$\sim $kock/GW/notes/psi-notes.pdf},
	Title = {Notes on psi classes},
	Year = {2001}}

@article{li02,
    AUTHOR = {Li, J.},
     TITLE = {A degeneration formula of {GW}-invariants},
   JOURNAL = {J. Differential Geom.},
  FJOURNAL = {Journal of Differential Geometry},
    VOLUME = {60},
      YEAR = {2002},
    NUMBER = {2},
     PAGES = {199--293}}

@book{mirror,
	Address = {Providence, RI},
	Author = {Hori, Kentaro and Katz, Sheldon and Klemm, Albrecht and Pandharipande, Rahul and Thomas, Richard and Vafa, Cumrun and Vakil, Ravi and Zaslow, Eric},
	Isbn = {0-8218-2955-6},
	Mrclass = {14J32 (14N35 32Q25 81T30 81T60)},
	Mrnumber = {MR2003030 (2004g:14042)},
	Mrreviewer = {Marcos Mari{\~n}o},
	Note = {With a preface by Vafa},
	Pages = {xx+929},
	Publisher = {American Mathematical Society},
	Series = {Clay Mathematics Monographs},
	Title = {Mirror symmetry},
	Volume = {1},
	Year = {2003}}

@incollection{r:adm,
    AUTHOR = {Cavalieri, Renzo},
     TITLE = {Hodge-type integrals on moduli spaces of admissible covers},
 BOOKTITLE = {The interaction of finite-type and {G}romov-{W}itten
              invariants ({BIRS} 2003)},
    SERIES = {Geom. Topol. Monogr.},
    VOLUME = {8},
     PAGES = {167--194},
 PUBLISHER = {Geom. Topol. Publ., Coventry},
      YEAR = {2006}}

@unpublished{r:tqft,
	Author = {Renzo Cavalieri},
	Note = {Preprint: mathAG/0512225},
	Title = {A {TQFT} for intersection numbers on moduli spaces of admissible covers},
	Year = {2005}}

@unpublished{r:dd2,
	Author = {Renzo Cavalieri},
	Note = {Preprint:mathAG/0608590},
	Title = {Generating Functions for {Hurwitz-Hodge} Integrals},
	Year = {2006}}

@article{p:kappa,
    AUTHOR = {Pandharipande, Rahul},
     TITLE = {The {$\varkappa$} ring of the moduli of curves of compact type},
   JOURNAL = {Acta Math.},
  FJOURNAL = {Acta Mathematica},
    VOLUME = {208},
      YEAR = {2012},
    NUMBER = {2},
     PAGES = {335--388}}

@article {cjm:thn,
    AUTHOR = {Cavalieri, Renzo and Johnson, Paul and Markwig, Hannah},
     TITLE = {Tropical {H}urwitz numbers},
   JOURNAL = {J. Algebraic Combin.},
  FJOURNAL = {Journal of Algebraic Combinatorics. An International Journal},
    VOLUME = {32},
      YEAR = {2010},
    NUMBER = {2},
     PAGES = {241--265},
      ISSN = {0925-9899},
     CODEN = {JAOME7},
   MRCLASS = {14N10 (14T05)},
  MRNUMBER = {2661417 (2011m:14089)},
MRREVIEWER = {Ethan G. Cotterill},
       DOI = {10.1007/s10801-009-0213-0},
       URL = {http://dx.doi.org/10.1007/s10801-009-0213-0},
}

@article {cjm:wcfdhn,
    AUTHOR = {Cavalieri, Renzo and Johnson, Paul and Markwig, Hannah},
     TITLE = {Wall crossings for double {H}urwitz numbers},
   JOURNAL = {Adv. Math.},
  FJOURNAL = {Advances in Mathematics},
    VOLUME = {228},
      YEAR = {2011},
    NUMBER = {4},
     PAGES = {1894--1937},
      ISSN = {0001-8708},
     CODEN = {ADMTA4},
   MRCLASS = {14N10 (05C78 05C90 14H30 14T05)},
  MRNUMBER = {2836109},
       DOI = {10.1016/j.aim.2011.06.021},
       URL = {http://dx.doi.org/10.1016/j.aim.2011.06.021},
}

@article {jpt:elsv,
    AUTHOR = {Johnson, P. and Pandharipande, R. and Tseng, H.-H.},
     TITLE = {Abelian {H}urwitz-{H}odge integrals},
   JOURNAL = {Michigan Math. J.},
  FJOURNAL = {Michigan Mathematical Journal},
    VOLUME = {60},
      YEAR = {2011},
    NUMBER = {1},
     PAGES = {171--198},
      ISSN = {0026-2285},
   MRCLASS = {14D23 (14H10 14H30)},
  MRNUMBER = {2785870 (2012c:14025)},
MRREVIEWER = {Arvid Siqveland},
       DOI = {10.1307/mmj/1301586310},
       URL = {http://dx.doi.org/10.1307/mmj/1301586310},
}

@article{br16,
    AUTHOR = {Buryak, Alexandr and Rossi, Paolo},
     TITLE = {Recursion relations for double ramification hierarchies},
   JOURNAL = {Comm. Math. Phys.},
  FJOURNAL = {Communications in Mathematical Physics},
    VOLUME = {342},
      YEAR = {2016},
    NUMBER = {2},
     PAGES = {533--568}}

@article{BR2016q,
    AUTHOR = {Buryak, Alexandr and Rossi, Paolo},
     TITLE = {Double ramification cycles and quantum integrable systems},
   JOURNAL = {Lett. Math. Phys.},
  FJOURNAL = {Letters in Mathematical Physics},
    VOLUME = {106},
      YEAR = {2016},
    NUMBER = {3},
     PAGES = {289--317}}

@article{OkounkovPandharipande2006Virasoro,
  author  = {Okounkov, Andrei and Pandharipande, Rahul},
  title   = {Virasoro constraints for target curves},
  journal = {Inventiones Mathematicae},
  volume  = {163},
  number  = {1},
  pages   = {47--108},
  year    = {2006}}

@incollection{Looijenga1993,
     author = {Looijenga, Eduard},
     title = {Intersection theory on {Deligne-Mumford} compactifications (after {Witten} and {Kontsevich})},
     booktitle = {S\'eminaire Bourbaki},
     series = {Ast\'erisque},
     pages = {187--212},
     note = {Expos\'e 768},
     year = {1993},
     publisher = {Soci\'et\'e math\'ematique de France}}

@incollection{Getzler1998,
  author    = {Getzler, Ezra},
  title     = {Topological recursion relations in genus 2},
  booktitle = {Integrable Systems and Algebraic Geometry ({Kobe}/{Kyoto}, 1997)},
  pages     = {73--106},
  publisher = {World Scientific Publishing},
  address   = {River Edge, NJ},
  year      = {1998}}

@article{spe24,
    AUTHOR = {Spelier, Pim},
     TITLE = {Polynomiality of the double ramification cycle},
   JOURNAL = {Algebr. Geom.},
  FJOURNAL = {Algebraic Geometry},
    VOLUME = {13},
      YEAR = {2026}}

@misc{pix23,
    AUTHOR = {Aaron Pixton},
    TITLE = {{DR} cycle polynomiality and related results},
    YEAR = {2023},
    NOTE = {Available on author's webpage}}

@incollection{Intro2SFT,
    AUTHOR = {Eliashberg, Yakov and Givental, Alexander and Hofer, Helmut},
     TITLE = {Introduction to symplectic field theory},
 BOOKTITLE = {GAFA 2000 (Tel Aviv, 1999)},
 PUBLISHER = {Geom. Funct. Anal.},
  FJOURNAL = {Geometric and Functional Analysis},
      YEAR = {2000},
     PAGES = {560--673}}

@misc{dz01,
      title={Normal forms of hierarchies of integrable {PDE}s, {F}robenius manifolds and {G}romov - {W}itten invariants}, 
      author={Boris Dubrovin and Youjin Zhang},
      note = {arXiv:math/0108160},
      year={2001}}

@misc{bls24,
    author = {Blot, X. and Lewa\'nski, D. and Shadrin, S.},
    title = {On the strong {DR}/{DZ} conjecture}, 
    year = {2024},
    note = {arXiv:2405.12334. To appear in Geom. Topol.}}

@incollection{wit90,
    AUTHOR = {Witten, Edward},
     TITLE = {Two-dimensional gravity and intersection theory on moduli space},
 BOOKTITLE = {Surveys in differential geometry ({C}ambridge, {MA}, 1990)},
     PAGES = {243--310},
 PUBLISHER = {Lehigh Univ.},
   ADDRESS = {Bethlehem, PA},
      YEAR = {1991}}

@article{kon92,
    AUTHOR = {Kontsevich, Maxim},
     TITLE = {Intersection theory on the moduli space of curves and the matrix {A}iry function},
   JOURNAL = {Comm. Math. Phys.},
  FJOURNAL = {Communications in Mathematical Physics},
    VOLUME = {147},
      YEAR = {1992},
    NUMBER = {1},
     PAGES = {1--23}}

@article{br_bihamiltonian,
    AUTHOR = {Buryak, Alexandr and Rossi, Paolo},
     TITLE = {Bihamiltonian structure of the {DR} hierarchy in the semisimple case},
   JOURNAL = {Comm. Math. Phys.},
  FJOURNAL = {Communications in Mathematical Physics},
    VOLUME = {406},
      YEAR = {2025},
    NUMBER = {9},
     PAGES = {Paper No. 205, 20}}

@article{magri_recursion,
    AUTHOR = {Magri, Franco},
     TITLE = {A simple model of the integrable {H}amiltonian equation},
   JOURNAL = {J. Math. Phys.},
  FJOURNAL = {Journal of Mathematical Physics},
    VOLUME = {19},
      YEAR = {1978},
    NUMBER = {5},
     PAGES = {1156--1162}}

@book{riordan,
    AUTHOR = {Riordan, John},
     TITLE = {Combinatorial identities},
      NOTE = {Reprint of the 1968 original},
 PUBLISHER = {Robert E. Krieger Publishing Co.},
   ADDRESS = {Huntington, NY},
      YEAR = {1979},
     PAGES = {xii+256}}

@article{op_trilogy1,
    AUTHOR = {Okounkov, Andrei and Pandharipande, Rahul},
     TITLE = {Gromov-{W}itten theory, {H}urwitz theory, and completed cycles},
   JOURNAL = {Ann. of Math. (2)},
  FJOURNAL = {Annals of Mathematics. Second Series},
    VOLUME = {163},
      YEAR = {2006},
    NUMBER = {2},
     PAGES = {517--560}}

@article{hs21,
    AUTHOR = {Holmes, David and Schmitt, Johannes},
     TITLE = {Infinitesimal structure of the pluricanonical double ramification locus},
   JOURNAL = {Compos. Math.},
  FJOURNAL = {Compositio Mathematica},
    VOLUME = {157},
      YEAR = {2021},
    NUMBER = {10},
     PAGES = {2280--2337}}

@article{BDGR18,
    AUTHOR = {Buryak, A. and Dubrovin, B. and Gu\'er\'e, J. and Rossi, P.},
     TITLE = {Tau-structure for the double ramification hierarchies},
   JOURNAL = {Comm. Math. Phys.},
  FJOURNAL = {Communications in Mathematical Physics},
    VOLUME = {363},
      YEAR = {2018},
    NUMBER = {1},
     PAGES = {191--260}}

@article{brs21,
    AUTHOR = {Buryak, Alexandr and Rossi, Paolo and Shadrin, Sergey},
     TITLE = {Towards a bihamiltonian structure for the double ramification hierarchy},
   JOURNAL = {Lett. Math. Phys.},
  FJOURNAL = {Letters in Mathematical Physics},
    VOLUME = {111},
      YEAR = {2021},
    NUMBER = {1},
     PAGES = {Paper No. 13, 40}}

@article{asKdV00,
    AUTHOR = {Avramidi, Ivan and Schimming, Rainer},
     TITLE = {A new explicit expression for the {K}orteweg-de {V}ries hierarchy},
   JOURNAL = {Math. Nachr.},
  FJOURNAL = {Mathematische Nachrichten},
    VOLUME = {219},
      YEAR = {2000},
     PAGES = {45--64}}

@article{br2025,
    AUTHOR = {Blot, X. and Rossi, P.},
     TITLE = {Meromorphic Differentials, Twisted DR Cycles and Quantum Integrable Hierarchies},
   JOURNAL = {Comm. Math. Phys.},
  FJOURNAL = {Communications in Mathematical Physics},
    VOLUME = {406},
      YEAR = {2025},
    NUMBER = {281},
     PAGES = {1--29}}
\end{document}